\documentclass[12pt]{article}

\usepackage{amsmath,amsfonts}
\usepackage{caption}
\usepackage{accents}
\usepackage{multirow}

\newtheorem{thm}{Theorem}
\newtheorem{ex}{Example}

\newtheorem{lem}[thm]{Lemma}

\newtheorem{assm}{Assumption}
\newtheorem{rem}{Remark}

\def\qed{\hfill\vrule height6pt width5pt depth2pt}

\newcommand\Tstrut{\rule{0pt}{3ex}}            
\newcommand\Bstrut{\rule[-1.5ex]{0pt}{0pt}} 

\newcommand{\ubar}[1]{\underaccent{\bar}{#1}}
\newcommand{\utilde}[1]{\underaccent{\tilde}{#1}}

\newcommand{\be}{\begin{eqnarray*}}
\newcommand{\ben}{\begin{eqnarray}}
\newcommand{\ee}{\end{eqnarray*}}
\newcommand{\een}{\end{eqnarray}}
\newcommand{\al}{\begin{align*}}
\newcommand{\eal}{\end{align*}}
\newcommand{\aln}{\begin{align}}
\newcommand{\ealn}{\end{align}}

\newcommand{\bbe}{\mathbb{E}}
\newcommand{\bbp}{\mathbb{P}}

\newcommand{\bbr}{\mathbb{R}}

\newcommand{\acal}{\mathcal{A}}

\newcommand{\fcal}{\mathcal{F}}

\newcommand{\scal}{\mathcal{S}}
\newcommand{\scalout}{\mathcal{S}_{\mathrm{out}}}
\newcommand{\scalin}{\mathcal{S}_{\mathrm{in}}}

\newcommand{\wcal}{\mathcal{W}}

\newcommand{\pcal}{\mathcal{P}}
\newcommand{\zcal}{\mathcal{Z}}

\def\Om{\Omega}

\let\<\langle
\let\>\rangle
\let\phi\varphi

\begin{document}

\author{Mihail Zervos\thanks{Department of Mathematics,
London School of Economics, Houghton Street, London
WC2A 2AE, UK, Email: mihalis.zervos@gmail.com}, \
Carlos Oliveira\thanks{Department of Mathematics and
CEMAT, Instituto Superior T\'ecnico, Universidade de
Lisboa, Av. Rovisco Pais, 1049-001 Lisboa, Portugal,
Email: carlosmoliveira@tecnico.ulisboa.pt}
\ and Kate Duckworth}

\title{{\sc An Investment Model with Switching Costs \\
and the Option to Abandon} \\ \vspace{3mm}
\normalsize (This is the complete version with proofs
of the paper that is published in \\
{\em Mathematical Methods of Operations Research\/}.)}

\maketitle

\begin{abstract}
\noindent
We develop a complete analysis of a general
entry-exit-scrapping model.
In particular, we consider an investment project that operates
within a random environment and yields a payoff rate that is a
function of a stochastic economic indicator such as the price of
or the demand for the project's output commodity.
We assume that the investment project can operate in two
modes, an ``open'' one and a ``closed'' one.
The transitions from one operating mode to the other one
are costly and immediate, and form a sequence of decisions
made by the project's management. We also assume that the
project can be permanently abandoned at a discretionary
time and at a constant sunk cost.
The objective of the project's management is to maximise
the expected discounted payoff resulting from the project's
management over all switching and abandonment strategies.
We derive the explicit solution to this stochastic control
problem that involves impulse control as well as discretionary
stopping.
It turns out that this has a rather rich structure and the
optimal strategy can take eight qualitatively different forms,
depending on the problemÕs data.

\mbox{}

\end{abstract}

\section{Introduction}

Optimal sequential switching is an area of stochastic
control that emerged from financial economics in the
context of real options (see Dixit and Pindyck~\cite{DP94}
and Trigeorgis~\cite{T96}).
Its numerous applications include the optimal scheduling
of production in a real asset such as a power plant that
can operate in distinct modes, say ``open'' and ``closed'',
as well as the optimal timing of sequentially investing
and disinvesting, e.g., in a given stock.
The references
Bayraktar and Egami \cite{BE10},
Brekke and {\O}ksendal \cite{BO94},
Carmona and Ludkovski \cite{CL08},
Djehiche, Hamad\`{e}ne and Popier \cite{DHP10},
Duckworth and Zervos \cite{DZ01},
El~Asri \cite{EA10},
El~Asri and Hamad\`{e}ne~\cite{EAH09},
Elie and Kharroubi \cite{EK14},
Gassiat, Kharroubi and Pham \cite{GKP12},
Guo and Tomecek \cite{GuoTom08},
Hamad\`{e}ne and Jeanblanc \cite{HamJean07},
Hamad\`{e}ne and Zhang \cite{HamZhang10},
Johnson and Zervos~\cite{JZ10},
Korn, Melnyk and Seifried~\cite{KMS17},
Lumley and Zervos \cite{LZ01},
Ly~Vath and Pham \cite{LyPham07},
Martyr~\cite{Ma16},
Pham \cite{Pham04},
Pham, Ly~Vath and Zhou \cite{PLZ},
Ren\'{e}, Campi, Langren\'{e} and Pham \cite{RCLP},
Song, Yin and Zhang \cite{SYZ09},
Tang and Yong \cite{TY93},
Tsekrekos and Yannacopoulos \cite{TY},
Zhang and Zhang~\cite{ZZ08}, and
Zhang \cite{Z15}
provide an alphabetically ordered list of important
contributions in the area.

In this paper, we derive the complete solution to a
problem of optimal sequential switching that incorporates
an additional permanent abandonment option.
The model that we study goes back to Brennan and
Schwartz~\cite{BrS85} who considered a firm's
decisions to operate, mothball or abandon a mine
producing a natural resource.
A special case of the model is extensively analysed
in Dixit and Pindyck~\cite[Section 7.2]{DP94} using
heuristic arguments and numerical examples in the context
of several real options applications.

To fix ideas, we consider an investment project that
operates within a random environment and yields a
payoff rate that is a function of a stochastic economic
indicator such as the price of or the demand for the
project's output commodity.
We model this economic indicator by the geometric Brownian
motion given by
\ben
dX_t = bX_t \, dt + \sqrt{2} \sigma X_t \, dW_t , \quad
X_0 = x > 0 , \label{SDE}
\een
where $b$ and $\sigma \neq 0$ are given constants
and $W$ is a standard Brownian motion.
We assume that the investment project can operate in
two modes, an ``open'' one and a ``closed'' one.
The transitions from one operating mode to the other one
are immediate and form a sequence of decisions made by the
project's management.
We use a process $Z$ with values in $\{ 0, 1 \}$ to model
such a sequence of decisions.
In particular, we assume that $Z_t = 1$ (resp., $Z_t = 0$)
if the project is ``open'' (resp., ``closed'') at time $t$. 
We also denote by $z \in \{0, 1 \}$ the project's mode at time
$0$, so that $Z_0 = z$.
The stopping times at which the jumps of $Z$ occur are the
intervention times at which the project's operating mode
is changed.
We assume that the project can be permanently
abandoned at a stopping time $\tau$, which is an additional
decision variable.
With each admissible strategy $(Z,\tau)$, we associate
the performance criterion
\begin{align}
J_{z,x} (Z, \tau) = \bbe \biggl[ & \int_0^\tau e^{-rs} h (X_s)
Z_s \, ds \nonumber \\
& - \sum _{j = 1}^\infty e^{-r T_j^1} K_1 {\bf 1}
_{\{ T_j^1 \leq \tau \}} - \sum _{j = 1}^\infty e^{-r T_j^0}
K_0 {\bf 1} _{\{ T_j^0 \leq \tau \}} - e^{-r \tau} K \biggr] ,
\label{J}
\end{align}
where $(T_j^1)$ (resp., $(T_j^0)$) is the sequence
of times at which $Z$ jumps from 0 to 1 (resp., from 1
to 0).
Here, $h : \mbox{} ]0, \infty[ \mbox{} \rightarrow \bbr$ models
the running payoff resulting from the investment project while
this is in its ``open'' operating mode.\footnote{Using a trivial
re-parametrisation, we can allow for the project to yield a
constant payoff rate while it is in its ``closed'' mode
(see Remark~\ref{rem:C}).}
The constants $K_1 > 0$ and $K_0 > 0$ are the costs
resulting from ``switching'' the project from its ``closed''
mode to its ``open'' one and vice versa, whereas
$K \in \bbr$ is the cost resulting from the decision to
permanently abandon it.
Note that we allow for $K$ to be negative, which
corresponds to a situation where capital can be
recovered at abandonment.\footnote{For the same reason,
it would make sense in some economic applications to
allow for at least $K_0$ to be negative, as long as
$K_1 + K_0 > 0$.
However, such a relaxation would add most significant
complexity and would result in a substantially longer
paper.}
Also, on the event $\{ T_j^\ell = \tau \}$, $\ell = 1, 0$,
a cost of $K_\ell + K$ is incurred at time $T_j^\ell$,
which corresponds to the possibility that the project's
operating mode can be switched just before the project
is permanently abandoned.\footnote{Although this
setting is convenient for the problem's formulation,
switching followed by immediate abandonment is
never optimal due to the strict positivity of $K_\ell$,
$\ell = 1, 0$.}
The objective is to maximise the performance criterion
$J_{z,x}$ over the set $\Pi_z$ of all admissible strategies
$(Z, \tau)$.
Accordingly, we define the value function $v$ by
\ben
v(z,x) = \sup _{(Z,\tau) \in \Pi_z} J_{z,x} (Z,\tau) ,
\quad \text{for } (z,x) \in \{ 0,1 \} \times \mbox{}
]0, \infty[ . \label{v}
\een

The related special case that arises if $X = W$,
$h(x) = x$ and $K>0$ was solved by Zervos~\cite{Z03}.
Although the analysis of this related problem has
shed some light on the qualitative nature of the optimal
strategy, its impact on the real options theory has been
limited by the rather unrealistic assumptions that the
underlying economic indicator is a standard Brownian
motion rather than a geometric Brownian motion and
that the running payoff function $h$ is linear.
The existence of an optimal strategy in a more general
context with finite time horizon was established by
Djehiche and Hamad\`{e}ne \cite{DH09} using systems
of Snell envelopes and viscosity solutions.
Despite its fundamental mathematical importance, this
result is of rather limited practical use because it does not
provide a qualitative characterisation of the optimal
strategy or a genuinely practical way of implementing it.

We derive the complete solution to the problem that we
study in an explicit form by solving its
Hamilton-Jacobi-Bellman (HJB) equation that takes
the form of a pair of coupled quasi-variational
inequalities.
In particular, we identify the five regions that partition
the state space $\{ 0,1 \} \times \mbox{} ]0, \infty[$
and characterise the optimal strategy, namely, the
``production'' region, the ``waiting'' region, the
``switch in'' region, the ``switch out'' region and the
``abandonment'' region.
It turns out that the qualitative nature of the problem's
solution is surprisingly rich and can take eight different
forms, depending on the problem data.
We illustrate the results derived using the choice
\ben
h(x) = c + x^\vartheta , \quad x > 0 , \label{h-ex}
\een
for some constants $c \in \bbr$, $\vartheta \in \mbox{}
]0,n[$\,\footnote{The
inequality $\vartheta < n$, where $n$ is defined by
(\ref{n-defn}), is essential for the value function to
be finite.},
and some related numerical calculations
(see Examples~\ref{EXf}--\ref{EXl}).

The value that may be added by waiting before
implementing a certain investment decision is a
central feature of the real options theory.
In some of the cases that arise in our analysis,
{\em value may be added by waiting before
choosing one of two investment actions of a
qualitatively different nature, one partially
reversible and one totally irreversible\/}.
To the best of our knowledge, such a possibility has
not been appreciated in the real options literature.
For instance, in Case~II.3 in Section~\ref{SS-II}
(see also Figure~6), the part of the ``production''
region identified by the set $\{ 1 \} \times \mbox{}
]\delta, \gamma[$ separates the ``abandonment''
region from the ``switch out'' region.
In this case, if the initial condition of the state process
is in this part of the state space, then it is optimal
to take no action before committing to {\em either\/}
enter a perpetual cycle of operating the investment
project by optimally switching it between its two
modes {\em or\/} permanently abandoning the
project, depending on whether the economic
indicator $X$ first rises to the level $\gamma$
or first drops to the level $\delta$.
Furthermore, the investment project has infinite
lifetime if the initial condition of the state process
is in $\{ 1 \} \times [\gamma, \infty[ \mbox{} \cup
\{ 0 \} \times \mbox{} ]0, \infty[$ and finite lifetime
with strictly positive probability otherwise.
The situation becomes more dramatic in
Case~III.2 in Section~\ref{SS-III} (see also Figure~8).
In this case, the part of the ``production''
region identified by the set $\{ 1 \} \times \mbox{}
]\delta, \gamma[$ separates the ``abandonment''
region from the ``switch out'' region, while the
whole ``waiting'' region $\{ 0 \} \times \mbox{}
]\zeta, \alpha[$ separates the ``abandonment''
region from the ``switch in'' region.
If the initial condition of the state process
is in this part of the ``production'' region
(resp., in the ``waiting'' region), then it is
optimal to take no action before committing
to {\em either\/} switch the investment project to
its ``closed'' mode {\em or\/} permanently
abandon it (resp., {\em either\/} switch the
investment project to its ``open'' mode {\em or\/}
permanently abandon it).
Contrary to the previous case, the investment
project's lifetime is always finite with strictly
positive probability, and with probability~1 if
$\mu - \sigma^2 \leq 0$.

The paper is organised as follows.
We formulate the stochastic optimisation problem that
we solve in Section~\ref{pr-form}.
In Section~\ref{HJB-eqns}, we consider the problem's
HJB equation, we discuss how it characterises the
five regions that determine the optimal strategy and
we recall some related implications of the assumptions
we make.
We present the explicit solution to the stochastic control
problem in Section~\ref{solution}.
Here, we organise the eight cases that arise in three
groups based on the analytical affinity of the different
cases.
To simplify the exposition of our main results,
we collect most proofs in two appendixes.

\section{Problem formulation} \label{pr-form}

We build the model that we study on a filtered probability space
$(\Om, \fcal, (\fcal_t), \bbp)$ satisfying the usual conditions and
supporting a standard one-dimensional $(\fcal_t)$-Brownian
motion $W$.
We denote by $\zcal$ the family of all  $(\fcal_t)$-adapted finite 
variation c\`{a}gl\`{a}d processes $Z$ with values in $\{0,1\}$,
and by $\scal$ the set of all $(\fcal_t)$-stopping times.

As we have discussed in the introduction, we consider an
investment project that operates within a random environment
and yields a payoff rate that is a function of a stochastic
economic indicator that is modelled by the geometric
Brownian motion given by (\ref{SDE}).
We assume that the investment project can operate in
two modes, an ``open'' one and a ``closed'' one.
We use a process $Z \in \zcal$ to model such a sequence of
decisions: $Z_t = 1$ (resp., $Z_t = 0$) if the project is ``open''
(resp., ``closed'') at time $t$. 
We also denote by $z \in \{0, 1 \}$ the project's mode at time
$0$, so that $Z_0 = z$.
The stopping times at which the jumps of $Z$ occur are the
intervention times at which the project's operating mode
is changed.
If we define recursively
\begin{gather}
T_1^1 = \inf \left\{ t \geq 0 \mid \ \Delta Z_t = 1 \right\} ,
\quad
T_1^0 = \inf \left\{ t \geq 0 \mid \ \Delta Z_t = -1 \right\} ,
\nonumber \\
T_{j+1}^1 = \inf \left\{ t > T_j^1 \mid \ \Delta Z_t = 1 \right\}
\quad \text{and} \quad
T_{j+1}^0 = \inf \left\{ t > T_j^0 \mid \ \Delta Z_t = -1 \right\} ,
\quad \text{for } j \geq 1 , \nonumber
\end{gather}
where $\Delta Z_t = Z_{t+} - Z_t$ and we adopt the usual
convention that $\inf \emptyset = \infty$, then $T_j^1$
(resp., $T_j^0$) are the $(\fcal_t)$-stopping  times at
which the project is switched from ``closed'' to ``open''
(resp., from ``open'' to ``closed'').
We also assume that the project can be permanently
abandoned at an $(\fcal_t)$-stopping time $\tau$.
We define the set of all admissible strategies to be
\be
\Pi_z = \bigl\{ (Z,\tau) \mid \ Z \in \zcal, \ Z_0 = z,
\text{ and } \tau \in \scal \bigr\} .
\ee
With each admissible strategy $(Z,\tau) \in \Pi_z$, we
associate the performance criterion given by (\ref{J}).
The objective is to maximise the performance criterion
$J_{z,x}$ over $\Pi_z$.
Accordingly, we define the value function $v$ by
(\ref{v}).

For the resulting optimisation problem to be well-posed
in the sense that there are no integrability problems and
there are no admissible strategies with payoff equal to
$\infty$, we make the following assumption.

\begin{assm} \label{Assumption} {\rm
The running payoff function $h : \mbox{} ]0, \infty[ \mbox{}
\rightarrow \bbr$ is right-continuous and increasing,
$\lim _{x \rightarrow \infty} h(x) = \infty$, and
\ben
\bbe \left[ \int _0^\infty e^{-rt} \bigl| h(X_t) \bigr| \, dt \right]
< \infty \label{DOM}
\een
for every initial condition $x > 0$.
Furthermore, $K_1 , K_0 > 0$ and $K \in \bbr$.
} \mbox{}\hfill\qed \end{assm}

\begin{rem} \label{rem:C} {\rm
To simplify the exposition, we have assumed that the investment
project yields zero payoff while it is in its ``closed''
mode.
In view of the calculation 
\begin{align*}
J_{z,x} (Z,\tau) = \bbe \biggl[ & \int_0^\tau e^{-rs} \bigl[ \bar{h}
(X_s) Z_s - C (1-Z_s) \bigr] \, ds \\
& - K_1 \sum _{j = 1}^\infty e^{-r T_j^1} {\bf 1}
_{\{ T_j^1 \leq \tau \}} - K_0 \sum _{j = 1}^\infty e^{-r T_j^0}
{\bf 1} _{\{ T_j^0 \leq \tau \}} - e^{-r \tau} \bar{K} \biggr]
+ \frac{C}{r} ,
\end{align*} 
where $C$ is a constant, $\bar{h} = h - C$ and $\bar{K}
= K + \frac{C}{r}$, we can see that allowing for a constant
payoff rate while the project is in its ``closed''
mode can be accommodated trivially in the model that we
study.
} \mbox{} \hfil\qed \end{rem}

\section{The Hamilton-Jacobi-Bellman (HJB) equation}
\label{HJB-eqns}

In view of standard stochastic control theory that has been
developed and used in references we have discussed in the
introduction, we expect that the value function of the problem
we study is given by
\ben
v(1, \cdot) = w_1 \quad \text{and} \quad v(0, \cdot) = w_0 ,
\label{eqn:v=w}
\een
where the functions $w_1, w_0 : \mbox{} ]0, \infty[ \mbox{}
\rightarrow \bbr$ satisfy the coupled quasi-variational
inequalities
\begin{align} 
\max \Bigl\{ \sigma^2 x^2 w_1'' (x) + bx w_1' (x) - rw_1 (x)
+ h(x) , \ w_0 (x) - w_1 (x) - K_0 , \  -w_1 (x) - K \Bigr\} & = 0 ,
\label{HJBa} \\
\max \Bigl\{ \sigma^2 x^2 w_0'' (x) + bx w_0' (x) - rw_0 (x) ,
\  w_1 (x) - w_0 (x) - K_1 , \ -w_0 (x) - K \Bigr\} & = 0 ,
\label{HJBb}
\end{align}
as well as appropriate growth conditions (see
Zervos~\cite[Theorem~1]{Z03} for a general verification
theorem).
In view of the heuristics explaining the structure of this
HJB equation, the state space $\{ 0,1 \} \times \mbox{}
]0, \infty[$ splits into five pairwise disjoint regions\footnote{In
the description of the five possible regions, we characterise
subsets of $]0,\infty[$ as open or closed relative to the
topology on $]0,\infty[$ that is the trace of the usual
topology on $\bbr$, for instance, $]0,a] = \mbox{} ]0,\infty[
\mbox{} \setminus \mbox{} ]a, \infty[$ and $[a, \infty[ \mbox{}
= \mbox{} ]0,\infty[ \mbox{} \setminus \mbox{} ]0, a[$ are
closed sets.}:
\smallskip

\noindent {\bf (i)}
{\em The ``production'' region $\{ 1 \} \times \pcal$,
where $\pcal$ is an open subset of $]0,\infty[$\/}.
Whenever the project is in its ``open'' mode and the
process $X$ takes values in $\pcal$, it is optimal to
keep the project in its ``open'' mode, which is associated
with production.
In particular, $\pcal$ is the set in which the function $w_1$
satisfies the ODE
\ben
\sigma^2 x^2 w''(x) + bx w'(x) - rw(x) + h(x) = 0 . \label{ODEb}
\een

\noindent {\bf (ii)}
{\em The ``waiting'' region $\{ 0 \} \times \wcal$,
where $\wcal$ is an open subset of $]0,\infty[$\/}.
If the project is in its ``closed'' mode and the process $X$
takes values in $\wcal$, then it is optimal to take no action,
namely, keep the project on standby.
The set $\wcal$ is characterised by the requirement that
$w_0$ satisfies the ODE
\ben
\sigma^2 x^2 w''(x) + b x w'(x) - rw(x) = 0 . \label{ODEa}
\een

\noindent {\bf (iii)}
{\em The ``switch out'' region $\{ 1 \} \times \scalout$,
where $\scalout$ is a closed subset of $]0,\infty[$\/}.
If the project is in its ``open'' mode, then it is optimal
to switch it to its ``closed'' mode as soon as $X$
takes values in $\scalout$.
The set $\scalout$ is characterised by the identity
\ben
w_1 (x) = w_0 (x) - K_0 \quad \text{for all } x \in \scalout .
\label{Sout-region}
\een

\noindent {\bf (iv)}
{\em The ``switch in'' region $\{ 0 \} \times \scalin$,
where $\scalin$ is a closed subset of $]0,\infty[$\/}.
It is optimal to switch the project from its ``closed'' to its
``open'' mode as soon as $X$ takes values in $\scalin$.
In this case,
\ben
w_0 (x) = w_1 (x) - K_1 \quad \text{for all } x \in \scalin .
\label{Sin-region}
\een

\noindent {\bf (v)}
{\em The ``abandonment'' region $\{ 0 \} \times \acal_0
\cup \{ 1 \} \times \acal_1$, where $\acal_0$, $\acal_1$ are
closed subsets of $]0,\infty[$\/}.
It is optimal to abandon permanently the project as soon
as the state process hits the abandonment region.
Accordingly,
\ben
w_i (x) = -K \quad \text{for all } x \in \acal_i \text{ and }
i=0,1 . \label{ab-region}
\een

\noindent
The tactics associated with these regions exhaust all possible
control actions.
Therefore,
\be
\pcal \cup \scalout \cup \acal_1 = \wcal \cup \scalin
\cup \acal_0 = \mbox{} ]0, \infty[ .
\ee

We will solve the control problem that we study by
identifying these regions and deriving appropriate
explicit solutions to the HJB equation
(\ref{HJBa})--(\ref{HJBb}).
To this end, we will use the following facts.
It is well-known that the general solution to the
Euler's ODE (\ref{ODEa}) is given by
\ben
w(x) = Ax^m + Bx^n , \label{ODEa-sol}
\een
for some constants $A, B \in \bbr$, where the constants
$m < 0 < n$ are defined by
\ben
m, n = \frac{1}{2\sigma^2} \left[ \sigma^2 - b \mp
\sqrt{(b - \sigma^2)^2 + 4\sigma^2r} \right] .
\label{n-defn}
\een
If $h : \mbox{} ]0, \infty[ \mbox{} \rightarrow \bbr$ is a function
satisfying the integrability condition in (\ref{DOM}), then a
particular solution to the ODE (\ref{ODEb}) is the function
$R_h : \mbox{} ]0, \infty[ \mbox{} \rightarrow \bbr$ given by
\begin{align} 
R_h(x) & = \frac{1}{\sigma^2 (n-m)} \left[ x^m \int _0^x
s^{-m-1} h(s) \, ds + x^n \int _x^\infty s^{-n-1} h(s) \, ds
\right] \nonumber \\
& = \bbe \left[ \int _0^\infty e^{-rs} h(X_s) \, ds \right] .
\label{Rh}
\end{align}
A straightforward calculation reveals that
\ben
R_h' (x) = \frac{1}{\sigma^2 (n-m)} \left[ m x^{m-1} \int _0^x
s^{-m-1} h(s) \, ds + n x^{n-1} \int _x^\infty s^{-n-1} h(s)
\, ds \right] . \label{Rh'(x)}
\een
Furthermore, for a choice of $h$ as in Assumption~\ref{Assumption},
\begin{gather}
R_h \text{ is increasing} , \label{Rh1} \\
h(0) := \lim_{x \downarrow 0} h(x) = r \lim_{x \downarrow 0}
R_h (x) \quad \text{and} \quad
\lim_{x \rightarrow \infty} R_h (x) = \infty , \label{Rh2} \\
\lim_{T \rightarrow \infty} e^{-rT} \, \bbe \Bigl[ \bigl| R_h (X_T)
\bigr| \Bigr] = 0 \label{Rh3} \\
\text{and} \quad
\bbe \left[ \int_0^T e^{-2rt} X_t^2 \bigl| R_h' (X_t) \bigr| ^2
\, dt \right] < \infty \quad \text{for all } T>0 . \label{Rh4} 
\end{gather} 
All of these claims regarding the function $R_h$ as well as
several more general results can be found in Knudsen,
Meister and Zervos~\cite[Section~4]{KMZ98}, and
Johnson and Zervos~\cite{JZ07}.

\begin{ex} \label{EXf} {\rm
If $h$ is the function given by (\ref{h-ex}), then
Assumption~\ref{Assumption} holds true if and only if
$\vartheta \in \mbox{} ]0, n[$,
in which case,
\be
R_h (x) = - \frac{x^\vartheta}{\sigma^2 \vartheta^2
+ (b - \sigma^2) \vartheta - r} + \frac{c}{r} .
\ee
We will illustrate our results numerically for the choices
\be
b = 0 , \quad \sigma = 1 , \quad r = 2 , \quad
\vartheta = 1 \quad \text{and} \quad K_1 = K_0 =
\frac{1}{2} , 
\ee
which are associated with
\be
m = -1, \quad n = 2 \quad \text{and} \quad
R_h (x) = \frac{1}{2} x + \frac{c}{2} . 
\ee
} \end {ex}

\section{The solution to the control problem} \label{solution}

We now derive the solution to the stochastic control problem
formulated in Section~\ref{pr-form} by identifying the
sets $\pcal$, $\wcal$, $\scalout$, $\scalin$, $\acal_1$,
$\acal_0$ we have discussed in the previous section and
deriving appropriate solutions to the HJB equation
(\ref{HJBa})--(\ref{HJBb}) using (\ref{ODEb})--(\ref{ab-region}).
To this end, we first note that, if the investment project is in
its ``open'' mode at time $0$ and is never switched to its
``closed'' mode or abandoned, then it will yield a total expected
discounted payoff equal to $R_h(x)$ (see (\ref{Rh})).
On the other hand, if the project is ``closed'' at time $0$ and
is never switched to its ``open'' operating mode or abandoned,
then it will yield 0 total expected discounted payoff. 
Since $R_h$ is increasing and $\lim _{x \rightarrow \infty}
R_h (x) = \infty$ (see (\ref{Rh1}) and (\ref{Rh2})), it should
be optimal to operate the project in its ``open'' mode whenever
the process $X$ takes sufficiently high values.
It follows that there exists $M>0$ such that
\be
]M, \infty[ \mbox{} \subseteq \pcal \quad \text{and} \quad
]M, \infty[ \mbox{} \subseteq \scalin .
\ee
If $\acal_1 \neq \emptyset$ (resp, $\acal_0 \neq \emptyset$),
then $\acal_1 = \mbox{} ]0, \delta]$ (resp., $\acal_0 = \mbox{}
]0, \zeta]$) for some $\delta > 0$ (resp., $\zeta > 0$) because
$R_h$ is increasing.
Furthermore, in view of the smoothness of a solution to
the HJB equation (\ref{HJBa})--(\ref{HJBb}) that is required
to identify it with the control problem's value function and
the analysis in the previous section, we expect that
the ``abandonment'' region does not have any common
boundary points with either the ``switch in'' region or the
``switch out'' region.

In light of these observations, we will show that the
production and the waiting regions $\pcal$ and $\wcal$
have the general forms
\ben
\pcal = \mbox{} ]\delta, \gamma[ \mbox{} \cup \mbox{}
]\beta, \infty[ \quad \text{and} \quad
\wcal = \mbox{} ]\zeta, \alpha[ , \label{PW-gen}
\een
for some $0 \leq \delta \leq \gamma \leq \beta < \infty$
and $0 \leq \zeta \leq \alpha < \infty$ (see Figures~1-8),
where we adopt the usual convention that, e.g., $]0,0[
\mbox{} = \emptyset$.
In view of the solutions to the ODEs (\ref{ODEb}),
(\ref{ODEa}) given in the previous section, the solution
to the HJB equation (\ref{HJBa})--(\ref{HJBb}) should
be such that
\ben
w_1 (x) = \left. \begin{cases} R_h (x) , & \text{for all } x \in
\mbox{} ]0, \infty[ , \text{ if } \delta = \gamma = \beta = 0 \\
Ax^m + R_h (x) , & \text{for all } x \in \mbox{} ]\beta, \infty[ ,
\text{ if } \gamma < \beta \text{ or } 0 < \delta = \gamma =
\beta \\ \Gamma_1 x^m + \Gamma_2 x^n + R_h (x) , &
\text{for all } x \in \mbox{} ]\delta, \gamma[ , \text{ if } 0 <
\delta < \gamma < \beta \end{cases} \right\} \label{w1-gen}
\een
and
\ben
w_0 (x) = \left. \begin{cases} Bx^n , & \text{for all }
x \in \mbox{} ]0, \alpha[ , \text{ if } \zeta = 0 < \alpha
\\ \Delta_1 x^m + \Delta_2 x^n , & \text{for all } x \in \mbox{}
]\zeta, \alpha[ , \text{ if } 0 < \zeta < \alpha \end{cases}
\right\} , \label{w0-gen}
\een
for some constants $A$, $\Gamma_1$, $\Gamma_2$,
$B$, $\Delta_1$ and $\Delta_2$ because these are the
only choices that are consistent with the requirements of
the verification theorem that we will use to identify the
solution to (\ref{HJBa})--(\ref{HJBb}) with the control
problem's value function.

To determine free-boundary points such as $\delta$,
$\gamma$, $\beta$, $\zeta$, $\alpha$ appearing in
(\ref{PW-gen}) and constants such as $A$, $\Gamma_1$,
$\Gamma_2$, $B$, $\Delta_1$, $\Delta_2$ appearing in
(\ref{w1-gen})--(\ref{w0-gen}), we will use the $C^1$
continuity that we expect the functions $w_1$, $w_0$
to have.
In particular, we will require that $w_1$, $w_0$ should
be $C^1$ at every boundary point separating any two of
the five regions.
Using the expressions (\ref{Rh}), (\ref{Rh'(x)}) and the
identity $\sigma^2 mn = - r$, we will then derive
appropriate systems of equations for the unknown
parameters.
We will only provide the results of these calculations
because they are straightforward to replicate.

We have organised the presentation of the possible cases
arising by splitting them in three groups.
Group I includes the cases in which it is not optimal to
switch or abandon the project if this is in its ``open''
mode.
Group II contains all cases where it may be optimal to
switch or abandon the project if this is in its ``open'' mode
but abandonment is not optimal if the project is in its
``closed'' mode.
Finally, Group~III includes all remaining cases. 
To make the presentation easier to follow, we develop
the proofs in Appendix~II.

\subsection{Group I: taking action is not optimal whenever the
project is in its ``open'' operating mode ($\pmb{\pcal = \mbox{}
]0,\infty[}$)} \label{SS-I}

All cases in this group are such that $\pcal = \mbox{}
]0,\infty[$ and are associated with a solution to the HJB
equation (\ref{HJBa})--(\ref{HJBb}) such that
\ben
w_1 (x) = R_h (x) \quad \text{for all } x > 0 .
\label{w1,I}
\een

\noindent
{\bf Case I.1 (Figure~1)}
In this case, it is optimal to immediately switch the investment
project to its ``open'' mode if it is originally ``closed''.
Accordingly,
\be
\pcal = \scalin = \mbox{} ]0,\infty[ \quad \text{and} \quad
\wcal = \scalout = \acal_0 = \acal_1 = \emptyset ,
\ee
and the functions $w_1$ and $w_0$ given by (\ref{w1,I})
and
\ben
w_0 (x) = R_h(x) - K_1, \quad \text{for } x>0 , \label{w0,I1}
\een
should satisfy the HJB equation (\ref{HJBa})--(\ref{HJBb}).

\begin{picture}(160,65)

\put(35,30){\begin{picture}(80,40)
\put(-20,0){\vector(1,0){120}}
\put(-20,15){\vector(1,0){120}}
\put(-30,-0.6){\scriptsize $z=0$}
\put(-30,14.4){\scriptsize $z=1$}
\put(103,-0.6){\scriptsize $x$}
\put(103,14.4){\scriptsize $x$}

\put(90,-5){\vector(0,1){10}}
\put(80,-5){\vector(0,1){10}}
\put(70,-5){\vector(0,1){10}}
\put(60,-5){\vector(0,1){10}}
\put(50,-5){\vector(0,1){10}}
\put(40,-5){\vector(0,1){10}}
\put(30,-5){\vector(0,1){10}}
\put(20,-5){\vector(0,1){10}}
\put(10,-5){\vector(0,1){10}}
\put(0,-5){\vector(0,1){10}}
\put(-10,-5){\vector(0,1){10}}
\put(25,-10){\scriptsize $w_0(x)=R_h(x)-K_1$}
\put(33,-15){\scriptsize (switch in)}
\put(30,20.5){\scriptsize (production)}
\put(28.5,25.5){\scriptsize $w_1(x)=R_h(x)$}
\end{picture}}

\put(35,8){\small{{\bf Figure 1.} Illustration of the regions
determining}}
\put(37,3){\small{the optimal strategy in the context of
Case I.1}}
\end{picture}

\begin{lem} \label{lem I.1}
The increasing functions $w_1$, $w_0$ defined by (\ref{w1,I}),
(\ref{w0,I1}) satisfy the HJB equation (\ref{HJBa})--(\ref{HJBb})
if and only if
\be
\max\{ rK_1, \, rK_1 - rK \} \leq h(0) .
\ee 
\end{lem}

\begin{ex} {\rm
If $h$ is the function given by (\ref{h-ex}) and the
problem data is as in Example~\ref{EXf}, then this
case characterises the optimal strategy if and only if
$K \in \bbr$ and $\max\{ 1, \, 1 - rK \} \leq c$.
} \end {ex}

\noindent
{\bf Case I.2 (Figure~2)}
In this case, it is optimal to switch the investment project to its
``open'' mode if it is originally ``closed'' as long as the process
$X$ takes sufficiently high values.
In particular, there exists a boundary point $\alpha > 0$ such that,
if the project starts in its ``closed'' mode, then it is optimal to wait
for all long as $X$ takes values strictly less than $\alpha$ and
switch the project to its ``open'' mode as soon as $X$ takes a
value exceeding $\alpha$.
Accordingly,
\be
\pcal = \mbox{} ]0, \infty[ , \quad \wcal = \mbox{} ]0, \alpha[ ,
\quad \scalin = [\alpha, \infty[ \quad \text{and} \quad
\scalout = \acal_0 = \acal_1 = \emptyset .
\ee
In view of (\ref{Sin-region}) and (\ref{w1-gen})--(\ref{w0-gen}),
the functions $w_1$ and $w_0$ given by (\ref{w1,I}) and
\ben
w_0 (x) = \left. \begin{cases}
Bx^n , & \text{if } x < \alpha \\
R_h(x) - K_1, & \text{if } x \geq \alpha 
\end{cases} \right\} \label{w0,I2}
\een
should satisfy the HJB equation (\ref{HJBa})--(\ref{HJBb}).

\begin{picture}(160,65)

\put(35,30){\begin{picture}(80,40)
\put(-20,0){\vector(1,0){120}}
\put(-20,15){\vector(1,0){120}}
\put(-30,-0.6){\scriptsize $z=0$}
\put(-30,14.4){\scriptsize $z=1$}
\put(103,-0.6){\scriptsize $x$}
\put(103,14.4){\scriptsize $x$}
\put(20,0){\circle*{1.3}}
\put(19,-4){$\alpha$}
\put(50,-5){\vector(0,1){10}}
\put(60,-5){\vector(0,1){10}}
\put(70,-5){\vector(0,1){10}}
\put(80,-5){\vector(0,1){10}}
\put(90,-5){\vector(0,1){10}}
\put(40,-5){\vector(0,1){10}}
\put(30,-5){\vector(0,1){10}}
\put(51,-10){\scriptsize $w_0(x)=R_h(x)-K_1$}
\put(58,-15){\scriptsize (switch in)}
\put(-10,-10){\scriptsize $w_0(x)=Bx^n$}
\put(-8,-15){\scriptsize (waiting)}
\put(30,20.5){\scriptsize (production)}
\put(28.5,25.5){\scriptsize $w_1(x)=R_h(x)$}
\end{picture}}

\put(35,8){\small{{\bf Figure 2.} Illustration of the regions
determining}}
\put(37,3){\small{the optimal strategy in the context of
Case I.2}}
\end{picture}

\noindent
The requirement that $w_0$ should be $C^1$ at $\alpha$
yields the expressions
\begin{gather}
B = \frac{1}{\sigma^2 (n-m)} \int _\alpha^\infty s^{-n-1}
\bigl[ h(s) - rK_1 \bigr] \, ds \label{parameter B,I2} \\
\text{and} \quad
\int_0^\alpha s^{-m-1} \bigl[ h(s) - rK_1 \bigr] \, ds = 0 .
\label{eq I.2}
\end{gather}

\begin{lem}\label{lem I.2}
Equation (\ref{eq I.2}) has a unique solution ${\alpha}>0$
and the functions $w_1$, $w_0$  defined by (\ref{w1,I}),
(\ref{w0,I2}), for $B>0$ given by (\ref{parameter B,I2}), are
increasing and satisfy the HJB equation (\ref{HJBa})--(\ref{HJBb})
if and only if
\be
0 \leq K \quad \text{and} \quad \max \{ -rK_0, \, -rK \} \leq h(0)
< rK_1 .
\ee 
\end{lem}

\begin{ex} {\rm
If $h$ is the function given by (\ref{h-ex}), then
(\ref{parameter B,I2}) and (\ref{eq I.2}) are equivalent to
\be
\alpha = \left( - \frac{(\vartheta - m) (rK_1 - c)}{m} \right)
^{-1/m} \quad \text{and} \quad
B = \frac{\alpha^{-n}}{\sigma^2 (n-m)} \left(
\frac{\alpha^\vartheta}{n - \vartheta} - \frac{rK_1 - c}{n}
\right) .
\ee
If the problem data is as in Example~\ref{EXf}, then this
case characterises the optimal strategy if and only if
$0 \leq K$ and $\max \{ {-1} , \, - rK \} \leq c < 1$.
In particular, if $c = \frac{1}{2}$, then
\be
\alpha = 1 \quad \text{and} \quad B = \frac{1}{4} .
\ee

} \end{ex}

\noindent
{\bf Case I.3 (Figure~3)}
This case differs from the previous one by the fact that 
abandoning the investment project if it is in its ``closed''
mode and the process $X$ takes values below a given
threshold level $\zeta$ becomes optimal.
Accordingly,
\be
\pcal = \mbox{} ]0, \infty[ , \quad \acal_0 = \mbox{} ]0, \zeta]
, \quad \wcal = \mbox{} ]\zeta, \alpha[ , \quad \scalin =
[\alpha, \infty[ \quad \text{and} \quad \scalout = \acal_1 =
\emptyset ,
\ee
and, in view of (\ref{Sin-region})--(\ref{ab-region}) and
(\ref{w1-gen})--(\ref{w0-gen}), the required solution to the
HJB equation (\ref{HJBa})--(\ref{HJBb}) should be given
by the function $w_1$ defined by (\ref{w1,I}) and the
function $w_0$ defined by
\ben
w_0 (x) = \left. \begin{cases}
-K , & \text{if } x \leq \zeta \\
\Delta_1 x^m + \Delta_2 x^n, & \text{if } \zeta
< x < \alpha \\ R_h (x) - K_1, & \text{if } x \geq
\alpha  \end{cases} \right\} . \label{w0,I3}
\een

\begin{picture}(160,65)

\put(35,30){\begin{picture}(80,40)
\put(-20,0){\vector(1,0){120}}
\put(-20,15){\vector(1,0){120}}
\put(-30,-0.6){\scriptsize $z=0$}
\put(-30,14.4){\scriptsize $z=1$}
\put(103,-0.6){\scriptsize $x$}
\put(103,14.4){\scriptsize $x$}
\put(50,0){\circle*{1.3}}
\put(10,0){\circle*{1.3}}
\put(49,-4){$\alpha$}
\put(9,-4){$\zeta$}
\put(60,-5){\vector(0,1){10}}
\put(70,-5){\vector(0,1){10}}
\put(80,-5){\vector(0,1){10}}
\put(90,-5){\vector(0,1){10}}
\put(63,-10){\scriptsize $w_0(x)=R_h(x)-K_1$}
\put(71,-15){\scriptsize (switch in)}
\put(18,-10){\scriptsize $w_0(x)=\Delta_1 x^m+\Delta_2 x^n$}
\put(27,-15){\scriptsize (waiting)}
\put(-15,-10){\scriptsize $w_0(x)=-K$}
\put(-16.5,-15){\scriptsize (abandonment)}
\put(30,20.5){\scriptsize (production)}
\put(28.5,25.5){\scriptsize $w_1(x)=R_h(x)$}
\linethickness{0.5mm}
\put(-20,0){\qbezier(0,0)(15,0)(30,0)}
\end{picture}}

\put(35,8){\small{{\bf Figure 3.} Illustration of the regions
determining}}
\put(37,3){\small{the optimal strategy in the context of
Case I.3}}
\end{picture}

\noindent
To determine the free-boundary points $\zeta$, $\alpha$
and the parameters $\Delta_1$, $\Delta_2$, we require
that $w_0$ should be $C^1$, which yields the expressions
\begin{gather}
f_1 (\zeta, \alpha) := m \int _0^\alpha s^{-m-1} \bigl[
h(s) - rK_1 \bigr] \, ds - rK \zeta^{-m} = 0 , \label{eq I.3.1} \\
f_2 (\zeta, \alpha) := n \int _\alpha^\infty s^{-n-1} \bigl[
h(s) - rK_1 \bigr] \, ds + rK \zeta^{-n} = 0 , \label{eq I.3.2} \\
\Delta_1 = \frac{rK \zeta^{-m}}{\sigma^2 m(n-m)}
\quad \text{and} \quad
\Delta_2 = -\frac{rK \zeta^{-n}}{\sigma^2 n(n-m)}
. \label{P I.3.5}
\end{gather}

\begin{lem}\label{lem I.3}
The system of equations (\ref{eq I.3.1})--(\ref{eq I.3.2})
has a unique solution $(\zeta, \alpha)$ such that
$0 < \zeta < \alpha$ and the functions $w_1$, $w_0$
defined by (\ref{w1,I}), (\ref{w0,I3}), for $\Delta_1 > 0$,
$\Delta_2 > 0$ given by (\ref{P I.3.5}), are increasing and
satisfy the HJB equation  (\ref{HJBa})--(\ref{HJBb})
if and only if
\be
K<0 \quad \text{and} \quad -rK \leq h(0) < rK_1 - rK .
\ee
\end{lem}

\begin{ex} {\rm
If $h$ is the function given by (\ref{h-ex}), then the
system of equations (\ref{eq I.3.1})--(\ref{eq I.3.2})
takes the form
\begin{align}
(rK_1 - c) \alpha^{-m} + \frac{m}{\vartheta - m}
\alpha ^{\vartheta - m} - rK \zeta ^{-m} & = 0 ,
\nonumber \\
(rK_1 - c) \alpha^{-n} - \frac{n}{n - \vartheta}
\alpha ^{- (n - \vartheta)} - rK \zeta ^{-n} & = 0 .
\nonumber
\end{align}
If the problem data is as in Example~\ref{EXf}, then this
case characterises the optimal strategy if and only if
$K<0$ and $- rK \leq c < 1 -rK$.
In particular, if $K = - \frac{1}{2}$ and $c = 1$, then
\be
\zeta = 2^{- \frac{1}{3}}, \quad \alpha = 2^{\frac{1}{3}}
\quad \text{and} \quad \Delta_1 = \Delta_2 =
2^{- \frac{1}{3}} \times 3^{-1} .
\ee
} \end{ex}

\subsection{Group II: taking action may be optimal if the project
is in its ``open'' mode but abandonment is not optimal whenever
the project is in its ``closed'' operating mode ($\pmb{\pcal \neq
\mbox{} ]0, \infty[}$ and $\pmb{\acal_0 = \emptyset}$)}
\label{SS-II}

We now consider cases that complement the ones in the
previous group and are characterised by the non-optimality
of abandonment whenever the project is in its ``closed''
mode.
In all of these cases, $\wcal = \mbox{} ]0, \alpha[$ and
$\scalin = [\alpha, \infty[$.
Otherwise, the cases are differentiated by the arrangement
of the optimal tactics whenever the project is in its ``open''
mode.
\medskip

\noindent 
{\bf Case II.1 (Figure~4)}
In this case, sequential switching of the investment project from
``open'' to ``closed'' and vice versa is optimal, and abandonment
is not part of the optimal strategy.
Whenever the project is in its ``open'' (resp., ``closed'') mode,
it is optimal to stay there for as long as the process $X$
takes values above (resp., below) a given threshold $\beta$
(resp., $\alpha$) and switch to its ``closed'' (resp., ``open'')
mode as soon as $X$ takes values below (resp., above) the
threshold $\beta$ (resp., $\alpha$).
Of course, for such a strategy to be well-defined, we must
have $\beta < \alpha$.
Accordingly,
\be
\scalout = \mbox{} ]0, \beta], \quad \pcal = \mbox{} ]\beta, \infty[
, \quad \wcal = \mbox{} ]0, \alpha[ , \quad \scalin = [\alpha, \infty[
\quad \text{and} \quad \acal_0 = \acal_1 = \emptyset .
\ee
In view of (\ref{Sout-region})--(\ref{Sin-region}) and
(\ref{w1-gen})--(\ref{w0-gen}), we can see that the required
solution to the HJB equation (\ref{HJBa})--(\ref{HJBb})
should be given by the functions defined by
\begin{align}
w_1 (x) & = \left. \begin{cases} Bx^n - K_0 , & \text{if }
x \leq \beta  \\ Ax^m + R_h (x) , & \text{if } x > \beta 
\end{cases} \right\} \label{w0,w1,II1-2} \\
\text{and} \quad w_0 (x) & = \left. \begin{cases}
Bx^n , & \text{if } x < \alpha  \\
Ax^m + R_h (x) - K_1 , & \text{if } x \geq \alpha 
\end{cases} \right\} . \label{w0,w1,II1-1}
\end{align}

\begin{picture}(160,65)

\put(35,30){\begin{picture}(80,40)
\put(-20,0){\vector(1,0){120}}
\put(-20,15){\vector(1,0){120}}
\put(-30,-0.6){\scriptsize $z=0$}
\put(-30,14.4){\scriptsize $z=1$}
\put(103,-0.6){\scriptsize $x$}
\put(103,14.4){\scriptsize $x$}

\put(50,0){\circle*{1.3}}
\put(49,-4){$\alpha$}
\put(60,-5){\vector(0,1){10}}
\put(70,-5){\vector(0,1){10}}
\put(80,-5){\vector(0,1){10}}
\put(90,-5){\vector(0,1){10}}

\put(20,15){\circle*{1.3}}
\put(19,17){$\beta$}
\put(-10,20){\vector(0,-1){10}}
\put(0,20){\vector(0,-1){10}}
\put(10,20){\vector(0,-1){10}}
\put(57,-10){\scriptsize $w_0(x)=Ax^m+R_h(x)-K_1$}
\put(68,-15){\scriptsize (switch in)}
\put(-10,-8){\scriptsize $w_0(x)=Bx^n$}
\put(-7,-13){\scriptsize (waiting)}
\put(57,20){\scriptsize (production)}
\put(50,25){\scriptsize $w_1(x)=Ax^m+R_h(x)$}
\put(-7.5,22.5){\scriptsize (switch out)}
\put(-12,27.5){\scriptsize $w_1(x)=Bx^n-K_0$}
\end{picture}}

\put(35,8){\small{{\bf Figure 4.} Illustration of the regions
determining}}
\put(37,3){\small{the optimal strategy in the context of
Case II.1}}
\end{picture}

\noindent
To determine the free-boundary points $\beta$, $\alpha$ and
the parameters $A$, $B$, we once again require that the
functions $w_1$, $w_0$ should be $C^1$, which yields
the expressions
\begin{align}
A & = -\frac{1}{\sigma^2 (n-m)} \int _0^\beta s^{-m-1} \bigl[
h(s) + rK_0 \bigr] \, ds , \label{P II.1.1} \\   
B & = \frac{1}{\sigma^2 (n-m)} \int_\alpha^\infty s^{-n-1}
\bigl[ h(s) - rK_1 \bigr] \, ds , \label{P II.1.2}
\end{align}
and the system of equations 
\begin{align}
m \int _\beta^\alpha s^{-m-1} h(s) \, ds + rK_0 \beta^{-m}
+ rK_1 \alpha^{-m} & = 0 , \label{eq II.1.1} \\
n \int _\beta^\alpha s^{-n-1} h(s) \, ds + rK_0 \beta^{-n}
+ rK_1\alpha^{-n} & = 0 . \label{eq II.1.2}
\end{align}

\begin{lem}\label{lem II.1}
The system of equations (\ref{eq II.1.1})--(\ref{eq II.1.2})
has a unique solution $(\beta, \alpha)$ such that
$0 < \beta < \alpha$ and the functions $w_1$, $w_0$
defined by (\ref{w0,w1,II1-2}), (\ref{w0,w1,II1-1}), for $A>0$,
$B>0$ given by (\ref{P II.1.1}), (\ref{P II.1.2}), are increasing
and satisfy the HJB equation (\ref{HJBa})--(\ref{HJBb}) if and
only if
\be
K_0 \leq K \quad \text{and} \quad h(0) < -rK_0 .
\ee 
\end{lem}

\begin{ex} {\rm
If $h$ is the function given by (\ref{h-ex}), then the
system of equations (\ref{eq II.1.1})--(\ref{eq II.1.2})
takes the form
\begin{align}
(rK_1 - c) \alpha^{-m} + (rK_0 + c) \beta ^{-m}
+ \frac{m}{\vartheta - m} \left( \alpha ^{\vartheta - m}
- \beta ^{\vartheta - m} \right) & = 0 ,
\nonumber \\
(rK_1 - c) \alpha^{-n} + (rK_0 + c) \beta ^{-n}
- \frac{n}{n - \vartheta} \left( \alpha ^{- (n - \vartheta)}
- \beta ^{- (n - \vartheta)} \right) & = 0 , \nonumber
\end{align}
while
\be
A = \frac{\beta ^{-m}}{\sigma^2 (n-m)} \left(
\frac{rK_0 + c}{m} - \frac{\beta ^\vartheta}{\vartheta - m}
\right) \quad \text{and} \quad
B = \frac{\alpha ^{-n}}{\sigma^2 (n-m)} \left(
\frac{\alpha ^\vartheta}{n - \vartheta} - \frac{rK_1 - c}{n} 
\right) . \nonumber
\ee
If the problem data is as in Example~\ref{EXf}, then this
case characterises the optimal strategy if and only if
$\frac{1}{2} \leq K$ and $c < {-1}$.
In particular, if $c = -2$, then
\be
\beta \simeq 0.537 , \quad \alpha \simeq 5.866 ,
\quad A \simeq 0.131 \quad \text{and} \quad B \simeq
0.042 .
\ee
} \end{ex}

\noindent
{\bf Case II.2 (Figure~5)}
Abandoning the project if this is in its ``open'' mode and
the state process $X$ takes values below a given
threshold $\delta_\dagger$ instead of switching it to its
``closed'' mode is the difference between this case and
the previous one.\footnote{We use the notation
$\delta_\dagger$ rather than the simpler $\delta$
because this point will appear in assumptions that
we will make in later cases.}
Accordingly,
\be
\acal_1 = \mbox{} ]0, \delta_\dagger], \quad \pcal =
\mbox{} ]\delta_\dagger, \infty[ , \quad \wcal = \mbox{}
]0, \alpha[ , \quad \scalin = [\alpha, \infty[ \quad \text{and}
\quad \scalout = \acal_0 = \emptyset ,
\ee
and the functions defined by
\begin{align}
w_1 (x) & = \left. \begin{cases} -K , &  \text{if } x \leq
\delta_\dagger  \\ Ax^m + R_h(x) , & \text{if } x >
\delta_\dagger \end{cases} \right\} \label{w0,w1,II2-2} \\
\text{and} \quad w_0 (x) & = \left. \begin{cases}
Bx^n , & \text{if } x < \alpha  \\
Ax^m + R_h(x) - K_1, & \text{if } x \geq \alpha 
\end{cases} \right\} \label{w0,w1,II2-1}
\end{align}
should provide a solution to the HJB equation
(\ref{HJBa})--(\ref{HJBb}).

\begin{picture}(160,65)

\put(35,30){\begin{picture}(80,40)
\put(-20,0){\vector(1,0){120}}
\put(-20,15){\vector(1,0){120}}
\put(-30,-0.6){\scriptsize $z=0$}
\put(-30,14.4){\scriptsize $z=1$}
\put(103,-0.6){\scriptsize $x$}
\put(103,14.4){\scriptsize $x$}

\put(50,0){\circle*{1.3}}
\put(49,-4){$\alpha$}
\put(60,-5){\vector(0,1){10}}
\put(70,-5){\vector(0,1){10}}
\put(80,-5){\vector(0,1){10}}
\put(90,-5){\vector(0,1){10}}

\put(20,15){\circle*{1.3}}
\put(19,17){$\delta_\dagger$}
\put(57,-10){\scriptsize $w_0(x)=Ax^m+R_h(x)-K_1$}
\put(69,-15){\scriptsize (switch in)}
\put(-10,-8){\scriptsize $w_0(x)=Bx^n$}
\put(-7,-13){\scriptsize (waiting)}
\put(52,20){\scriptsize (production)}
\put(45,25){\scriptsize $w_1(x)=Ax^m+R_h(x)$}
\put(-11.5,20){\scriptsize (abandonment)}
\put(-10,25){\scriptsize $w_1(x)=-K$}
\linethickness{0.5mm}
\put(-20,15){\qbezier(0,0)(20,0)(40,0)}
\end{picture}}

\put(35,8){\small{{\bf Figure 5.} Illustration of the regions
determining}}
\put(37,3){\small{the optimal strategy in the context of
Case II.2}}
\end{picture}

\noindent
Requiring that $w_1$, $w_0$ should be $C^1$, we obtain the
expressions
\begin{align}
A & = - \frac{1}{\sigma^2 (n-m)} \int _0^{\delta_\dagger}
s^{-m-1} \bigl[ h(s) + rK \bigr] \, ds , \label{P II.2.5M} \\
B & = \frac{1}{\sigma^2 (n-m)} \int _\alpha^\infty s^{-n-1}
\bigl[ h(s) - rK_1 \bigr] \, ds , \label{P II.2.8M}
\end{align}
and the system of equations
\begin{align}
\int _{\delta_\dagger}^\infty s^{-n-1} \bigl[ h(s) + rK
\bigr] \, ds & = 0 , \label{P II.2.6M} \\
f(\delta,\alpha) := m \int _{\delta_\dagger}^\alpha
s^{-m-1} \bigl[ h(s) - rK_1 \bigr] \, ds + r (K_1 + K)
\delta_\dagger^{-m} & = 0 . \label{P II.2.9M}
\end{align}
The following result involves the point
\ben
K_0^\star = -K_1 - \frac{m \hat{x}^m}{r} \int _{\hat{x}}
^\alpha s^{-m-1} \bigl[ h(s) - rK_1 \bigr] \, ds ,
\label{K0-star-defn}
\een
where $\hat{x}$ solves the equation
\ben
m \hat{x}^m \int _{\hat{x}}^\alpha s^{-m-1} \bigl[
h(s) - rK_1 \bigr] \, ds  - n \hat{x}^n \int _{\hat{x}}^\alpha
s^{-n-1} \bigl[ h(s) - rK_1 \bigr] \, ds = 0 . \label{xhat-eqn}
\een

\begin{lem} \label{lem II.2}
The system of equations (\ref{P II.2.6M})--(\ref{P II.2.9M})
has a unique solution $(\delta_\dagger, \alpha)$ such that
$0 < \delta_\dagger < \alpha $, while equation (\ref{xhat-eqn})
has a unique solution $\hat{x} \in \mbox{} ]\delta_\dagger ,
\alpha[$.
Given these solutions, the functions $w_1$,
$w_0$ defined by (\ref{w0,w1,II2-2}), (\ref{w0,w1,II2-1}), for
$A>0$, $B>0$ given by (\ref{P II.2.5M}), (\ref{P II.2.8M}),
are increasing and satisfy the HJB equation
(\ref{HJBa})--(\ref{HJBb}) if and only if
\be
0 \leq K
\ee
and
\be
\Bigl( K < K_0  \text{ and } -rK_0 \leq h(0) < -rK
\Bigr) \quad \text{or} \quad \Bigl( K_0^\star \leq K_0
\text{ and } h(0) < -rK_0 \Bigr) ,
\ee
where $K_0^\star \in \mbox{} ]K, - r^{-1} h(0)[$, which
depends on all problem data except $K_0$, is defined
by (\ref{K0-star-defn}).
\end{lem}

\begin{ex} \label{ex-ddag} {\rm
If $h$ is the function given by (\ref{h-ex}), then the
system of equations (\ref{P II.2.6M})--(\ref{P II.2.9M})
takes the form
\begin{gather}
\delta_\dagger^\vartheta = - \frac{(n - \vartheta) (c+rK)}{n}
, \nonumber \\
(rK_1 - c) \alpha^{-m} + \frac{m}{\vartheta - m} \alpha
^{\vartheta - m} + (rK + c) \delta_\dagger^{-m}
- \frac{m}{\vartheta - m} \delta_\dagger^{\vartheta - m}
= 0 , \nonumber
\end{gather}
while
\be
A = - \frac{\vartheta (rK + c) \delta_\dagger^{-m}}
{r (\vartheta - m)}
\quad \text{and} \quad
B = \frac{\alpha ^{-n}}{\sigma^2 (n-m)} \left(
\frac{\alpha ^\vartheta}{n - \vartheta} - \frac{rK_1 - c}{n} 
\right) . \nonumber
\ee
The critical point $K_0^\star$ defined by (\ref{K0-star-defn})
admits the expression
\begin{gather}
K_0^\star = -K_1 + \frac{rK_1 - c}{r} \left[ 1 - \left(
\frac{\hat{x}}{\alpha} \right) ^m \right] - \frac{m}
{r (\vartheta - m)}  \hat{x}^m ( \alpha ^{\vartheta - m}
-  \hat{x}^{\vartheta - m} ) \in \mbox{} ]K, - r^{-1} c[ ,
\nonumber
\end{gather}
where $\hat{x}$ is the unique solution to the equation
\begin{gather}
(rK_1 - c) \left[ \left( \frac{\hat{x}}{\alpha} \right) ^m
- \left( \frac{\hat{x}}{\alpha} \right) ^n \right] +
\frac{m}{\vartheta - m} \alpha ^{\vartheta - m}
\hat{x}^m + \frac{n}{n - \vartheta} \alpha ^{-(n-\vartheta)}
\hat{x}^n - \frac{\vartheta (n-m)}{(n - \vartheta) (\vartheta - m)}
\hat{x}^\vartheta = 0 . \nonumber
\end{gather}
If the problem data is as in Example~\ref{EXf}, then this
case characterises the optimal strategy if and only if
either ($0 \leq K < \frac{1}{2}$ and ${-1} \leq c < -rK$)
or ($0 \leq K$, $K_0^\star \leq \frac{1}{2}$ and $c < -1$).
If $K = 0$ and $c = -1$, then
\be
\delta_\dagger = \frac{1}{2} , \quad \alpha =
2 + \frac{\sqrt{13}}{2} , \quad A = \frac{1}{8}
\quad \text{and} \quad B \simeq 0.065 .
\ee
while, if $K = \frac{1}{4}$ and $c = -2$, then
\be
\hat{x} \simeq 0.808 , \quad K_0^\star = 0.276 ,
\quad \delta_\dagger = \frac{3}{4} , \quad \alpha =
3 + \frac{\sqrt{117}}{4} , \quad A = \frac{9}{32}
\quad \text{and} \quad B \simeq 0.043 .
\ee
} \end{ex}

\noindent
{\bf Case II.3 (Figure~6)}
The last case in this group is a hybrid of the previous two.
If the investment project is initially in its ``open'' mode and
the initial value $x$ of the process $X$ is greater than
a threshold $\gamma$ or it is initially in its ``closed'' mode,
then it is optimal to follow the same strategy as in Case~II.1,
which is determined by two thresholds $\beta < \alpha$
such that $\gamma < \beta$.
In this case, the project is sequentially switched from ``open''
to ``closed'' and vice versa, and it is never abandoned.
On the other hand, if the project is initially in its ``open''
mode and the initial value $x$ of $X$ is strictly less than
$\gamma$, then it is optimal to abandon the project as
soon as $X$ falls below another threshold $\delta < \gamma$
before hitting $\gamma$.
Otherwise, it is optimal to switch the project to its ``closed''
mode if $X$ rises to $\gamma$ before hitting $\delta$,
and then maintain the sequential switching strategy defined
by $\beta$ and $\alpha$.
Accordingly,
\begin{gather}
\acal_1 = \mbox{} ]0, \delta] , \quad \pcal = \mbox{} ]\delta,
\gamma[ \mbox{} \cup \mbox{} ]\beta, \infty[ , \quad \scalout
= [\gamma, \beta] , \nonumber \\
\wcal = \mbox{} ]0, \alpha[ , \quad \scalin = [\alpha, \infty[
\quad \text{and} \quad \acal_0 = \emptyset . \nonumber
\end{gather}
In view of (\ref{Sout-region})--(\ref{ab-region}) and
(\ref{w1-gen})--(\ref{w0-gen}), we can see that the required
solution to the HJB equation (\ref{HJBa})--(\ref{HJBb})
should be given by the functions defined by
\begin{align}
w_1 (x) & = \left. \begin{cases}
-K , & \text{if } x \leq \delta \\
\Gamma_1 x^m + \Gamma_2 x^n + R_h(x) , & \text{if }
\delta < x < \gamma \\
Bx^n - K_0 , & \text{if } \gamma \leq x \leq \beta \\
Ax^m + R_h (x) , & \text{if } x > \beta
\end{cases} \right\} \label{w0,w1,II3-2M} \\
\text{and} \quad
w_0 (x) & = \left. \begin{cases}
Bx^n , & \text{if } x < \alpha  \\
Ax^m + R_h(x) - K_1, & \text{if } x \geq \alpha
\end{cases} \right\} . \label{w0,w1,II3-1M}
\end{align}

\begin{picture}(160,65)

\put(35,30){\begin{picture}(80,40)
\put(-20,0){\vector(1,0){120}}
\put(-20,15){\vector(1,0){120}}
\put(-30,-0.6){\scriptsize $z=0$}
\put(-30,14.4){\scriptsize $z=1$}
\put(103,-0.6){\scriptsize $x$}
\put(103,14.4){\scriptsize $x$}

\put(73,0){\circle*{1.3}}
\put(72,-4){$\alpha$}
\put(81,-5){\vector(0,1){10}}
\put(91,-5){\vector(0,1){10}}
\put(0,15){\circle*{1.3}}
\put(-1,17){$\delta$}
\put(32,15){\circle*{1.3}}
\put(30.5,17){$\gamma$}
\put(40,20){\vector(0,-1){10}}
\put(50,20){\vector(0,-1){10}}
\put(57,15){\circle*{1.3}}
\put(56,17){$\beta$}

\put(72,-10){\scriptsize $w_0(x)=w_1(x)-K_1$}
\put(79.5,-15){\scriptsize (switch in)}
\put(10,-7){\scriptsize $w_0(x)=Bx^n$}
\put(13,-12){\scriptsize (waiting)}
\put(-21,22.5){\scriptsize (abandonment)}
\put(-19.5,27.5){\scriptsize $w_1(x)=-K$}
\put(37,22.5){\scriptsize (switch out)}
\put(32,27.5){\scriptsize $w_1(x) = w_0(x) - K_0$}
\put(7.5,17.5){\scriptsize (production)}
\put(1,22.5){\scriptsize $\Gamma_1x^m+\Gamma_2x^n+R_h(x)$}
\put(10,27.5){\scriptsize $w_1(x)=$}
\put(71,22.5){\scriptsize (production)}
\put(65,27.5){\scriptsize $w_1(x)=Ax^m+R_h(x)$}
\linethickness{0.5mm}
\put(-20,15){\qbezier(0,0)(10,0)(20,0)}
\end{picture}}

\put(35,8){\small{{\bf Figure 6.} Illustration of the regions
determining}}
\put(37,3){\small{the optimal strategy in the context of
Case II.3}}
\end{picture}

\noindent
To determine $\Gamma_1$, $\Gamma_2$, $A$, $B$, $\delta$,
$\gamma$, $\beta$ and $ \alpha$ we require that $w_1$, $w_0$
should be $C^1$ at the free-boundary points $\delta$, $\gamma$,
$\beta$ and $\alpha$.
In view of this requirement, we can verify that $\delta$, $\gamma$,
$\beta$ and $ \alpha$ should satisfy the equations (\ref{eq II.1.1}),
(\ref{eq II.1.2}),
\begin{align}
F_1(\delta,\gamma) := m \int _\delta^\gamma s^{-m-1}
\bigl[ h(s) + rK_0 \bigr] \, ds + r(K - K_0) \delta^{-m} & = 0
\label{P II.3.14} \\
\text{and} \quad
F_2(\delta,\gamma) := n \int _\delta^\gamma s^{-n-1}
\bigl[ h(s) + rK_0 \bigr] \, ds + r(K - K_0) \delta^{-n}
\hspace{3mm} & \nonumber \\
\mbox{} + n \int _\beta^\infty s^{-n-1} \bigl[ h(s) + rK_0 \bigr]
\, ds & = 0 , \label{P II.3.15}
\end{align}
while $A$, $B$, $\Gamma_1$ and $\Gamma_2$ should be
given by (\ref{P II.1.1}), (\ref{P II.1.2}),
\begin{align}
\Gamma_1 & = - \frac{1}{\sigma^2 (n-m)} \int _0^\gamma
s^{-m-1} \bigl[ h(s) + rK_0 \bigr] \, ds \label{P II.3.10} \\
\text{and} \quad
\Gamma_2 & = - \frac{1}{\sigma^2 (n-m)} \int _\gamma^\beta
s^{-n-1} \bigl[ h(s) + rK_0 \bigr] \, ds. \label{P II.3.11}
\end{align}

\begin{lem}\label{lem II.3}
The system of equations (\ref{eq II.1.1}), (\ref{eq II.1.2}),
(\ref{P II.3.14}) and (\ref{P II.3.15}) has a unique solution
$(\delta, \gamma, \beta, \alpha)$ such that $0 < \delta <
\gamma < \beta < \alpha$ and the functions $w_1$, $w_0$
defined by (\ref{w0,w1,II3-2M}), (\ref{w0,w1,II3-1M}), for
$A>0$, $B>0$, $\Gamma_1>0$, $\Gamma_2>0$ given by
(\ref{P II.1.1}), (\ref{P II.1.2}), (\ref{P II.3.10}), (\ref{P II.3.11}),
are increasing and satisfy the HJB equation
(\ref{HJBa})--(\ref{HJBb}) if and only if
\be
0 \leq K , \quad h(0) < -rK_0 \quad \text{and}
\quad K < K_0 < K_0^\star ,
\ee
where $K_0^\star \in \mbox{} ]K, - r^{-1} h(0) \bigr[$,
which depends on all problem data except $K_0$,
is as in Lemma~\ref{lem II.2}.
\end{lem}

\noindent
We note that the conditions of this result can all hold
true only if $h(0) < 0$.

\begin{ex} {\rm
If $h$ is the function given by (\ref{h-ex}), then the
system of equations (\ref{P II.3.14})--(\ref{P II.3.15})
takes the form
\begin{align}
(rK_0 + c) (\gamma ^{-m} - \delta ^{-m}) -
\frac{m}{\vartheta - m} (\gamma ^{\vartheta - m}
- \delta ^{\vartheta - m}) + r(K_0-K) \delta ^{-m} & = 0 ,
\nonumber \\
(rK_0 + c) (\delta ^{-n} - \gamma ^{-n} + \beta^{-n})
+ \frac{n}{n - \vartheta} (\delta ^{- (n - \vartheta)} -
\gamma ^{-(n - \vartheta)} + \beta ^{- (n - \vartheta)})
- r(K_0-K) \delta ^{-n} & = 0 , \nonumber
\end{align}
while
\begin{gather}
\Gamma_1 = \frac{\gamma ^{-m}}{\sigma^2 (n-m)} \left(
\frac{rK_0 + c}{m} - \frac{\gamma ^\vartheta}{\vartheta - m}
\right) \nonumber \\
\text{and} \quad
\Gamma_2 = - \frac{1}{\sigma^2 (n-m)} \left[
\frac{rK_0 + c}{n} \left( \gamma^{-n} - \beta ^{-n} \right)
+ \frac{1}{n - \vartheta} \left( \gamma ^{-(n - \vartheta)}
- \beta ^{-(n - \vartheta)} \right) \right] . \nonumber
\end{gather}
If the problem data is as in Example~\ref{EXf}, then this
case characterises the optimal strategy if and only if
$0 \leq K < \frac{1}{2}$, $c < {-1}$ and $ \frac{1}{2}
< K_0^\star$, where $K_0^\star$ is as in
Example~\ref{ex-ddag}.
In particular, if $K = \frac{5}{11}$ and $c = -4$, then
\begin{gather*}
\hat{x} \simeq 1.706 , \quad K_0^\star \simeq 0.524 ,
\quad \delta \simeq 0.279 , \quad \gamma \simeq 1.348
, \quad \beta \simeq 1.740 , \quad \alpha \simeq 9.194 , \\
A \simeq 1.235 , \quad B \simeq 0.026 , \quad \Gamma_1
\simeq 1.045 \quad \text{and} \quad \Gamma_2 \simeq
0.054 .
\end{gather*}
} \end{ex}

\subsection{Group III: the remaining cases} \label{SS-III}

We now consider the remaining cases.
These are characterised by the fact that it may be optimal
to abandon the investment project when this is in its ``closed''
mode.
\medskip

\noindent 
{\bf Case III.1 (Figure~7)}
This case is the modification of Case~II.2 (see Figure~5) that
arises if abandonment when the project is in its ``closed''
mode becomes part of the optimal tactics.
In this case,
\be
\acal_1 = \mbox{} ]0, \delta_\dagger], \quad \pcal = \mbox{}
]\delta_\dagger, \infty[ , \quad \acal_0 = \mbox{} ]0, \zeta] ,
\quad \wcal = \mbox{} ]\zeta, \alpha[ , \quad \scalin =
[\alpha, \infty[ \quad \text{and} \quad \scalout = \emptyset ,
\ee
and the functions defined by
\begin{align}
w_1 (x) & = \left. \begin{cases}
-K , & \text{if } x \leq \delta_\dagger \\
Ax^m + R_h (x) , & \text{if } x \geq \delta_\dagger 
\end{cases} \right\} \label{w0,w1,III1-1M} \\
\text{and} \quad
w_0 (x) & = \left. \begin{cases}
-K , & \text{if } x \leq \zeta \\
\Delta_1 x^m + \Delta_2 x^n, & \text{if } \zeta \leq
x \leq \alpha \\
Ax^m + R_h (x) - K_1 , & \text{if } x  \geq \alpha 
\end{cases} \right\} \label{w0,w1,III1-0M}
\end{align}
should provide a solution to the HJB equation
(\ref{HJBa})--(\ref{HJBb}).

\begin{picture}(160,65)

\put(35,30){\begin{picture}(80,40)\label{figuere II.2}
\put(-20,0){\vector(1,0){120}}
\put(-20,15){\vector(1,0){120}}
\put(-30,-0.6){\scriptsize $z=0$}
\put(-30,14.4){\scriptsize $z=1$}
\put(103,-0.6){\scriptsize $x$}
\put(103,14.4){\scriptsize $x$}

\put(50,0){\circle*{1.3}}
\put(49,-4){$\alpha$}
\put(60,-5){\vector(0,1){10}}
\put(70,-5){\vector(0,1){10}}
\put(80,-5){\vector(0,1){10}}
\put(90,-5){\vector(0,1){10}}
\put(20,15){\circle*{1.3}}
\put(19,17){$\delta_\dagger$}
\put(10,0){\circle*{1.3}}
\put(9,-4){$\zeta$}
\put(10,0){\line(0,1){5}}
\put(10,5){\line(1,0){15}}
\put(25,5){\vector(0,-1){5}}

\put(57,-10){\scriptsize $w_0(x)=Ax^m+R_h(x)-K_1$}
\put(68,-15){\scriptsize (switch in)}
\put(15,-10){\scriptsize $w_0(x) = \Delta_1 x^m+ \Delta_2 x^n$}
\put(24,-15){\scriptsize (waiting)}
\put(-15,-10){\scriptsize $w_1(x)=-K$}
\put(-16.5,-15){\scriptsize (abandonment)} 
\put(51.5,20){\scriptsize (production)}
\put(45,25){\scriptsize $w_1(x)=Ax^m+R_h(x)$}
\put(-11,20){\scriptsize (abandonment)}
\put(-10,25){\scriptsize $w_1(x)=-K$}
\linethickness{0.5mm}
\put(-20,15){\qbezier(0,0)(20,0)(40,0)}
\linethickness{0.5mm}
\put(-20,0){\qbezier(0,0)(15,0)(30,0)}
\end{picture}}

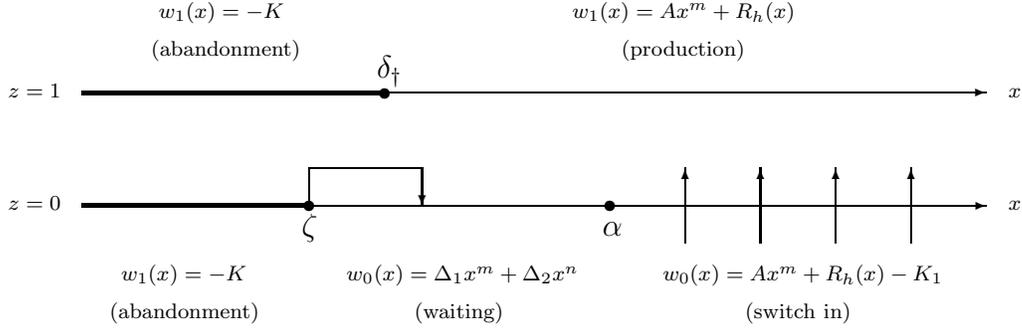
\put(15,8){\small{{\bf Figure 7.} Illustration of the regions
determining the optimal strategy}}
\put(15,3){\small{ in the context of Case III.1
($\zeta$ can be smaller as well as larger than
$\delta_\dagger$)}}
\end{picture}

\noindent
To determine $A$, $\Delta_1$, $\Delta_2$, $\delta_\dagger$,
$\zeta$ and $ \alpha$ we require that $w_1$, $w_0$
should be $C^1$ at the free-boundary points $\delta_\dagger$,
$\zeta$ and $\alpha$.
In view of this requirement, we can verify that $\delta_\dagger$,
$\zeta$ and $\alpha$ should satisfy the system of equations
\begin{align}
G_1 (\delta_\dagger, \zeta, \alpha) := m \int _{\delta_\dagger}
^\alpha s^{-m-1} \bigl[ h(s) - rK_1 \bigr] \, ds + r (K_1 + K)
\delta_\dagger^{-m} - rK \zeta^{-m} & = 0 \label{P III.1.14} \\
\text{and} \quad
G_2 (\delta_\dagger, \zeta, \alpha) := -n \int _{\delta_\dagger}
^\alpha s^{-n-1} \bigl[ h(s) - rK_1\bigr] \, ds - r (K_1 + K)
\delta_\dagger^{-n} + rK \zeta^{-n} & = 0 , \label{P III.1.15}
\end{align}
where $\delta_\dagger$ is given by (\ref{P II.2.6M}),
while, $A$, $\Delta_1$ and $\Delta_2$ should be given by
(\ref{P II.2.5M}),
\begin{align}
\Delta_1 & = A + \frac{1}{\sigma^2 (n-m)} \int _0^\alpha
s^{-m-1} \bigl[ h(s) - rK_1 \bigr] \, ds =
\frac{rK \zeta^{-m}}{\sigma^2 m (n-m)} \label{P III.1.12} \\
\text{and} \quad
\Delta_2 &= \frac{1}{\sigma^2 (n-m)} \int _\alpha^\infty
s^{-n-1} \bigl[ h(s) - rK_1 \bigr] \, ds =
- \frac{rK \zeta^{-n}}{\sigma^2 n (n-m)} . \label{P III.1.13}
\end{align}
The following result involves the equation
\ben
G_2 \bigl( \delta_\dagger, \delta_\dagger, \alpha
(K_1) ; K_1 \bigr) = 0 , \label{lem7-z-d-rel3-m}
\een
for $K_1$, in which we make explicit the dependence
of $\alpha$ and $G_2$ on $K_1$ (note that $\delta_\dagger$
does not depend on $K_1$).
Also, it involves 
the point
\ben
K_0^\dagger = -K_1 - \frac{n \hat{x}^n}{r} \int _{\hat{x}}
^\alpha s^{-n-1} \bigl[ h(s) - rK_1 \bigr] \, ds ,
\label{K0-dag-defn}
\een
where $\hat{x}$ solves the equation
\ben
m \hat{x}^m \int _{\hat{x}}^\alpha s^{-m-1} \bigl[
h(s) - rK_1 \bigr] \, ds  - n \hat{x}^n \int _{\hat{x}}^\alpha
s^{-n-1} \bigl[ h(s) - rK_1 \bigr] \, ds = 0 . \label{xhat-eqn2}
\een

\begin{lem} \label{lem III.1}
The system of equations (\ref{P II.2.6M}), (\ref{P III.1.14})
and (\ref{P III.1.15}) has a unique solution $(\delta_\dagger,
\zeta, \alpha)$ such that $0 < \delta_\dagger \wedge \zeta
\leq \delta_\dagger \vee \zeta < \alpha$.
If $h(\delta_\dagger) < 0$, then there exists a unique
solution $K_1^\dagger > 0$ to (\ref{lem7-z-d-rel3-m})
that depends on all of the problem data except $K_1$, 
$K_0$.
If $h(\delta_\dagger) < 0$ and $K < K_1^\dagger$, then
equation (\ref{xhat-eqn2}) has a unique solution
$\hat{x} \in \mbox{} ]\delta_\dagger, \alpha[$ and
the point $K_0^\dagger > 0$ depends on all of the
problem data except $K_0$.
Furthermore, $\lim _{K_1 \uparrow K_1^\dagger}
K_0^\dagger (K_1) = 0$, and the free-boundary points
$\zeta$ and $\delta_\dagger$, which do not depend on
$K_0$, are such that
\begin{align}
& 0 < \zeta < \delta_\dagger \quad \text{if } h(\delta_\dagger)
< 0 \text{ and } K_1 < K_1^\dagger , \label{lem7-d>z} \\
& 0 < \zeta = \delta_\dagger \quad \text{if } h(\delta_\dagger)
< 0 \text{ and } K_1 = K_1^\dagger \label{lem7-d=z} \\
\text{and} \quad
& 0 < \delta_\dagger < \zeta \quad \text{if } h(\delta_\dagger)
\geq 0 \text{ or } \Bigl( h(\delta_\dagger)
< 0 \text{ and } K_1 > K_1^\dagger \Bigr) .
\label{lem7-d<z}
\end{align}
The functions $w_1$, $w_0$ defined by (\ref{w0,w1,III1-1M}),
(\ref{w0,w1,III1-0M}), for $A>0$, $\Delta_1 > 0$, $\Delta_2 > 0$
given by (\ref{P II.2.5M}), (\ref{P III.1.12}), (\ref{P III.1.13}), are
increasing and satisfy the HJB equation
(\ref{HJBa})--(\ref{HJBb}) if and only if
\be
K < 0, \quad h(0)< -rK
\ee
and
\begin{gather}
\Bigl( -rK_0 \leq h(0) \Bigr) \quad \text{or} \quad
\Bigl( h(0) < -rK_0 \text{ and } h(\delta_\dagger) \geq 0
\Bigr) \nonumber \\
\text{or} \quad \Bigl( h(0) < -rK_0 , \  h(\delta_\dagger)
< 0 \text{ and } K_1 \geq K_1^\dagger \Bigr) \nonumber \\
\text{or} \quad \Bigl( h(0) < -rK_0 , \ h(\delta_\dagger)
< 0 , \ K_1 < K_1^\dagger \text{ and } K_0 \geq
K_0^\dagger \Bigr) . \nonumber
\end{gather}
\end{lem}

\begin{ex} {\rm
If $h$ is the function given by (\ref{h-ex}), then the
system of equations (\ref{P III.1.14})--(\ref{P III.1.15})
takes the form
\begin{align}
(rK_1 - c) (\alpha^{-m} - \delta _\dagger^{-m}) +
\frac{m}{\vartheta - m} (\alpha ^{\vartheta - m} - \delta
_\dagger^{\vartheta - m}) + r(K_1 + K) \delta _\dagger^{-m}
- rK \zeta^{-m} & = 0 , \nonumber \\
(rK_1 - c) (\alpha^{-n} - \delta _\dagger^{-n}) -
\frac{n}{n - \vartheta} (\alpha ^{-(n - \vartheta)} - \delta
_\dagger^{- (n - \vartheta)}) + r(K_1 + K) \delta _\dagger^{-n}
- rK \zeta^{-n} & = 0 , \nonumber
\end{align}
where $\delta_\dagger$ admits the expression given
in Example~\ref{ex-ddag}.
The equation (\ref{lem7-z-d-rel3-m}) that the critical point
$K_1^\dagger$ satisfies if $h (\delta_\dagger) < 0$
takes the form
\be
c \bigl[ \alpha^{-n} (K_1) - \delta_\dagger^{-n} \bigr]
+ \frac{n}{n - \vartheta} \bigl[ \alpha^{-(n-\vartheta)} (K_1)
- \delta_\dagger^{-(n-\vartheta)} \bigr] - K_1 \alpha^{-n}
(K_1) = 0 ,
\ee
while that critical point $K_0^\star$ defined by
(\ref{K0-dag-defn})  if $h (\delta_\dagger) < 0$
admits the expression
\begin{gather}
K_0^\star = -K_1 + \frac{rK_1 - c}{r} \left[ 1 - \left(
\frac{\hat{x}}{\alpha} \right) ^n \right] + \frac{n}
{r (n-\vartheta)}  \hat{x}^n ( \alpha ^{-(n-\vartheta)}
-  \hat{x}^{-(n-\vartheta)} ) ,
\nonumber
\end{gather}
where $\hat{x}$ is the unique solution to the equation
\begin{gather}
(rK_1 - c) \left[ \left( \frac{\hat{x}}{\alpha} \right) ^m
- \left( \frac{\hat{x}}{\alpha} \right) ^n \right] +
\frac{m}{\vartheta - m} \alpha ^{\vartheta - m}
\hat{x}^m + \frac{n}{n - \vartheta} \alpha ^{-(n-\vartheta)}
\hat{x}^n - \frac{\vartheta (n-m)}{(n - \vartheta) (\vartheta - m)}
\hat{x}^\vartheta = 0 . \nonumber
\end{gather}
If the problem data is as in Example~\ref{EXf}, then this
case characterises the optimal strategy if and only if
($K < 0$ and ${-1} \leq c < -rK$) or
($K < 0$ and $rK \leq c < {-1}$) or
($K < 0$, $c < \min \{ {-1}, rK \}$ and $K_1^\dagger \leq
\frac{1}{2}$) or
($K < 0$, $c < \min \{ {-1}, rK \}$, $\frac{1}{2} <
K_1^\dagger$ and $K_0^\dagger \leq \frac{1}{2}$).
If $K = - \frac{1}{2}$ and $c = 0$, then
\be
\delta_\dagger = \frac{1}{2} , \quad \zeta \simeq 1.283
, \quad \alpha \simeq 2.678 , \quad A = \frac{1}{8} , \quad
\Delta_1 \simeq 0.428 \quad \text{and} \quad \Delta_2
\simeq 0.101 ,
\ee
if $K = - 1$ and $c = -\frac{3}{2}$, then
\be
\delta_\dagger = \frac{7}{4} , \quad \zeta \simeq 2.625
, \quad \alpha \simeq 5.250 , \quad A = \frac{49}{32} ,
\quad \Delta_1 \simeq 1.750 \quad \text{and} \quad
\Delta_2 \simeq 0.048 ,
\ee
if $K = - \frac{1}{2}$ and $c = -\frac{3}{2}$, then
\begin{gather*}
K_1^\dagger \simeq 0.007 , \quad \delta_\dagger =
\frac{5}{4} , \quad \zeta \simeq 1.798 , \quad \alpha
\simeq 4.771 , \\
A = \frac{25}{32} , \quad \Delta_1 \simeq 0.599
\quad \text{and} \quad \Delta_2 \simeq 0.052 ,
\end{gather*}
while, if $K = - \frac{1}{2}$ and $c = -4$, then
\begin{gather*}
K_1^\dagger \simeq 0.595 , \quad \hat{x} \simeq
2.542 , \quad K_0^\dagger
\simeq 5 \times 10^{-4} , \quad \delta_\dagger =
\frac{5}{2} , \quad \zeta \simeq 2.440 , \quad \alpha
\simeq 8.336 , \\
A = \frac{25}{8} , \quad \Delta_1 \simeq 0.813 \quad
\text{and} \quad \Delta_2 \simeq 0.028 .
\end{gather*}
} \end{ex}

\noindent
{\bf Case III.2 (Figure~8)}
This case is the modification of Case~II.3 that arises when
it is optimal to abandon the project when this is in its ``closed''
mode and the process $X$ takes sufficiently low values.
In this case,
\begin{gather}
\acal_1 = \mbox{} ]0, \delta] , \quad \pcal = \mbox{} ]\delta,
\gamma[ \mbox{} \cup \mbox{} ]\beta, \infty[ , \quad \scalout
= [\gamma, \beta] , \nonumber \\
\acal_0 = \mbox{} ]0, \zeta] , \quad \wcal = \mbox{} ]\zeta,
\alpha[ \quad \text{and} \quad \scalin = [\alpha, \infty[ ,
\nonumber
\end{gather}
and the required solution to the HJB equation
(\ref{HJBa})--(\ref{HJBb}) should be given by the functions
\begin{align}
w_1 (x) & = \left. \begin{cases}
-K , & \text{if } x \leq \delta \\
\Gamma_1 x^m + \Gamma_2 x^n + R_h (x) , & \text{if }
\delta < x < \gamma \\
\Delta_1 x^m + \Delta_2 x^n - K_0, & \text{if } \gamma
\leq x \leq \beta \\
Ax^m + R_h (x),  & \text{if } x > \beta 
\end{cases} \right\} \label{w0,III.2} \\
\text{and} \quad
w_0 (x) & = \left. \begin{cases}
-K, & \text{if } x \leq \zeta\\
\Delta_1 x^m + \Delta_2 x^n , & \text{if } \zeta < x
< \alpha \\
Ax^m + R_h (x) - K_1, & \text{if } x  \geq \alpha 
\end{cases} \right\} . \label{w1,III.2}
\end{align}

\begin{picture}(160,65)

\put(35,30){\begin{picture}(80,40)
\put(-20,0){\vector(1,0){120}}
\put(-20,15){\vector(1,0){120}}
\put(-30,-0.6){\scriptsize $z=0$}
\put(-30,14.4){\scriptsize $z=1$}
\put(103,-0.6){\scriptsize $x$}
\put(103,14.4){\scriptsize $x$}

\put(73,0){\circle*{1.3}}
\put(72,-4){$\alpha$}
\put(81,-5){\vector(0,1){10}}
\put(91,-5){\vector(0,1){10}}
\put(5,15){\circle*{1.3}}
\put(4,17){$\delta$}
\put(37.5,15){\circle*{1.3}}
\put(36,17){$\gamma$}
\put(0,0){\circle*{1.3}}
\put(-1,-4){$\zeta$}
\put(43.5,20){\vector(0,-1){10}}
\put(52,20){\vector(0,-1){10}}
\put(58,15){\circle*{1.3}}
\put(57,17){$\beta$}

\put(73,-10){\scriptsize $w_0(x)=w_1(x)-K_1$}
\put(80,-15){\scriptsize (switch in)}
\put(20,-10){\scriptsize $w_0(x)=\Delta_1 x^m+\Delta_2 x^n$}
\put(29,-15){\scriptsize (waiting)}
\put(-20.5,22.5){\scriptsize (abandonment)}
\put(-19,27.5){\scriptsize $w_1(x)=-K$}
\put(40,22.5){\scriptsize (switch out)}
\put(35,27.5){\scriptsize $w_1(x) = w_0 (x) - K_0$}
\put(11.5,17.5){\scriptsize (production)}
\put(6,22.5){\scriptsize $\Gamma_1x^m+\Gamma_2x^n+R_h(x)$}
\put(15,27.5){\scriptsize $w_1(x)=$}
\put(70,22.5){\scriptsize (production)}
\put(65,27.5){\scriptsize $w_1(x)=Ax^m+R_h(x)$}
\put(-20.5,-15){\scriptsize (abandonment)}
\put(-19,-10){\scriptsize $w_0(x)=-K$}
\linethickness{0.5mm}
\put(-20,15){\qbezier(0,0)(10,0)(25,0)}
\linethickness{0.5mm}
\put(-20,0){\qbezier(0,0)(10,0)(20,0)}
\end{picture}}

\put(35,8){\small{{\bf Figure 8.} Illustration of the regions
determining}}
\put(36,3){\small{the optimal strategy in the context of Case III.2}}
\end{picture}

\noindent
Once again, we specify $\Gamma_1$, $\Gamma_2$,
$A$, $\Delta_1$, $\Delta_2$, $\zeta$, $\delta$, $\gamma$,
$\beta$ and $\alpha$ by requiring that the functions $w_1$,
$w_0$ should be $C^1$.
This requirement implies that the free-boundary points
$\zeta$, $\delta$, $\gamma$, $\beta$ and $\alpha$ should
satisfy the system of equations given by (\ref{eq II.1.1}),
(\ref{eq II.1.2}), 
\begin{align}
G_3 (\delta, \gamma, \beta) & := n \int _\delta^\infty s^{-n-1}
\bigl[ h(s) + rK \bigr] \, ds - n \int _\gamma^\beta s^{-n-1}
\bigl[ h(s) + rK_0 \bigr] \, ds \nonumber \\
& = 0 , \label{P III.2.20'} \\
G_4 (\zeta, \beta) & := n \int _\beta^\infty s^{-n-1} \bigl[
h(s) + rK_0 \bigr] \, ds + r K \zeta^{-n} \nonumber \\
& = 0 \label{P III.2.30'} \\
\text{and} \quad
G_5 (\zeta, \delta, \gamma) & := m \int _0^\gamma s^{-m-1}
\bigl[ h(s) + rK_0 \bigr] \, ds - m \int _0^\delta s^{-m-1}
\bigl[ h(s) + rK \bigr]  \, ds -rK \zeta ^{-m} \nonumber \\
& = 0 , \label{P III.2.29'}
\end{align}
while the constants $\Gamma_1$, $\Gamma_2$, $A$,
$\Delta_1$, $\Delta_2$ should be given by
\begin{gather}
\Gamma_1 = - \frac{1}{\sigma^2 (n-m)} \int_0^\delta
s^{-m-1} \bigl[ h(s) + rK \bigr] \, ds , \label{III.2.Gam1} \\
\Gamma_2 = - \frac{1}{\sigma^2 (n-m)} \int _\delta^\infty
s^{-n-1} \bigl[ h(s) + rK \bigr] \, ds , \label{III.2.Gam2} \\
\Delta_1 = \frac{rK \zeta^{-m}}{\sigma^2 m(n-m)} ,
\qquad
\Delta_2 = - \frac{rK \zeta^{-n}}{\sigma^2 n(n-m)}
\label{III.2.D12} \\
\text{and} \quad
A = \Delta_1 - \frac{1}{\sigma^2 (n-m)} \int _0^\alpha
s^{-m-1} \bigl[ h(s) - rK_1 \bigr] \, ds . \label{III.2.A}
\end{gather}

\begin{lem} \label{lem III.2}
The system of equations (\ref{eq II.1.1}), (\ref{eq II.1.2}),
(\ref{P III.2.20'}), (\ref{P III.2.30'}), (\ref{P III.2.29'}) has
a unique solution $(\delta, \gamma, \beta, \alpha)$ such
that $0 < \delta < \gamma < \beta < \alpha $ and  the
functions $w_1$, $w_0$ defined by (\ref{w0,III.2}),
(\ref{w1,III.2}), for $\Gamma_1 > 0$, $\Gamma_2 > 0$,
$A > 0$, $\Delta_1 > 0$, $\Delta_2 > 0$ given by
(\ref{III.2.Gam1})--(\ref{III.2.A}), are increasing and
satisfy the HJB equation (\ref{HJBa})--(\ref{HJBb})
if and only if
\be
K < 0 , \quad h(0) < -rK_0 , \quad h(\delta_\dagger)
< 0 , \quad K_1 < K_1^\dagger \quad \text{and} \quad
K_0 < K_0^\dagger ,
\ee
where $\delta_\dagger > 0$ is the unique
solution to (\ref{P II.2.6M}), and $K_1^\dagger > 0$
(resp., $K_0^\dagger > 0$), which depends on all problem
data except $K_1$, $K_0$ (resp., $K_0$)
is as in Lemma~\ref{lem III.1}.
\end{lem}

\begin{ex} \label{EXl} {\rm
If $h$ is the function given by (\ref{h-ex}), then the
system of equations (\ref{P III.1.14})--(\ref{P III.1.15})
takes the form
\begin{align}
(rK + c) \delta^{-n} + \frac{n}{n - \vartheta} (\delta
^{- (n - \vartheta)} - \gamma ^{-(n - \vartheta)} + \beta
^{- (n - \vartheta)}) + (rK_0 + c) (\beta ^{-n} - \gamma ^{-n})
& = 0 , \nonumber \\
(rK_0 + c) \beta^{-n} + \frac{n}{n - \vartheta} \beta
^{-(n - \vartheta)} + rK \zeta^{-n} & = 0 , \nonumber \\
(rK + c) \delta^{-m} - (rK_0 + c) \gamma^{-m} -
\frac{m}{\vartheta - m} (\delta ^{\vartheta - m} - \gamma
^{\vartheta - m}) - rK \zeta^{-m} & = 0 , \nonumber
\end{align}
while
\begin{gather*}
\Gamma_1 = - \frac{\delta^{-m}}{\sigma^2 (n-m)} \left(
- \frac{rK+c}{m} + \frac{\delta^\vartheta}{\vartheta - m}
\right) , \quad
\Gamma_2 = - \frac{\delta^{-n}}{\sigma^2 (n-m)} \left(
\frac{rK+c}{n} + \frac{\delta^\vartheta}{n-\vartheta} \right)
\\
\text{and} \quad
A = \Delta_1 - \frac{\alpha^{-m}}{\sigma^2 (n-m)} \left(
\frac{rK_1 - c}{m} + \frac{\alpha^\vartheta}{\vartheta - m}
\right) .
\end{gather*}
If the problem data is as in Example~\ref{EXf}, then this
case characterises the optimal strategy if and only if
$K < 0$, $c < \min \{ {-1}, rK \}$, $\frac{1}{2} <
K_1^\dagger$ and $\frac{1}{2} < K_0^\dagger$.
If $K = - \frac{1}{50}$ and $c = -13$, then
\begin{gather*}
K_1^\dagger \simeq 392.048 , \quad \hat{x} \simeq
9.756, \quad K_0^\dagger \simeq 0.501 , \\
\zeta \simeq 0.806 , \quad \delta \simeq 6.514 , \quad
\gamma \simeq 7.924 , \quad \beta \simeq 7.942 ,
\quad \alpha \simeq 22.275 , \\
\Gamma_1 \simeq 21.242 , \quad \Gamma_2 \simeq
5 \times 10^{-5} , \quad \Delta_1 \simeq 0.011 , \quad
\Delta_2 \simeq 0.010 \quad \text{and} \quad
A \simeq 21.266 .
\end{gather*}
} \end{ex}

\subsection{The main result}

The following table summarises the conditions on
the problem data that determine the optimality of
each of the cases that we have studied in
Sections~\ref{SS-I}-\ref{SS-III}.
An inspection of the table reveals that these mutually
exclusive conditions exhaust the whole range of
possible problem data.
Therefore, Lemmas~\ref{lem I.1}-\ref{lem III.2} provide
a complete solution to the HJB equation
(\ref{HJBa})--(\ref{HJBb}).

\begin{center}
\begin{tabular}{||c|c|c|c||}
\hline\hline
\multicolumn{2}{||c|}{Conditions on $K_1 > 0$, $K_0 > 0$,
$K \in \bbr$ and $h(\cdot)$} & Case & $w_1$, $w_0$
\Tstrut\Bstrut \\ \hline \hline
\multirow{9}{*}{$0 \leq K$}
& $r K_1 \leq h(0)$
& I.1, Lemma~\ref{lem I.1} & (\ref{w1,I}), (\ref{w0,I1})
\Tstrut\Bstrut \\ \cline{2-4}
& $\max\{- rK_0, \, -rK\} \leq h(0) < rK_1$
& I.2, Lemma~\ref{lem I.2} & (\ref{w1,I}), (\ref{w0,I2})
\Tstrut\Bstrut \\ \cline{2-4}
& $K_0 \leq K$ and $h(0) < -rK_0$
& II.1, Lemma~\ref{lem II.1} & (\ref{w0,w1,II1-1}), (\ref{w0,w1,II1-2})
\Tstrut\Bstrut \\ \cline{2-4}
&$K < K_0$ and $-rK_0 \leq h(0) < -rK$
& II.2, Lemma~\ref{lem II.2} & (\ref{w0,w1,II2-2}), (\ref{w0,w1,II2-1})
\Tstrut\Bstrut \\ \cline{2-4}
&$K < K_0^\star \leq K_0$ and $h(0) < -rK_0$
& II.2, Lemma~\ref{lem II.2} & (\ref{w0,w1,II2-2}), (\ref{w0,w1,II2-1})
\Tstrut\Bstrut \\ \cline{2-4} 
&$K < K_0 < K_0^\star$ and $h(0) < -rK_0$
& II.3, Lemma~\ref{lem II.3} & (\ref{w0,w1,II3-2M}), (\ref{w0,w1,II3-1M})
\Tstrut\Bstrut \\ \hline
\multirow{11}{*}{$K < 0$}
& $rK_1-rK \leq h(0)$
& I.1, Lemma~\ref{lem I.1} & (\ref{w1,I}), (\ref{w0,I1})
\Tstrut\Bstrut \\ \cline{2-4}
& $-rK \leq h(0) < rK_1 - rK$
& I.3, Lemma~\ref{lem I.3} & (\ref{w1,I}), (\ref{w0,I3})
\Tstrut\Bstrut \\ \cline{2-4}
& $-rK_0 \leq h(0) < -rK$
& III.1, Lemma~\ref{lem III.1} & (\ref{w0,w1,III1-1M}),
(\ref{w0,w1,III1-0M})
\Tstrut\Bstrut \\ \cline{2-4}
& $h(0) < -rK_0$ and & &
\Tstrut \\
&$h(\delta_\dagger) \geq 0$ or $\bigl( h(\delta_\dagger)
< 0$ and $K_1 \geq K_1^\dagger
\bigr)$ & III.1, Lemma~\ref{lem III.1} & (\ref{w0,w1,III1-1M}),
(\ref{w0,w1,III1-0M})
\Tstrut \\ 
& or $\bigl( h(\delta_\dagger) < 0 , \ K_1 <
K_1^\dagger$ and $K_0 \geq K_0^\dagger \bigr)$
& &
\Tstrut\Bstrut \\ \cline{2-4}
& $h(0)< -rK_0$,
& \multirow{2}{*}{III.2, Lemma~\ref{lem III.2}} &
\multirow{2}{*}{(\ref{w0,III.2}), (\ref{w1,III.2})}
\Tstrut \\
&$h(\delta_\dagger) < 0$, $ K_1 < K_1^\dagger$
and $K_0 < K_0^\dagger$ & &
\Tstrut\Bstrut \\ \hline\hline
\end{tabular}
\end{center}
\smallskip

\begin{thm} \label{Thm}
Consider the stochastic optimal control problem
formulated in Section~\ref{pr-form} and suppose that
Assumption~\ref{Assumption} holds true.
The value function $v$ is given by (\ref{eqn:v=w}), where
$w_1$, $w_0$ are as in Lemmas~\ref{lem I.1}-\ref{lem III.2}.
In each of the possible cases arising, the optimal strategy
$(Z^\circ, \tau^\circ)$ is as discussed in the proof below.
\end{thm}
{\bf Proof.} 
Given any initial condition $(z,x) \in \{ 0,1 \} \times \mbox{}
]0, \infty[$ and any strategy $(Z,\tau) \in \Pi _z$, the
monotone convergence theorem and (\ref{DOM}) in
Assumption~\ref{Assumption} imply that $\lim
_{m \rightarrow \infty} J_{z,x} (Z, \tau \wedge T_m) =
J_{z,x} (Z, \tau)$  for every sequence of times $(T_m)$
such that $T_m \rightarrow \infty$.
By construction, there exists a constant $C>0$ such that
\be
\bigl| w (z,x) \bigr| \leq C \bigl( 1 + \bigl| R_h (x) \bigr|
\bigr) \quad \text{and} \quad
\bigl| w_x (z,x) \bigr| \leq C \bigl( 1 + \bigl| R_h' (x)
\bigr| \bigr) \quad \text{for all } x>0 ,
\ee
where $w(z,x) = z w_1 (x) + (1-z) w_0 (x)$.
These estimates, (\ref{Rh3}) and (\ref{Rh4}) imply that
\be
\lim _{T \rightarrow \infty} \bbe \Bigl[ e^{-rT} \bigl|
w(Z_T, X_T) \bigr| \Bigr] = 0 ,
\ee
and that the process $M$ defined by
\be
M_T = \int _0^T e^{-rt} X_t w_x (Z_t, X_t) \, dW_t
\ee
is a square integrable martingale for every switching
strategy $Z \in \zcal$.
Furthermore, $w_1$, $w_0$ are $C^1$ as well as $C^2$
outside a finite set, and they satisfy the HJB equation
(\ref{HJBa})--(\ref{HJBb}) in the classical sense.
In view of these observations, we can see that
Theorem~1 in Zervos~\cite{Z03} implies that $w = v$
as long as there exists an optimal strategy
$(Z^\circ, \tau^\circ)$, namely, a switching strategy
$Z^\circ \in \zcal$ such that
\begin{align}
\sigma^2 X_t^2 w_{xx} (Z_t^\circ, & X_t)
+ b X_t w_x (Z_t^\circ, X_t) - r w(Z_t^\circ, X_t)
+ Z_t^\circ h(X_t) = 0 , \nonumber \\
& \bigl[ w(1, X_t) - w(0, X_t) - K_1 \bigr] (\Delta
Z_t^\circ)^+ = 0 \nonumber \\
\text{and} \quad
& \bigl[ w(0, X_t) - w(1, X_t) - K_0 \bigr] (\Delta
Z_t^\circ)^- = 0 , \nonumber
\end{align}
for all $t \leq \tau^\circ$, where
\be
\tau^\circ = \inf \left\{ t \geq 0 \mid \ w(Z_t^\circ , X_t)
= - K \right\} .
\ee
Such a switching strategy is constructed in Duckworth
and Zervos~\cite[Theorem~5]{DZ01} and
Zervos~\cite[Theorem~1]{Z03} for Cases~I.1, I.2, II.1,
II.2 and II.3.
For the remaining cases, it can be constructed using
similar arguments.
\mbox{}\hfill$\Box$

\section{Conclusion}

In this paper, we considered a general entry-exit-scrapping
model with positive switching costs.
We fully characterised the optimal switching and abandonment
strategy by deriving an explicit solution to the control problem's
HJB equation.
It turned out that the optimal strategy can take eight
qualitatively different forms, depending on the problem data.
The analysis of these cases gives rise to the observation
that value may be added by waiting before choosing between
two investment actions of a qualitatively different nature
(one partially reversible and one totally irreversible). 
Furthermore, it suggests that having  ``waiting''
regions to separate regions of the state space associated
with different types of actions should be a generic rather than an
exceptional property of the optimal strategy in real option
models.

\section*{Acknowledgments}
We thank an anonymous referee and an associate editor
for comments and suggestions that enhanced our original
manuscript.
\bigskip

\noindent
The research of Carlos Oliveira was supported by
Funda\c{c}\~ao para a Ci\^encia e Tecnologia through the
grant SFRH/BD/102186/2014.

\section*{Appendix I: auxiliary results}

We first prove the following results that we will use to
streamline the proofs of
Lemmas~\ref{lem I.1}-\ref{lem III.2} in the main paper.

\begin{lem} \label{auxiliar_1}
Suppose that the function $h$ satisfies the requirements
of Assumption~\ref{Assumption}.
Given any constants $\nu \geq 0$ and $L$,
\ben
\lim _{x \rightarrow \infty} \int _\nu^x s^{-m-1} \bigl[
h(s) + L \bigr] \, ds = \infty . \label{eq A1}
\een
Also, given any $\nu \in \mbox{} ]0,\infty]$ and $L$ such
that $h(0)+L<0$,
\ben
\lim _{x \downarrow 0} \int _x^\nu s^{-n-1} \bigl[ h(s) + L
\bigr] \, ds = - \infty . \label{eq A2}
\een
\end{lem}
{\bf Proof.}
Since $h$ is increasing  and $\lim _{x \rightarrow \infty}
h(x) = \infty$, there exist constants $x_1 > \nu$, $M>0$
such that $h(x) + L > M$ for all $x \geq x_1$.
Therefore,
\be
\lim _{x \rightarrow \infty} \int _\nu^x s^{-m-1} \bigl[
h(s) + L \bigr] \, ds \geq \lim _{x \rightarrow \infty} \left[
\int _\nu^{x_1} s^{-m-1} \bigl[ h(s) + L \bigr] \, ds +
\frac{M}{m} x_1^{-m} - \frac{M}{m} x^{-m}\right] = \infty
\ee
because $m<0$.
The assumption $h(0) + L < 0$ implies that there exist
constants $x_2 > 0$ and $\varepsilon > 0$ such that
$h(x) + L \leq -\varepsilon$ for all $x \leq x_2$.
It follows that
\be
\lim _{x \downarrow 0} \int _x^\nu s^{-n-1} \bigl[ h(s) + L
\bigr] \, ds \leq \lim _{x \downarrow 0} \left[ \int _{x_2}^\nu
s^{-n-1} \bigl[ h(s) + L \bigr] \, ds + \frac{\varepsilon}{n}
x_2^{-n} - \frac{\varepsilon}{n} x^{-n} \right] = -\infty
\ee
because $n>0$.
\mbox{}\hfill$\Box$
\bigskip

We have assumed that $h$ is increasing and
right-continuous rather than strictly increasing and
continuous.
Therefore, the statements as well as the proofs of the
following two results have to take into account carefully
the possible jumps or intervals of constancy of $h$.

\begin{lem} \label{auxiliar_2}
Suppose that the function $h$ satisfies the requirements
of Assumption~\ref{Assumption}.
Fix any constants $\nu$ and $L$ such that $h(0) + L < 0$
and
\be
\nu > \inf \bigl\{ x > 0 \mid \ h(x) + L > 0 \bigr\} \geq \inf
\bigl\{ x > 0 \mid \ h(x) + L \geq 0 \bigr\} =: \ubar{\nu} ,
\ee
and consider the function
$q : \mbox{} ]0, \nu] \rightarrow \bbr$ defined by
\be
q(x) = \frac{mx^{m-1}}{\sigma^2 (n-m)} \int _x^\nu s^{-m-1}
\bigl[ h(s) + L \bigr] \, ds - \frac{nx^{n-1}}{\sigma^2 (n-m)}
\int _x^\nu s^{-n-1} \bigl[ h(s) + L \bigr] \, ds .
\ee
There exists a unique $\hat{x} \in [0, \ubar{\nu}[$ such
that
\be
q(x) \left. \begin{cases} > 0 & \text{for all } x \in
\mbox{} ]0, \hat{x}[ , \text{ if } \hat{x} > 0 \\ < 0 &
\text{for all } x \in \mbox{} ]\hat{x}, \nu[ \end{cases}
\right\} .
\ee
\end{lem}
{\bf Proof.}
We define
$p(x) = \sigma^2 x^{-m+1} q(x)$, we note that
\ben
p(\nu) = 0 \quad \text{and} \quad p(x) < 0 \quad
\text{for all } x \in [\ubar{\nu}, \nu[ , \label{aux-proof-I}
\een
and we use the integration by parts formula to calculate
\begin{align}
p'(x) & = x^{n-m-1} \left( x^{-n} \bigl[ h(x) + L \bigr]
- n \int _x^\nu s^{-n-1} \bigl[ h(s) + L \bigr] \, ds \right)
\nonumber \\
& = x^{n-m-1} \left( \nu^{-n} \bigl[ h(\nu) + L \bigr]
- \int _{]x, \nu]} s^{-n} \, dh(s) \right)
=: x^{n-m-1} u(x) , \nonumber
\end{align}
where we have taken into account that $h$ is
right-continuous.
The function $u$ is increasing because $h$ is an
increasing function.
Combining this observation with the assumption that
$h(\nu) + L > 0$, we can see that, if we define
\be
\ubar{x} := \inf \bigl\{ x \in \mbox{} ]0, \nu] \mid \ u(x)
\geq 0 \bigr\} \leq \inf \bigl\{ x \in \mbox{} ]0, \nu] \mid
\ u(x) > 0 \bigr\} =: \bar{x} ,
\ee
then $\ubar{x}, \bar{x} \in [0, \nu]$ and
\be
p'(x) \left. \begin{cases} < 0 & \text{for all } x \in \mbox{}
]0, \ubar{x}[ , \text{ if } \ubar{x} > 0 \\ = 0 & \text{for all }
x \in [\ubar{x}, \bar{x}] , \text{ if } \bar{x} > 0 \\ > 0 &
\text{for all } x \in \mbox{} ]\bar{x}, \nu[ , \text{ if }
\bar{x} < \nu \end{cases} \right\} .
\ee
These inequalities and (\ref{aux-proof-I}) imply that
$\bar{x} < \nu$.
Furthermore, if $\hat{x} = \inf \bigl\{ x \in \mbox{}
]0, \nu] \mid \ p(x) < 0 \bigr\}$, then $\hat{x} > 0$ if and
only if $\lim _{x \downarrow 0} p(x) > 0$, and, if
$\hat{x} > 0$, then $p'(\hat{x}) < 0$.
The required conclusions follow from these observations
and the fact that $p(x)$ and $q(x)$ have the same sign.
\mbox{}\hfill$\Box$

\begin{lem}\label{auxiliar_2'}
Suppose that the function $h$ satisfies the requirements
of Assumption~\ref{Assumption}.
Fix any constants $\nu > 0$ and $L$ such that $h(\nu)
+ L <0$, and consider the function $q: \mbox{} ]\nu, \infty[
\mbox{} \rightarrow \bbr$ defined by
\be
q(x) = \frac{mx^{m-1}}{\sigma^2 (n-m)} \int _\nu^x s^{-m-1}
\bigl[ h(s) + L \bigr] \, ds - \frac{nx^{n-1}}{\sigma^2 (n-m)}
\int _\nu^x s^{-n-1} \bigr[ h(s) + L \bigr] \, ds .
\ee
There exists a unique $\hat{x} \in \mbox{} ]\bar{\nu}, \infty]$
such that
\be
q(x) \left. \begin{cases} > 0 & \text{for all } x \in
\mbox{} ]\nu, \hat{x}[ \\ < 0 & \text{for all } x
> \hat{x}, \text{ if } \hat{x} < \infty \end{cases} \right\} ,
\ee
where
\be
\bar{\nu} := \inf \bigl\{ x > 0 \mid \ h(x) + L > 0 \bigr\}
\geq \inf \bigl\{ x > 0 \mid \ h(x) + L \geq 0 \bigr\} =:
\ubar{\nu} > \nu .
\ee
\end{lem}
{\bf Proof.}
We define $p(x) = \sigma^2 x^{-m+1} q(x)$ and we note that
the right continuity of $h$ implies that $\ubar{\nu} > \nu$.
A simple inspection of the definition of $p$ reveals that
\ben
p(\nu) = 0 \quad \text{and} \quad p(x) > 0 \text{ for all }
x \in \mbox{} ]\nu, \ubar{\nu}] . \label{aux-proof-II}
\een
Using the integration by parts formula, we calculate 
\begin{align}
p'(x) & = - x^{n-m-1} \left( x^{-n} \bigl[ h(x) + L \bigr]
+ n \int _\nu^x s^{-n-1} \bigl[ h(s) + L \bigr] \, ds \right)
\nonumber \\
& = - x^{n-m-1} \left( \nu^{-n} \bigl[ h(\nu) + L \bigr]
+ \int _{]\nu, x]} s^{-n} \, dh(s) \right) =: - x^{n-m-1}
u(x) , \nonumber
\end{align}
where we have taken into account that $h$ is
right-continuous.
The function $u$ is increasing because $h$ is an
increasing function.
In view of this observation and the assumption that
$h(\nu) + L < 0$, we can see that, if we define
\be
\ubar{x} := \sup \bigl\{ x \geq \nu \mid \ u(x)
< 0 \bigr\} \leq \sup \bigl\{ x \geq \nu \mid \ u(x)
\leq 0 \bigr\} = : \bar{x} ,
\ee
then $\ubar{x} , \bar{x} \in [\nu, \infty]$ and
\be
p'(x) \left. \begin{cases} > 0 & \text{for all } x \in \mbox{}
]\nu, \ubar{x}[ , \text{ if } \ubar{x} > \nu \\ = 0 &
\text{for all } x \in [\ubar{x} , \bar{x}] , \text{ if }
\bar{x} > \nu \\ < 0 & \text{for all } x \in \mbox{}
]\bar{x} , \infty[ , \text{ if } \bar{x} < \infty \end{cases}
\right\} .
\ee
These inequalities and (\ref{aux-proof-II}) imply that
$\ubar{x} > \nu$.
If we define $\hat{x} = \sup \bigl\{ x \geq \nu \mid \ p(x)
> 0 \bigr\}$, then $\hat{x} < \infty$ if and only if
$\lim _{x \rightarrow \infty} p(x) < 0$, and, if $\hat{x} < \infty$,
then $p'(\hat{x}) < 0$.
The required conclusions follow from these observations
and the fact that $p(x)$ and $q(x)$ have the same sign.
\mbox{}\hfill$\Box$

\begin{lem} \label{auxiliar_3}
Suppose that the function $h$ satisfies the requirements
of Assumption~\ref{Assumption}.
The function $f: \mbox{} ]0, \infty[ \mbox{} \rightarrow \bbr$
defined by $f(x) = x^{-m+1} R_h' (x)$ is strictly increasing.
\end{lem}
{\bf Proof.}
Using the expression (\ref{Rh'(x)}) for $R_h'$ and the
integration by parts formula, we calculate
\be
f'(x) = \frac{1}{\sigma^2} x^{n-m-1} \left( - x^{-n} h(x)
+n \int _x^\infty s^{-n-1} h(s) \, ds \right)
= \frac{1}{\sigma^2} x^{n-m-1} \int _x^\infty s^{-n} \, dh(s) ,
\ee 
and the claim follows because $h$ is increasing and
$\lim _{x \rightarrow \infty} h(x) = \infty$.
\mbox{}\hfill$\Box$
\bigskip

We will also need the following simple real analysis
result.

\begin{lem} \label{auxiliar_4}
Given points $0 < z_1 <z_2$ and $\kappa \in \bbr$, if
$f: [z_1, z_2] \rightarrow \bbr$ is any right-continuous
increasing function that is not identically 0 and is such
that
\ben
\int _{z_1}^{z_2} s^\kappa f(s) \, ds = 0 , \label{lem-f-ass}
\een
then
\ben
\int _{z_1}^{z_2} s^\mu f(s) \, ds < 0 \quad
\text{for all } \mu < \kappa . \label{lem-f-conc}
\een
\end{lem}
{\bf Proof.}
We define $y = \inf \bigl\{ x \in [z_1, z_2] \mid \
f(x) \geq 0 \bigr\}$, and we note that (\ref{lem-f-ass})
and the fact that $f$ is increasing and not identically 0
imply that $y \in \mbox{} ]z_1, z_2[$, $f(x) < 0$ for
all $x \in [z_1, y[$ and $f(x) \geq 0$ for all $x \in
[y, z_2]$.
In view of these observations, we can see that,
given any $\mu < \kappa$,
\be
0 = \int _{z_1}^{z_2} s^{\kappa - \mu} s^\mu f(s)
\, ds > y^{\kappa - \mu} \int _{z_1}^y s^\mu f(s) \, ds
+ y^{\kappa - \mu} \int _y^{z_2} s^\mu f(s) \, ds ,
\ee
and (\ref{lem-f-conc}) follows.
\mbox{}\hfill$\Box$
\bigskip

\section*{Appendix II: proof of Lemmas~\ref{lem I.1}--\ref{lem III.2}}

In each of the proofs, we mark with bold the first occurrence
of each of the conditions determining the optimality of the
case.
\bigskip

\noindent
{\bf Proof of Lemmas~\ref{lem I.1}, \ref{lem I.2} and~\ref{lem II.1}.}  
The functions $w_1$, $w_0$ defined by (\ref{w1,I}), (\ref{w0,I1})
satisfy the HJB equation
\begin{align} 
\max \left\{ \sigma^2 x^2 w_1'' (x) + bx w_1' (x) - rw_1 (x) +
h(x) , \ w_0 (x) - K_0 - w_1 (x) \right\} & = 0 , \label{HJB-DZ1} \\
\max \left\{ \sigma^2 x^2 w_0'' (x) + bx w_0' (x) - rw_0(x) ,
\ w_1 (x) - K_1 - w_0 (x) \right\} & = 0 \label{HJB-DZ2}
\end{align}
if and only if $\pmb{rK_1 \leq h(0)}$ (see Duckworth and
Zervos~\cite[Lemma 2]{DZ01}).
These functions will satisfy the HJB equation
(\ref{HJBa})--(\ref{HJBb}) if and only if
\ben
- w_1 (x) - K \leq 0 \quad \text{and} \quad -w_0 (x) - K
\leq 0 \quad \text{for all } x > 0 . \label{lem1-2-4-ineq}
\een
In view of (\ref{Rh2}), the first of these inequalities is true
if and only if $-rK \leq h(0)$, while the second one is true
if and only if $\pmb{rK_1 - rK \leq h(0)}$, and
Lemma~\ref{lem I.1} follows because $K_1 > 0$.
\smallskip

Equation (\ref{eq I.2}) has a unique solution $\alpha > 0$ and the
functions $w_1$, $w_0$ defined by (\ref{w1,I}), (\ref{w0,I2})
and (\ref{parameter B,I2}) are increasing and satisfy the HJB
equation (\ref{HJB-DZ1})--(\ref{HJB-DZ2}) if and only if
$\pmb{-rK_0 \leq h(0) < rK_1}$ (see Duckworth and
Zervos~\cite[Lemma 3]{DZ01}).
These functions will satisfy the HJB equation
(\ref{HJBa})--(\ref{HJBb}) if only if they satisfy
(\ref{lem1-2-4-ineq}).
The first of these inequalities holds true if and only if
$\pmb{-rK \leq h(0)}$, while the second one is true if
and only if $\pmb{K \geq 0}$ because $\lim
_{x \downarrow 0} w_0 (x) = 0$ and $w_0$ is increasing,
and Lemma~\ref{lem I.2} follows.
\smallskip

The system of equations (\ref{eq II.1.1})--(\ref{eq II.1.2}) has a
unique solution $(\alpha,\beta)$ such that $0 < \beta < \alpha$
if and only if  $\pmb{h(0) < -rK_0}$, in which case, the functions
$w_1$, $w_0$ defined by (\ref{w0,w1,II1-1})--(\ref{P II.1.2})
are increasing and satisfy the HJB equation
(\ref{HJB-DZ1})--(\ref{HJB-DZ2}) (see Duckworth and
Zervos~\cite[Lemma 3]{DZ01}).
Furthermore, the solution $(\alpha,\beta)$ is such that
\ben
\beta < \inf \bigl\{ x > 0 \mid \ h(x) + rK_0 \geq 0 \bigr\}
\quad \text{and} \quad \alpha > \sup \bigl\{ x > 0 \mid
\ h(x) - rK_1 \leq 0 \bigr\} . \label{DZ-ab-prop}
\een
The functions $w_1$, $w_0$ will satisfy the HJB equation
(\ref{HJBa})--(\ref{HJBb}) if and only if they satisfy
(\ref{lem1-2-4-ineq}).
Both of these inequalities will be true if and only if
$\pmb{K \geq K_0}$ because $w_1$, $w_0$ are increasing
and $\lim _{x \downarrow 0} w_1 (x) = \lim _{x \downarrow 0}
w_0 (x) - K_0 = -K_0$, and Lemma~\ref{lem II.1} follows.
\mbox{}\hfill$\Box$
\bigskip

\noindent
{\bf Proof of Lemma~\ref{lem I.3}.}
In view of the monotonicity of $h$, a simple inspection of
(\ref{eq I.3.1})--(\ref{eq I.3.2}) reveals that this system of
equations has no solution if $K=0$.
On the other hand, the functions $w_1$, $w_0$ defined by
(\ref{w1,I}), (\ref{w0,I3}) can satisfy the HJB equation
(\ref{HJBa})--(\ref{HJBb}) only if $\sigma^2 x^2 w_0'' (x)
+ bx w'_0 (x) - rw_0 (x) = rK \leq 0$ for all $x < \zeta$.
We therefore assume that $\pmb{K < 0}$ in what follows.

To establish conditions under which the system of equations
(\ref{eq I.3.1})--(\ref{eq I.3.2}) has a unique solution
$(\zeta,\alpha)$ such that $0 < \zeta < \alpha$ when
$K<0$, we define
\begin{gather}
\ubar{\alpha} := \inf \bigl\{ x > 0 \mid \ h(x) - rK_1 \geq 0
\bigr\} \leq \inf \bigl\{ x > 0 \mid \ h(x) - rK_1 > 0 \bigr\} =:
\bar{\alpha} \nonumber \\
\text{and} \quad
\ubar{\zeta} := \inf \bigl\{ x > 0 \mid \ h(x) + rK - rK_1
\geq 0 \bigr\} \geq \bar{\alpha} ,
\end{gather}
and we note that $\bar{\alpha} = \ubar{\zeta}$ if and only if
$\ubar{\zeta} = 0$.
If $\ubar{\zeta} > 0$, then the assumption that $h$ is
increasing and (\ref{eq A1}) in Lemma~\ref{auxiliar_1}
imply that there exists a unique $\hat{\zeta} > \ubar{\zeta}$
such that
\be
f_1(\zeta, \zeta) = m \int _0^\zeta s^{-m-1} \bigl[
h(s) + rK - rK_1 \bigr] \, ds \left. \begin{cases} > 0 ,
& \text{if } \zeta \in \mbox{} ]0, \hat{\zeta}[ \\ < 0 ,
& \text{if } \zeta \in \mbox{} ]\hat{\zeta}, \infty[
\end{cases} \right\} .
\ee
On the other hand,
\be
\text{if } \ubar{\zeta} = \bar{\alpha} = 0 \text{ then }
f_1(\zeta, \zeta) < 0 \text{ for all } \zeta > 0 .
\ee
Furthermore, the calculation
\be
\frac{\partial f_1 (\zeta, \alpha)}{\partial \alpha}
= m \alpha^{-m-1} \bigl[ h(\alpha) - rK_1 \bigr]
\ee
implies that (I) $f_1 (\zeta, \cdot)$ is strictly increasing
in $]\zeta, \ubar{\alpha}[$ and constant in $[\ubar{\alpha},
\bar{\alpha}]$, if $\bar{\alpha} > 0$ and $\zeta <
\bar{\alpha}$, and (II) $f_1 (\zeta, \cdot)$ is strictly
decreasing in $]\zeta \vee \bar{\alpha}, \infty[$.
Combining these observations with the fact that
$\lim _{\alpha \rightarrow \infty} f_1(\zeta, \alpha)
= -\infty$, which follows from (\ref{eq A1}) in
Lemma~\ref{auxiliar_1}, we can see that, given $\zeta
> 0$, there exists a unique $\alpha > \zeta$ such that
$f_1 (\zeta, \alpha) = 0$ if and only if
\ben
\ubar{\zeta} > 0 \quad \Leftrightarrow \quad
\pmb{h(0) + rK - rK_1 < 0} \label{lem3-ass-app}
\een
and $\zeta \in \mbox{} ]0, \hat{\zeta}[$.
It follows that, if the inequalities in (\ref{lem3-ass-app})
hold true, then there exists a unique function
$\ell : \mbox{} ]0, \hat{\zeta}[ \mbox{} \rightarrow \mbox{}
]\bar{\alpha}, \infty[$ such that
\be
\lim _{\zeta \uparrow \hat{\zeta}} \ell (\zeta) = \hat{\zeta}
, \quad \zeta < \ell (\zeta) \quad \text{and} \quad
f_1 \bigl( \zeta, \ell (\zeta) \bigr) = 0 \quad \text{for all }
\zeta \in \mbox{} ]0, \hat{\zeta}[ .
\ee
Furthermore, differentiating the identity $f_1\bigl( \zeta,
\ell (\zeta) \bigr) =0$ with respect to $\zeta$, we obtain
\ben
\ell' (\zeta) = - \frac{\zeta^{-m-1} rK}{\ell^{-m-1} (\zeta)
\bigl[ h \bigl( \ell (\zeta) \bigr) - rK_1 \bigr]} > 0 ,
\label{lem3-ell'}
\een
the inequality following because $K<0$ and $\ell (\zeta)
> \bar{\alpha}$.

In the presence of (\ref{lem3-ass-app}), we will show
that the system of equations (\ref{eq I.3.1})--(\ref{eq I.3.2})
has a unique solution $(\zeta,\alpha)$ such that $0 < \zeta
< \alpha$ if we prove that the equation $ f_2 \bigl( \zeta,
\ell (\zeta) \bigr) = 0$ has a unique solution $\zeta \in
\mbox{} ]0, \hat{\zeta}[$.
To this end, we note that
\be
\lim _{\zeta \uparrow \hat{\zeta}} f_2 \bigl( \zeta, \ell (\zeta)
\bigr) = f_2 \bigl( \hat{\zeta}, \hat{\zeta} \bigr) = n \int
_{\hat{\zeta}}^\infty s^{-n-1} \bigl[ h(s) + rK - rK_1\bigr]
\, ds > 0
\ee
and
\be
\lim _{\zeta \downarrow 0} f_2 \bigl( \zeta, \ell (\zeta)
\bigr) < \lim _{\zeta \downarrow 0} \left( n \int
_{\bar{\alpha}}^\infty s^{-n-1} \bigl[ h(s) - rK_1\bigr]
\, ds + rK \zeta^{-n} \right) = - \infty .
\ee
Combining these calculations with
\begin{align}
\frac{d f_2 \bigl( \zeta, \ell (\zeta) \bigr)}{d \zeta} & = -n
\ell^{-n-1} (\zeta) \bigl[ h \bigl( l(\zeta) \bigr) - rK_1 \bigl]
\ell' (\zeta) - nrK \zeta^{-n-1} \nonumber \\
& \stackrel{(\ref{lem3-ell'})}{=} -nrK \zeta^{-m-1}
\bigl[ \ell^{m-n} (\zeta) - \zeta^{m-n} \bigr] > 0 , \nonumber
\end{align}
we can see that  that the equation $f_2 \bigl( \zeta,
\ell (\zeta) \bigr) = 0$ has a unique solution $\zeta \in
\mbox{} ]0, \hat{\zeta}[$, as required.

The $C^1$ functions $w_1$, $w_0$ defined by (\ref{w1,I}),
(\ref{w0,I3}) are increasing because $w'_1 (x) = R'_h (x)
\stackrel{(\ref{Rh1})}{\geq} 0$ for all $x>0$, $w'_0(x) = 0$
for all $x \in \mbox{} ]0,\zeta]$, $w'_0 (x) = R'_h (x) \geq 0$
for all $x \geq \alpha$, and
\begin{align}
w_0'(x) & = m \Delta_1 x^{m-1} + n \Delta_2 x^{n-1}
\nonumber \\
& \stackrel{(\ref{P I.3.5})}{=} - \frac{rK}{\sigma^2 (n-m) x}
\left[ \left( \frac{x}{\zeta} \right) ^n - \left( \frac{x}{\zeta} \right)
^m \right] > 0 \quad \text{for all } x \in \mbox{} ]\zeta, \alpha[ .
\nonumber
\end{align}
To show that these increasing functions provide a solution
to the HJB equation (\ref{HJBa})--(\ref{HJBb}), we still need
to prove that
\begin{align}
\sigma^2 x^2 w_0'' + bx w'_0 (x) - rw_0 (x) \leq 0 & \quad
\text{for all } x > \alpha, \label{lem3-HJBineq1} \\
w_0 (x) - w_1 (x) - K_0 \leq 0 & \quad \text{for all } x > 0 ,
\label{lem3-HJBineq2} \\
w_1 (x) - w_0 (x) - K_1 \leq 0 & \quad \text{for all } x \leq
\alpha , \label{lem3-HJBineq3} \\
\text{and} \quad - w_1 (x) - K \leq 0 & \quad \text{for all }
x \geq 0 . \label{lem3-HJBineq4}
\end{align}
The inequality (\ref{lem3-HJBineq1}) is equivalent to
$h(x) - rK_1 \geq 0$ for all $x > \alpha$, which is true
thanks to the fact that $\alpha > \bar{\alpha}$, where
$\bar{\alpha}$ is defined at the beginning of the proof,
and the assumption that $h$ is increasing.
The inequality (\ref{lem3-HJBineq2}) for $x \geq \alpha$
is equivalent to $K_1 + K_0 \geq 0$, which is true by
assumption.
In view of (\ref{Rh1})--(\ref{Rh2}), the inequality
(\ref{lem3-HJBineq4}) is equivalent to
$\pmb{-rK \leq h(0)}$.
Similarly, the inequality (\ref{lem3-HJBineq2}) for
$x \leq \zeta$ is equivalent to $-rK - rK_0 \leq h(0)$,
which is implied by $-rK \leq h(0)$.
The inequality (\ref{lem3-HJBineq3}) for $x < \zeta$
is equivalent to $w_1 (x) + K - K_1 \leq 0$ and will
follow immediately once we have established
(\ref{lem3-HJBineq3}) for $x \in [\zeta, \alpha]$
because $w_1$ is increasing.

To establish (\ref{lem3-HJBineq2})--(\ref{lem3-HJBineq3})
for $x \in [\zeta, \alpha]$, and complete the proof, we note
that these inequalities are equivalent to
\ben
-K_1 - K_0 \leq g_1 (x) \leq 0 \quad \text{for all } x
\in [\zeta, \alpha] ,  \label{lem3-final-ineq}
\een
where $g_1 (x) = w_0 (x) - w_1 (x) - K_0$.
Using (\ref{w1,I}), (\ref{w0,I3}) and (\ref{P I.3.5}) we
calculate
\begin{align}
g_1 (x) = \mbox{} & \frac{x^m}{\sigma^2 (n-m)} \int
_x^\alpha s^{-m-1} \bigl[ h(s) - rK_1 \bigr] \, ds
\nonumber \\
& - \frac{x^n}{\sigma^2 (n-m)} \int _x^\alpha s^{-n-1}
\bigl[ h(s) - rK_1 \bigr] \, ds - K_1 - K_0 \label{lem3-g1} \\
\text{and} \quad g_1' (x) = \mbox{} & \frac{mx^{m-1}}
{\sigma^2 (n-m)} \int_x^\alpha s^{-m-1} \bigl[ h(s)
- rK_1 \bigr] \, ds \nonumber \\ 
& - \frac{nx^{n-1}} {\sigma^2 (n-m)} \int _x^\alpha
s^{-n-1} \bigl[ h(s) - rK_1 \bigr] \, ds . \label{lem3-g1'}
\end{align}
These expressions, the fact that (\ref{lem3-HJBineq2})
holds true for all $x \leq \zeta$, and the $C^1$ continuity
of $w_1$, $w_0$ at $\zeta$ imply that
\be
g_1 (\zeta) \leq 0 , \quad g_1' (\zeta) = - w_1' (\zeta)
\leq 0 , \quad g_1 (\alpha) = -K_1 - K_0 < 0 \quad
\text{and} \quad g_1' (\alpha) = 0 .
\ee
Recalling that $\alpha > \bar{\alpha}$, where $\bar{\alpha}$
is defined at the beginning of the proof, we combine
the inequalities $g_1' (\zeta) \leq 0$ and $g_1' (\alpha)
= 0$ with Lemma~\ref{auxiliar_2} for $\nu = \alpha$,
$L = -rK_1$ and $q=g_1'$, to see that $g_1'(x) < 0$
for all $x \in \mbox{} ]\zeta, \alpha[$.
It follows that $g_1 (x)$ decreases from $g_1 (\zeta)
\leq 0$ to $g_1 (\alpha) = - K_1 - K_0$ as $x$ increases
from $\zeta$ to $\alpha$, and (\ref{lem3-final-ineq})
follows.
\mbox{}\hfill$\Box$
\bigskip

\noindent
{\bf Proof of Lemma~\ref{lem II.2}.}
Using (\ref{eq A2}) in Lemma~\ref{auxiliar_1} and the
assumptions that $h$ is increasing and $\lim
_{x \rightarrow \infty} h(x) = \infty$, we can see that
equation (\ref{P II.2.6M}) has a unique solution
$\delta_\dagger > 0$ if and only if $\pmb{h(0) + rK < 0}$.
Furthermore, the solution $\delta_\dagger$ is such that
\ben
h(x) + rK < 0 \quad \text{for all } x \leq \delta_\dagger .
\label{lem5-proof-1}
\een 
Before addressing the solvability of (\ref{P II.2.9M}),
we note that the functions $w_1$, $w_0$ defined by
(\ref{w0,w1,II2-2}), (\ref{w0,w1,II2-1}) can satisfy the
HJB equation (\ref{HJBa})--(\ref{HJBb}) only if
$w_0 (x) = Bx^n \geq -K$ for all $x \leq \alpha$.
This inequality cannot be true for $x$ arbitrarily close
to 0 if $-K > 0$.
Therefore, we assume in what follows that
\ben
\pmb{K \geq 0} \quad \Rightarrow \quad K + K_1 > 0 ,
\label{lem5-proof-1.5}
\een
the implication following because $K_1 > 0$.

To show that equation (\ref{P II.2.9M}) has a unique solution
$\alpha > \delta_\dagger$, we define
\ben
\delta_\dagger < \ubar{\alpha} := \inf \bigl\{ x > 0 \mid
\ h(x) - rK_1 \geq 0 \bigr\} \leq \inf \bigl\{ x > 0 \mid \
h(x) - rK_1 > 0 \bigr\} =: \bar{\alpha} .
\label{lem5-bar-alpha}
\een
Here, the first inequality follows because $h$
is right-continuous, and
(\ref{lem5-proof-1})--(\ref{lem5-proof-1.5}) imply
that $h(x) - rK_1 < h(x) + rK \leq h(\delta_\dagger)
+ rK < 0$ for all $x \leq \delta_\dagger$.
In view of the calculation
\be
\frac{\partial f (\delta_\dagger,\alpha)}{\partial \alpha}
= m \alpha^{-m-1} \bigl[ h(\alpha) - rK_1 \bigr] ,
\ee
we can see that $f(\delta_\dagger, \cdot)$ is strictly
increasing in $]\delta_\dagger, \ubar{\alpha}[$ and
strictly decreasing in $]\bar{\alpha}, \infty[$. 
Combining this observation with the calculation
$f(\delta_\dagger, \delta_\dagger) = r (K + K_1)
\delta_\dagger^{-m} > 0$ and the fact that
$\lim _{\alpha \rightarrow \infty} f(\delta_\dagger, \alpha)
= -\infty$, which follows from (\ref{eq A1}) in
Lemma~\ref{auxiliar_1}, we can see that equation
(\ref{P II.2.9M}) has a unique solution $\alpha >
\delta_\dagger$.
Furthermore, this solution satisfies 
\ben
\alpha > \bar{\alpha} := \inf \bigl\{ x > 0 \mid \ h(x)
- rK_1 > 0 \bigr\} . \label{lem5-proof-2}
\een

To streamline the proof, we establish the claims on
the solvability of (\ref{xhat-eqn}) below (see
(\ref{xhat-solv}) and the expression of $g_1'$ in
(\ref{lem3-g1'})).

In view of (\ref{lem5-proof-1}) and (\ref{lem5-proof-2}),
a simple inspection of the expressions (\ref{P II.2.5M})
and (\ref{P II.2.8M}) reveals that $A, B > 0$.
Also, the expression (\ref{Rh'(x)}) for $R_h'$ and the
fact that $\delta_\dagger$ satisfies equation
(\ref{P II.2.6M}) imply that $mA = - \delta_\dagger^{-m+1}
R'_h (\delta_\dagger)$.
Therefore,
\be
R'_h(x) + mA x^{m-1} = x^{m-1} \Bigl[ x^{1-m} R'_h (x)
- \delta_\dagger^{1-m} R'_h (\delta_\dagger) \Bigr]
> 0, \quad \text{for all } x > \delta_\dagger ,
\ee
the inequality following from Lemma~\ref{auxiliar_3}.
Using these results, it is straightforward to verify that
the functions $w_1$, $w_0$ defined by (\ref{w0,w1,II2-2}),
(\ref{w0,w1,II2-1}) are both increasing.

To complete the proof, we need to derive additional conditions
under which the functions $w_1$, $w_0$ are indeed solutions
to the HJB equation (\ref{HJBa})--(\ref{HJBb}).
In view of our analysis thus far, this amounts to establishing
the inequalities
\begin{align}
\sigma^2 x^2 w_1'' (x) + bx w'_1 (x) - rw_1 (x) + h(x)
\leq 0 & \quad \text{for all } x < \delta_\dagger ,
\label{lem5-HJBineq1} \\
\sigma^2 x^2 w_0'' (x) + bx w'_0 (x) - rw_0 (x) \leq 0 &
\quad \text{for all } x > \alpha , \label{lem5-HJBineq2} \\
w_0 (x) - w_1 (x) - K_0 \leq 0 & \quad \text{for all } x>0 ,
\label{lem5-HJBineq3} \\
w_1 (x) - w_0 (x) - K_1 \leq 0 & \quad \text{for all }
x \leq \alpha , \label{lem5-HJBineq4} \\
- w_1 (x) - K \leq 0 & \quad \text{for all } x \geq
\delta_\dagger , \label{lem5-HJBineq5} \\
\text{and} \quad - w_0 (x) - K \leq 0 & \quad \text{for all }
x>0 . \label{lem5-HJBineq6}
\end{align}
The inequality (\ref{lem5-HJBineq1}) is equivalent to
$h(x) + rK \leq 0$ for all $x < \delta_\dagger$, which is
true thanks to (\ref{lem5-proof-1}).
Similarly, (\ref{lem5-HJBineq2}) follows from (\ref{lem5-proof-2}).
The inequality (\ref{lem5-HJBineq3}) for $x \leq \delta_\dagger$
is equivalent to $Bx^n + K - K_0 \leq 0$, which can be true
only if $\pmb{K < K_0}$, and will follow immediately once
we establish (\ref{lem5-HJBineq3}) for $x \in [\delta_\dagger,
\alpha]$ because $x \mapsto Bx^n$ is strictly increasing.
Also, (\ref{lem5-HJBineq3}) for $x \geq \alpha$ is equivalent
to $K_0+K_1>0$, which is true by assumption.
For $x \leq \delta_\dagger$, the inequality (\ref{lem5-HJBineq4})
is equivalent to $Bx^n + K + K_1 \geq 0$, which follows from
(\ref{lem5-proof-1.5}) and the fact that $B>0$.
Furthermore, (\ref{lem5-HJBineq5}) and (\ref{lem5-HJBineq6})
hold true because $w_1$, $w_0$ are increasing, $w_1
(\delta_\dagger) = -K$, $\lim _{x \downarrow 0} w_0 (x)
= 0$ and $K \geq 0$.

The inequalities (\ref{lem5-HJBineq3}) and (\ref{lem5-HJBineq4})
for $x \in [\delta_\dagger, \alpha]$ are equivalent to
\ben
- K_1 - K_0 \leq g_1 (x) \leq 0 \quad \text{for all } x \in
[\delta_\dagger, \alpha] , \label{lem5-final-ineq}
\een
where $g_1 (x) =w_0 (x) - w_1 (x)-K_0$.
Using (\ref{w0,w1,II2-2})--(\ref{P II.2.8M}) and (\ref{P II.2.9M}),
we can verify that $g_1$ and $g_1'$ admit the expressions
given by (\ref{lem3-g1}) and (\ref{lem3-g1'}).
These expressions, the inequality (\ref{lem5-HJBineq4})
for $x \leq \delta_\dagger$, which we have established above,
and the $C^1$ continuity of $w_1$, $w_0$ at $\delta_\dagger$
imply that
\be
g_1 (\delta_\dagger) \geq - K_1 - K_0 , \quad
g_1' (\delta_\dagger) = nB \delta_\dagger^{-n-1} > 0 , \quad
g_1 (\alpha) = - K_1 - K_0 \quad \text{and} \quad
g_1' (\alpha) = 0 .
\ee
In view of these results, (\ref{lem5-proof-2}) and
Lemma~\ref{auxiliar_2} for $\nu = \alpha$, $L = -rK_1$ and
$q=g_1'$, we can see that there exists $\hat{x} \in \mbox{}
]\delta_\dagger, \ubar{\alpha}[$, where $\ubar{\alpha}$ is defined
by (\ref{lem5-bar-alpha}), such that
\ben
g_1' (x) \left. \begin{cases} > 0 & \text{for all } x \in
[\delta_\dagger, \hat{x}[ \\ < 0 & \text{for all } x \in
\mbox{} ]\hat{x}, \alpha[ \end{cases} \right\} .
\label{xhat-solv}
\een
It follows that $g_1$ has a unique maximum in
$[\delta_\dagger, \alpha]$ and (\ref{lem5-final-ineq})
holds true if and only if $g_1(\hat{x}) \leq 0$.
Using the expressions (\ref{lem3-g1}), (\ref{lem3-g1'})
of $g_1$, $g_1'$, equation (\ref{P II.2.9M}), and the identity
$\sigma^2 mn  = -r$, we calculate
\begin{align}
g_1 (\hat{x}) & = - \frac{m \hat{x}^m}{r} \int _{\hat{x}}
^\alpha s^{-m-1} \bigl[ h(s) - rK_1 \bigr] \, ds - K_1 - K_0
\nonumber \\
& = \frac{\hat{x}^m}{r} \left( m \int _{\delta_\dagger}^{\hat{x}}
s^{-m-1} \bigl[ h(s) + rK_0 \bigr] \, ds + r(K - K_0)
\delta_\dagger^{-m} \right) . \label{g1(hatx)}
\end{align}
The second of these expressions and the assumption
$K < K_0$ that we have made above imply that
$g_1(\hat{x}) < 0$ and (\ref{lem5-final-ineq})
holds true if
\be
\pmb{0 \leq h(0) + rK_0} .
\ee
On the other hand, the first expression in (\ref{g1(hatx)})
implies that $g_1(\hat{x}) \leq 0$ and (\ref{lem5-final-ineq})
holds true if
\ben
\pmb{h(0) + rK_0 < 0\quad \text{and} \quad
K_0 \geq -K_1 - \frac{m \hat{x}^m}{r} \int _{\hat{x}}
^\alpha s^{-m-1} \bigl[ h(s) - rK_1 \bigr] \, ds =: K_0^\star}
. \label{K0-star}
\een
A simple inspection of (\ref{P II.2.6M})--(\ref{P II.2.9M})
and (\ref{lem3-g1'}) that determine $\delta_\dagger$,
$\alpha$ and $\hat{x}$ reveals that these points do not
depend on $K_0$.
Therefore, $K_0^\star$ is independent of $K_0$.
In the context of (\ref{K0-star}),
\ben
\pmb{K < K_0^\star < - r^{-1} h(0)} . \label{K0-star-range}
\een
The second inequality here follows immediately from the
fact that the second identity in (\ref{g1(hatx)})
implies that $g_1 (\hat{x}) < 0$ for all $K_0 \geq
- r^{-1} h(0)$.
In view of the linear dependence of $g_1 (\hat{x})$ on
$K_0$, we can see that
\be
K_0^\star > K \quad \Leftrightarrow \quad
\bigl( \text{ if } K_0 = K , \text{ then } g_1 (\hat{x}) > 0
\bigr) .
\ee
Combining this observation with the fact that, if
$K = K_0$, then
\be
g_1 (\hat{x}) \geq g_1 (\delta_\dagger) = B
\delta_\dagger^n > 0 ,
\ee
we obtain the first inequality in (\ref{K0-star-range}).

For future reference, we note that the first expression in
(\ref{g1(hatx)}) implies that
\be
g_1 (\hat{x}) = - \frac{\hat{x}^m}{r} \left( m \int
_{\hat{x}}^\alpha s^{-m-1} h(s) \, ds + rK_1 \alpha^{-m}
+ rK_0 \hat{x}^{-m} \right) .
\ee
Combining this result with the fact that $g_1'(\hat{x}) = 0$
and (\ref{lem3-g1'}), we obtain
\be
g_1 (\hat{x}) = - \frac{\hat{x}^n}{r} \left( n \int
_{\hat{x}}^\alpha s^{-n-1} h(s) \, ds + rK_1 \alpha^{-n}
+ rK_0 \hat{x}^{-n} \right) .
\ee
Comparing these identities with
(\ref{eq II.1.1})--(\ref{eq II.1.2}), we can see that
\ben
K_0 = K_0^\star \quad \Leftrightarrow \quad
g_1 (\hat{x}) = 0  \quad \Leftrightarrow \quad
(\hat{x}, \alpha) \text{ is the solution to 
(\ref{eq II.1.1})--(\ref{eq II.1.2})} . \label{(xhat,a)=(b,a)}
\een
\mbox{}\hfill$\Box$
\bigskip

\noindent
{\bf Proof of Lemma~\ref{lem II.3}.}
In view of Lemma~\ref{lem II.1}, the system of equations
(\ref{eq II.1.1})--(\ref{eq II.1.2}) has a
unique solution $(\alpha,\beta)$ such that $0 < \beta < \alpha$
if and only if  $\pmb{h(0) < -rK_0}$.
To establish conditions under which there exists a unique pair
$(\delta,\gamma)$ satisfying the system of equations
(\ref{P II.3.14})--(\ref{P II.3.15}) and such that $0 < \delta <
\gamma < \beta$, we first note that  (\ref{DZ-ab-prop})
and the assumption that $h$ is increasing imply that
\ben
h(x) + rK_0 < 0 \quad \text{for all } x \leq \beta .
\label{lem6-h+rK0-ineq}
\een
In view of this observation, a simple inspection of (\ref{P II.3.14})
reveals that there are no $0 < \delta < \gamma < \beta$ such that
$F_1 (\delta,\gamma) = 0$ if $K \geq K_0$.
Therefore, we assume that $\pmb{K < K_0}$ in what
follows.
Given any $\gamma \in \mbox{} ]0, \beta]$, the calculations
\begin{gather}
\frac{\partial F_1 (\delta, \gamma)}{\partial \delta} = -m
\delta^{-m-1} \bigl[ h(\delta) + rK_0 + r(K - K_0) \bigr]
< 0 \quad \text{for all } \delta \in \mbox{} ]0, \gamma[ ,
\nonumber \\
\lim _{\delta \downarrow 0} F_1(\delta, \gamma) = m \int
_0^\gamma s^{-m-1} \bigl[ h(s) + rK_0 \bigr] \, ds
> 0 \quad \text{and} \quad F_1 (\gamma, \gamma) = r
(K - K_0) \gamma^{-m} < 0 \nonumber
\end{gather}
imply that there exists a unique $\delta \in \mbox{} ]0,
\gamma[$ such that $F_1(\delta, \gamma) = 0$.
It follows that there exists a unique mapping $\ell : \mbox{}
]0, \beta] \rightarrow \mbox{} ]0, \beta[$ such that
\be
\ell (\gamma) < \gamma \quad \text{and} \quad F_1
\bigl( \ell (\gamma) , \gamma \bigr) = 0 \quad \text{for all }
\gamma \in \mbox{} ]0, \beta] .
\ee
Differentiating the identity here with respect to $\gamma$,
we obtain
\ben
\ell' (\gamma) = \frac{\gamma^{-m-1} \bigl[ h(\gamma)
+ rK_0 \bigr]}{\ell^{-m-1} (\gamma) \bigl[ h \bigl( \ell
(\gamma) \bigr) + rK \bigr]} . \label{lem6-ell'}
\een

In view of these results, we can see that the system of
equations (\ref{P II.3.14})--(\ref{P II.3.15}) has a unique
solution $(\delta,\gamma)$ such that $0 < \delta < \gamma
< \beta$ if and only if the equation $F_2 \bigl( \ell (\gamma)
, \gamma \bigr) = 0$ has a unique solution $\gamma \in \mbox{}
]0, \beta[$.
To derive conditions under which this is indeed the case,
we use (\ref{lem6-h+rK0-ineq}) and (\ref{lem6-ell'}) to calculate 
\be
\frac{d F_2 \bigl( \ell (\gamma) , \gamma)}{d\gamma}
= - n \bigl[ h(\gamma) + rK_0 \bigr] \gamma^{-m-1}
\bigl[ \ell^{m-n} (\gamma) - \gamma^{m-n} \bigr]
> 0 \quad \text{for all } \gamma \in \mbox ]0, \beta[ .
\ee
Combining this result with the identity $\lim
_{\gamma \downarrow 0} F_2 \bigl( \ell (\gamma) , \gamma
\bigr) = - \infty$, which follows from (\ref{lem6-h+rK0-ineq})
and the assumption that $K<K_0$, we can see that the
equation $F_2 \bigl( \ell (\gamma) , \gamma \bigr) = 0$
has a unique solution $\gamma \in \mbox{} ]0, \beta[$ if
and only if 
\ben
F_2 \bigl( \ell (\beta) , \beta \bigr) = \int _{\ell (\beta)}^\infty
s^{-n-1} \bigl[ h(s) + rK \bigr] \, ds > 0 . \label{lem6-NS-cond}
\een

To derive necessary and sufficient conditions under which
this inequality holds true, we fix all other problem data
and we parametrise $\beta$ and $\ell (\beta)$ by
$K_0 \in \mbox{} ]K, - r^{-1} h(0)[$.
Differentiating the identities (\ref{eq II.1.1})--(\ref{eq II.1.2})
with respect to $K_0$, we calculate
\be
\frac{\partial \beta (K_0)}{\partial K_0} =
\frac{\sigma^2 (- m \alpha^{n-m} + n \beta^{n-m}) \beta}
{\bigl[ h(\beta) + rK_0 \bigr] (\alpha^{n-m} - \beta^{n-m})}
< 0 ,
\ee
the inequality following thanks to (\ref{lem6-h+rK0-ineq}).
Also, differentiating the identity
\ben
F_1 \bigl( \ell (\beta) , \beta \bigr) \equiv m \int
_{\ell (\beta)} ^\beta s^{-m-1} \bigl[ h(s) + rK_0 \bigr]
\, ds + r(K - K_0) \ell^{-m} (\beta) = 0
\label{lem6-F1(l(b),b)}
\een
with respect to $K_0$, we can see that
\be
\frac{\partial \ell \bigl( \beta (K_0) ; K_0 \bigr)}{\partial K_0}
= \frac{m \beta^{-m-1} \bigl[ h \bigl( \beta (K_0) \bigr) +
rK_0 \bigr] \frac{\partial \beta (K_0)}{\partial K_0} - r
\beta^{-m} (K_0)}
{m \ell^{-m-1} \bigl( \beta (K_0) ; K_0 \bigr) \Bigl[ h \bigl(
\ell \bigl( \beta (K_0) ; K_0 \bigr) \bigr) + rK \Bigr]} .
\ee
Using these results, we can differentiating the identity
(\ref{lem6-NS-cond}) with respect to $K_0$ to obtain
\begin{align}
& \frac{\partial F_2 \bigl( \ell \bigl( \beta (K_0) ; K_0 \bigr)
, \beta (K_0) ; K_0 \bigr)}{\partial K_0} \nonumber \\
& \qquad = - \ell^{-(n-m)} \bigl( \beta (K_0) ; K_0 \bigr)
\left( \beta^{-m-1} \bigl[ h \bigl( \beta (K_0) \bigr) + rK_0
\bigr] \frac{\partial \beta (K_0)}{\partial K_0} - \frac{r}{m}
\beta^{-m} (K_0) \right) \nonumber \\
& \qquad < 0 . \label{lem6-dF2/dK0}
\end{align}

Comparing equation (\ref{P II.2.6M}) and the
second expression in (\ref{g1(hatx)}) with the expression
(\ref{lem6-NS-cond}) and the identity (\ref{lem6-F1(l(b),b)}),
and taking into account (\ref{(xhat,a)=(b,a)}), 
we can see that
\be
F_2 \bigl( \ell (\beta) , \beta \bigr) = 0 \quad
\Leftrightarrow \quad \bigl( \ell (\beta) = \delta_\dagger
\text{ and } \beta = \hat{x} \bigr) \quad \Leftrightarrow
\quad K_0 = K_0^\star ,
\ee
where $\delta_\dagger$, $\hat{x}$ and $K_0^\star$
are as in the analysis that established (\ref{lem5-final-ineq})
in proof of Lemma~\ref{lem II.2}.
In view of this observation and (\ref{lem6-dF2/dK0}),
we can see that the inequality in (\ref{lem6-NS-cond})
holds true if and only if
$\pmb{K_0 \in \mbox{} ]K, K_0^\star[}$.

To proceed further, we first note that the restriction of
$w_1$ in $[\gamma, \infty[$ as well as the function $w_0$
are increasing thanks to Lemma~\ref{lem II.1}.
We also note that (\ref{lem6-h+rK0-ineq}) implies that
$\Gamma_1 > 0$ and $\Gamma_2 > 0$.
The function $w_1$ is constant in $]0, \delta]$.
Furthermore, it is increasing in $[\delta, \gamma]$ because
\begin{align*}
w'_1 (x) & = m \Gamma_1 x^{m-1} + n \Gamma_2 x^{n-1}
+ R_h' (x) \\
& = \frac{m x^{m-1}}{\sigma^2 (n-m)} \int _\delta^x
s^{-m-1} \bigl[ h(s) + rK \bigr] \, ds - \frac{n x^{n-1}}
{\sigma^2 (n-m)} \int _\delta^x s^{-n-1} \bigl[ h(s) +
rK \bigr] \, ds \nonumber \\
& > 0 \qquad \text{for all } x \in \mbox{} ]\delta, \gamma] ,
\end{align*}
the inequality following thanks to (\ref{lem6-h+rK0-ineq})
and the assumption $K< K_0$ that we have already made.

Since the restriction of $w_1$ in $[\gamma, \infty[$ and
the function $w_0$ satisfy the HJB equation
(\ref{HJB-DZ1})--(\ref{HJB-DZ2}), we will prove that
$w_1$, $w_0$ are indeed solutions to the HJB equation
(\ref{HJBa})--(\ref{HJBb}) if we show that
\begin{align}
\sigma^2 x^2 w_1'' (x) + bx w'_1 (x) - rw_1 (x) + h(x) \leq 0
& \quad \text{for all } x < \delta , \label{lem6-HJB-neq1} \\
w_0 (x) - w_1 (x) - K_0 \leq 0 & \quad \text{for all } x \leq
\gamma , \label{lem6-HJB-neq2} \\
w_1 (x) - w_0 (x) - K_1 \leq 0 & \quad \text{for all } x \leq
\gamma , \label{lem6-HJB-neq3} \\
-w_1 (x) - K \leq 0 & \quad \text{for all } x>\delta ,
\label{lem6-HJB-neq4} \\
\text{and} \quad -w_0 (x) - K \leq 0 & \quad \text{for all } x>0
\label{lem6-HJB-neq5} .
\end{align}
The inequality (\ref{lem6-HJB-neq1}) follows immediately
from (\ref{lem6-h+rK0-ineq}), the fact that $\delta < \beta$
and the assumption that $K<K_0$.
The inequality (\ref{lem6-HJB-neq4}) follows immediately
from the facts that $w_1$ is increasing and $w_1 (\delta)
= -K$, while the inequality (\ref{lem6-HJB-neq5}) is
equivalent to $\pmb{K \geq 0}$ because $w_0$ is also
increasing.
The inequality (\ref{lem6-HJB-neq3}) for $x \leq \delta$
is equivalent to $Bx^n + K + K_1 \geq 0$, which is true
because $B>0$.
For $x < \delta$, (\ref{lem6-HJB-neq2}) holds true if
$w_0 (\delta) - w_1 (\delta) - K_0 \equiv B \delta^n
+ K - K_0 \leq 0$ because $B>0$.
Therefore, (\ref{lem6-HJB-neq2}) for $x < \delta$ will
follow immediately once we establish it for $x \geq \delta$.

The inequalities (\ref{lem6-HJB-neq2}) and
(\ref{lem6-HJB-neq3}) for $x \in [\delta, \gamma]$
are equivalent to
\ben
- K_1 - K_0 \leq g_2 (x) \leq 0 \quad \text{for all }
x \in [\delta, \gamma] , \label{lem6-final-ineq}
\een
where $g_2 (x) = w_0 (x) - w_1 (x) - K_0$.
Using (\ref{P II.1.2}), (\ref{eq II.1.2}) and
(\ref{P II.3.10})--(\ref{P II.3.11}), we calculate
\begin{align}
g_2 (x) = \mbox{} & \frac{x^m}{\sigma^2 (n-m)} \int
_x^\gamma s^{-m-1} \bigl[ h(s) + rK_0 \bigr] \, ds
\nonumber \\
& - \frac{x^n}{\sigma^2 (n-m)} \int _x^\gamma
s^{-n-1} \bigl[ h(s) + rK_0 \bigr] \, ds \label{lem6-g2} \\
\text{and} \quad
g'_2 (x) = \mbox{} & \frac{mx^{m-1}}{\sigma^2 (n-m)}
\int _x^\gamma s^{-m-1} \bigl[ h(s) + rK_0 \bigr] \,ds
\nonumber \\
& - \frac{nx^{n-1}}{\sigma^2 (n-m)} \int _x^\gamma
s^{-n-1} \bigl[ h(s) + rK_0 \bigr] \, ds > 0 \quad
\text{for all } x \in [\delta, \gamma[ , \label{lem6-g2'}
\end{align}
the inequality following thanks to (\ref{lem6-h+rK0-ineq}).
Combining the fact that $g_2 (x)$ is strictly increasing
as $x$ increases from $\delta$ to $\gamma$ with the
inequality $g_2 (\delta) \geq -K_1 - K_0$, which follows
from (\ref{lem6-HJB-neq3}) for $x \leq \delta$ that we
have established above, and the identity
$g_2 (\gamma) = 0$, we obtain (\ref{lem6-final-ineq}).
\mbox{}\hfill$\Box$
\bigskip

\noindent
{\bf Proof of Lemma~\ref{lem III.1}.}
As we have seen at the beginning of the proof of
Lemma~\ref{lem II.2}, equation (\ref{P II.2.6M}) has a
unique solution $\delta_\dagger > 0$ if and only if
$\pmb{h(0) + rK < 0}$, in which case,
\ben
h(x) + rK < 0 \quad \text{for all } x \leq \delta_\dagger ,
\label{lem7-proof-1}
\een
and $A>0$.
We also note that the inequality 
$\sigma^2 x^2 w_0'' (x) + bx w_0' (x) - rw_0 (x) \leq 0$,
which is associated with the HJB equation (\ref{HJBb}),
can be true for $x < \zeta$ if and only if $K \leq 0$.
If $K=0$, then (\ref{P II.2.6M}) and
(\ref{P III.1.14})--(\ref{P III.1.15})  imply that $\alpha$
and $\delta_\dagger$ should satisfy
\be
\int _\alpha^\infty s^{-n-1} \bigl[ h(s) - rK_1 \bigr] \, ds
= 0 \quad \text{and} \quad
\int _{\delta_\dagger}^\alpha s^{-m-1} \bigl[ h(s) - rK_1
\bigr] \, ds = - \frac{rK_1}{m} \delta_\dagger^{-m} > 0 ,
\ee
which is not possible because $h$ is increasing.
We therefore assume that $\pmb{K < 0}$ in what
follows.
In particular, this assumption implies that
$\Delta_1 > 0$ and $\Delta_2 > 0$.

To establish the solvability of
(\ref{P III.1.14})--(\ref{P III.1.15}), we note that the
calculations
\be
\lim _{\zeta \downarrow 0} G_2 (\delta_\dagger, \zeta,
\alpha) = - \infty \quad \text{and} \quad
\frac{\partial G_2 (\delta_\dagger, \zeta, \alpha)}
{\partial \zeta} = - nrK \zeta^{-n-1} > 0
\ee
ensure that, given any $\alpha > \delta_\dagger$
fixed, the equation $G_2 (\delta_\dagger, \zeta, \alpha)
= 0$ has a unique solution $\zeta \in \mbox{} ]0, \alpha[$
if and only if 
\be
H_1 (\alpha) := G_2 (\delta_\dagger, \alpha, \alpha)
= - rK_1 \alpha^{-n} - n \int _{\delta_\dagger}^\alpha
s^{-n-1} \bigl[ h(s) + rK \bigr] \, ds > 0 .
\ee
In view of the calculations
\begin{gather}
H_1 (\delta_\dagger) = - rK_1 \delta_\dagger^{-n} < 0
, \quad \lim _{\alpha \rightarrow \infty} H (\alpha) = 0
\nonumber \\
\text{and} \quad
H_1' (\alpha) = - n \alpha^{-n-1} \bigl[ h(\alpha) + rK - rK_1
\bigr] \left. \begin{cases} > 0 , & \text{if } \alpha \in \mbox{}
]\delta_\dagger , \utilde{\alpha}[ \\ < 0 , & \text{if } \alpha \in
\mbox{} ]\tilde{\alpha} , \infty[ \end{cases} \right\} ,
\nonumber
\end{gather}
where
\begin{align}
\delta_\dagger < \utilde{\alpha} := \inf \bigl\{ x > 0 \mid
\ h(x) + rK - rK_1 \geq 0\} & \nonumber \\
\leq \inf \bigl\{ x>0 \mid \ h(x) + rK - rK_1 > 0 \bigr\}
& =: \tilde{\alpha} , \label{lem7-tilde(a)}
\end{align}
we can see that there exists a unique $\hat{\alpha} \in
\mbox{} ]\delta_\dagger, \utilde{\alpha}[$ such that 
\begin{equation}
H_1 (\alpha) \left. \begin{cases} < 0 , & \text{if } \alpha
\in [\delta_\dagger , \hat{\alpha} [ \\ > 0, & \text{if }
\alpha \in ]\hat{\alpha} , \infty[ \end{cases} \right\} .
\nonumber
\end{equation}
It follows that there exists a function $\ell : \mbox{}
]\hat{\alpha} , \infty[ \mbox{} \rightarrow \mbox ]0, \infty[$
such that
\begin{equation}
G_2 \bigl(\delta_\dagger, \ell (\alpha) , \alpha \bigr) =
0 \text{ and } \ell (\alpha) < \alpha \text{ for all } \alpha
\in \mbox{} ]\hat{\alpha}, \infty[ , \quad \text{and}
\quad \lim _{\alpha \downarrow \hat{\alpha}} \ell
(\alpha) = \hat{\alpha} . \label{lem7-ell(a)}
\end{equation} 
Furthermore,
\begin{equation}
\ell (\alpha) < \delta_\dagger \ \Leftrightarrow \ G_2
(\delta_\dagger, \delta_\dagger , \alpha) > 0 \quad
\text{and} \quad \ell (\alpha) = \delta_\dagger
\ \Leftrightarrow \ G_2 (\delta_\dagger, \delta_\dagger,
\alpha) = 0 . \label{lem7-ell(a)-delta}
\end{equation}
Differentiating the identity $G_2 \bigl(\delta_\dagger,
\ell (\alpha) , \alpha \bigr) = 0$ with respect to $\alpha$,
we obtain
\ben
\ell' (\alpha) = - \frac{1}{rK} \ell^{n+1} (\alpha)
\alpha^{-n-1} \bigl[ h(\alpha) -rK_1 \bigr] .
\label{lem7-ell'(a)}
\een

In view of the analysis thus far, we will show that
there exist unique $0 < \zeta < \alpha$ such that
(\ref{P III.1.14}) and (\ref{P III.1.15}) hold true if we
prove that there exists a unique $\alpha > \hat{\alpha}$
such that
$G_1 \bigl(\delta_\dagger, \ell (\alpha) , \alpha
\bigr) = 0$.
To this end, we note that
\be
\lim _{\alpha \downarrow \hat{\alpha}} G_1 \bigl(
\delta_\dagger, \ell (\alpha) , \alpha \bigr) = G_1 \bigl(
\delta_\dagger, \hat{\alpha}, \hat{\alpha} \bigr) = rK_1
\delta_\dagger^{-m} + m \int _{\delta_\dagger}
^{\hat{\alpha}} s^{-m-1} \bigl[ h(s) + rK - rK_1 \bigr]
\, ds > 0 ,
\ee
the inequality following because $\hat{\alpha} \in
\mbox{} ]\delta_\dagger, \utilde{\alpha}[$, where
$\utilde{\alpha}$ is defined by (\ref{lem7-tilde(a)}).
In view of the inequality $-rK \ell^{-m} (\alpha) <
-rK \alpha^{-m}$ and (\ref{eq A1}) in Lemma~\ref{auxiliar_1},
we can see that
\be
\lim _{\alpha \rightarrow \infty} G_1 \bigl( \delta_\dagger,
\ell (\alpha) , \alpha \bigr) \leq \lim _{\alpha \rightarrow \infty}
\left( rK_1 \delta_\dagger^{-m} + m \int _{\delta_\dagger}
^\alpha s^{-m-1} \bigl[ h(s) + rK - rK_1 \bigr] \, ds \right)
= - \infty .
\ee
Combining these results with the observation that
\begin{align*}
\frac{\partial G_1 \bigl( \delta_\dagger, \ell (\alpha),
\alpha \bigr)}{\partial \alpha} & = m \alpha^{-m-1} \bigl[
h(\alpha) - rK_1 \bigr] + mrK \ell^{-m-1} (\alpha) \ell'
(\alpha) \\
& = m \alpha^{-n-1} \bigl[ h(\alpha) - rK_1 \bigr] \bigl[
\alpha^{n-m} - \ell^{n-m} (\alpha) \bigr] \nonumber \\
& 
\left. \begin{cases} > 0, & \text{if } \hat{\alpha} <
\ubar{\alpha} \text{ and } \alpha \in \mbox{}
]\hat{\alpha} , \ubar{\alpha}[ \\ < 0, & \text{if } \alpha \in
\mbox{} ]\hat{\alpha} \vee \bar{\alpha}, \infty[ \end{cases}
\right\} ,
\end{align*}
where we have used (\ref{lem7-ell'(a)}) and the
definitions
\be
\ubar{\alpha} := \inf \bigl\{ x>0 \mid \ h(x) - rK_1 \geq 0
\bigr\} \leq \inf \bigl\{ x>0 \mid \ h(x) - rK_1 > 0 \bigr\}
=: \bar{\alpha} < \tilde{\alpha},
\ee
we can see that equation $G_1 \bigl(\delta_\dagger,
\ell (\alpha) , \alpha \bigr) = 0$ has a unique solution
$\alpha > \hat{\alpha}$ such that
\ben
h(x) - rK_1 > 0 \quad \text{for all } x \geq \alpha .
\label{lem7-h(a)-rK1}
\een

We now investigate under what conditions $\zeta \leq
\delta_\dagger$ or $\zeta > \delta_\dagger$.
To this end, we calculate
\ben
G_2 (\delta_\dagger, \delta_\dagger, \alpha) = - n \int
_{\delta_\dagger}^\alpha s^{-n-1} h(s) \, ds - rK_1
\alpha^{-n} . \label{lem7-G2(d,d,a)}
\een
In view of this expression and (\ref{lem7-ell(a)-delta}),
we can see that
\be
\pmb{h(\delta_\dagger) \geq 0 \quad \Rightarrow
\quad \delta_\dagger < \zeta} .
\ee
If $h(\delta_\dagger) < 0$, then we fix all other
problem data and we parametrise $G_1$, $G_2$,
$\zeta$ and $\alpha$ by $K_1 > 0$ (note that
$\delta_\dagger$ does not depend on $K_1$).
Differentiating (\ref{P III.1.14}) and (\ref{P III.1.15}),
and eliminating $\frac{\partial \zeta (K_1)}{\partial K_1}$,
we calculate
\ben
\frac{\partial \alpha (K_1)}{\partial K_1} =
\frac{\sigma^2 \alpha (K_1) \bigl[ n \zeta^{m-n} (K_1)
-m \alpha^{m-n} (K_1) \bigr]} {\bigl[ \zeta^{m-n} (K_1)
- \alpha^{m-n} (K_1) \bigr] \bigl[ h \bigl( \alpha (K_1)
\bigr) - rK_1 \bigr]} > 0 , \label{lem7-z-d-rel1}
\een
the inequality following thanks to (\ref{lem7-h(a)-rK1})
and the fact that $\zeta < \alpha$.
Furthermore, (\ref{lem7-h(a)-rK1}) implies that
\ben
\lim _{K_1 \rightarrow \infty} \alpha (K_1) = \infty .
\label{lem7-z-d-rel2}
\een
Using (\ref{lem7-z-d-rel1}), we calculate
\be
\frac{\partial G_2 \bigl( \delta_\dagger, \delta_\dagger,
\alpha (K_1) ; K_1 \bigr)}{\partial K_1} = - n \alpha^{-n-1}
(K_1) \bigl[ h \bigl( \alpha (K_1) \bigr) - rK_1 \bigl]
\frac{\partial \alpha (K_1)}{\partial K_1} - r
\alpha (K_1) < 0 .
\ee
In view of (\ref{lem7-z-d-rel2}), we can see that
\be
\lim _{K_1 \rightarrow \infty} G_2 \bigl( \delta_\dagger,
\delta_\dagger, \alpha (K_1) ; K_1 \bigr) \leq -n \int
_{\delta_\dagger}^\infty s^{-n-1} h(s) \, ds = rK \delta
_\dagger^{-n} < 0 ,
\ee
where we have also used (\ref{P II.2.6M}) and the
assumption $K<0$ that we have made above.
In light of (\ref{lem7-ell(a)-delta}), it follows that, if
\ben
\lim _{K_1 \downarrow 0} G_2 \bigl( \delta_\dagger,
\delta_\dagger, \alpha (K_1) ; K_1 \bigr) > 0 ,
\label{lem7-z-d-rel3}
\een
then
\begin{align}
& \pmb{h(\delta_\dagger) < 0} \nonumber \\
& \quad \pmb{\Rightarrow \
\text{there exists a unique } K_1^\dagger > 0
\text{ such that } \left. \begin{cases} \zeta < \delta_\dagger
& \text{for all } K_1 \in \mbox{} ]0, K_1^\dagger[ \\
\delta_\dagger < \zeta & \text{for all } K_1 \in \mbox{}
]K_1^\dagger, \infty[ \end{cases} \right\} } .
\label{K1-dagger}
\end{align}
This analysis also establishes
(\ref{lem7-d>z})--(\ref{lem7-d<z}) if (\ref{lem7-z-d-rel3})
holds true whenever $h(\delta_\dagger) < 0$.
For future reference, we note that
\ben
\text{if } h(\delta_\dagger) < 0 , \text{ then } G_2 \bigl(
\delta_\dagger, \delta_\dagger, \alpha (K_1^\dagger) ;
K_1^\dagger \bigr) = 0 . \label{lem7-G2(d,d,a,K1+)}
\een
Furthermore, $K_1^\dagger$ does not depend on
$K_1$ itself or $K_0$.

To prove that (\ref{lem7-z-d-rel3}) is indeed true,
we first note that the analysis of the solvability of
(\ref{P III.1.14})--(\ref{P III.1.15}) remains true
for $K_1 = 0$.
In particular, if we define $\zeta_0 = \lim
_{K_1 \downarrow 0} \zeta (K_1)$ and $\alpha_0
= \lim _{K_1 \downarrow 0} \alpha (K_1)$,
then
\be
\delta_\dagger < \alpha_0 , \quad & \zeta_0 < \alpha_0 ,
\quad \bar{\alpha}_0 := \inf \bigl\{ x>0 \mid \ h(x)
> 0 \bigr\} < \alpha_0 ,
\ee
\begin{align}
L_1 (\zeta_0, \alpha_0) & := \lim _{K_1 \downarrow 0}
G_1 \bigl( \delta_\dagger, \zeta (K_1) , \alpha (K_1) ;
K_1 \bigr) \nonumber \\
& \equiv m \int _{\delta_\dagger}^{\alpha_0} s^{-m-1}
h(s) \, ds + rK \delta_\dagger^{-m} - rK \zeta_0^{-m}
= 0 \nonumber \\
\text{and} \quad
L_2 (\zeta_0, \alpha_0) & := \lim _{K_1 \downarrow 0}
G_2 \bigl( \delta_\dagger, \zeta (K_1), \alpha (K_1) ;
K_1 \bigr) \nonumber \\
& \equiv -n \int _{\delta_\dagger}^{\alpha_0} s^{-n-1}
h(s) \, ds - rK \delta_\dagger^{-n} + rK \zeta_0^{-n} = 0 .
\nonumber
\end{align}
In view of Lemma~\ref{auxiliar_4}, the identities
here cannot be satisfied for $\zeta_ 0 = \delta_\dagger$.
To prove that $\zeta _0 < \delta_\dagger$, which is
equivalent to (\ref{lem7-z-d-rel3}), we argue by
contradiction and we assume that $\delta_\dagger
< \zeta_0$.
The calculations
\be
\frac{\partial L_1 (\zeta, \alpha)}{\partial \alpha} =
m \alpha^{-m-1} h(\alpha) \left. \begin{cases} \geq 0,
& \text{if } \alpha \leq \bar{\alpha}_0 \\ < 0, & \text{if }
\alpha > \bar{\alpha}_0 \end{cases} \right\}
\quad \text{and} \quad
\lim _{\alpha \rightarrow \infty} L_1 (\zeta, \alpha)
= - \infty ,
\ee
where the last one follows from (\ref{eq A1}) in
Lemma~\ref{auxiliar_1}, imply that, given any
$\zeta > \delta_\dagger$, there exists a unique
$\alpha > \zeta$ such that $L_1 (\zeta, \alpha) = 0$
if and only if
\ben
L_1 (\zeta, \zeta) = m \int _{\delta_\dagger}^\zeta
s^{-m-1} \bigl[ h(s) + rK \bigr] \, ds > 0 .
\label{lem7-lambda0}
\een
In view of our assumptions on $h$, (\ref{lem7-proof-1})
and the fact that $\lim _{\zeta \rightarrow \infty}
L_1 (\zeta, \zeta) = -\infty$, which follows
from (\ref{eq A1}) in Lemma~\ref{auxiliar_1}, we can
see that there exists a unique
\be
\zeta^\ddagger > \inf \bigl\{ x> 0 \mid \ h(x)
+rK > 0 \bigr\} > \delta_\dagger
\ee
such that (\ref{lem7-lambda0}) holds true if and only
if $\zeta \in \mbox{} ]\delta_\dagger, \zeta^\ddagger[$.
It follows that there exists a function $\lambda :
\mbox{} ]\delta_\dagger, \zeta^\ddagger[ \mbox{}
\rightarrow \mbox{} ]\bar{\alpha}_0, \infty[$ such that
\ben
L_1 \bigl( \zeta, \lambda (\zeta) \bigr) = 0
\text{ and } \zeta < \lambda (\zeta) \text{ for all }
\zeta \in \mbox{} ]\delta_\dagger, \zeta^\ddagger[ , \quad
\text{and} \quad \lambda (\zeta_0) = \alpha_0 .
\label{lem7-lambda-defn}
\een
Differentiating the identity $L_1 \bigl( \zeta, \lambda
(\zeta) \bigr) = 0$ with respect to $\zeta$, we
calculate
\be
\lambda' (\zeta) = - \frac{rK \zeta^{-m-1}}{\lambda
^{-m-1} (\zeta) h \bigl( \lambda (\zeta) \bigr)} .
\ee
Using this result, we obtain
\be
\frac{d L_2 \bigl( \zeta, \lambda (\zeta) \bigr)}
{d\zeta} = - nrK \zeta^{-m-1} \bigl[ \zeta ^{-(n-m)}
- \lambda ^{-(n-m)} (\zeta) \bigr] > 0 .
\ee
Combining this calculation with (\ref{lem7-lambda-defn})
and the limit
\be
\lim _{\zeta \downarrow \delta_\dagger} L_2 \bigl(
\zeta, \lambda (\zeta) \bigr) = - n \int _{\delta_\dagger}
^{\lambda (\delta_\dagger)} s^{-n-1} h(s) \, ds \geq 0 ,
\ee
where the inequality follows from the fact that
$L_1 \bigl( \delta_\dagger, \lambda (\delta_\dagger)
\bigr) = 0$ and Lemma~\ref{auxiliar_4}, we can see
that there exist no $\zeta_0 < \alpha_0$ such that
$\zeta_0 > \delta_\dagger$ and $L_1 (\zeta_0,
\alpha_0) = L_2 (\zeta_0, \alpha_0) = 0$, which
establishes the required contradiction.

To streamline the proof, we establish the claims on
the solvability of (\ref{xhat-eqn2}) below (see
(\ref{xhat-solv2}) and the expression of $g_1'$ in
(\ref{lem3-g1'})).

To show that $w_1$, $w_0$ are increasing,
it suffices to prove that $w_1$ is increasing in
$[\delta_\dagger, \infty[$ and $w_0$ is increasing in
$[\zeta, \alpha]$.
The first of these claims follows from the calculation
\be
w'_1 (x) = R'_h (x) + mA x^{m-1} = x^{m-1} \bigl[ x^{-m+1}
R_h' (x) - \delta_\dagger^{-m+1} R_h' (\delta_\dagger)
\bigr] > 0 \quad \text{for all } x > \delta_\dagger ,
\ee
where we have used the expression for $A$ given by
(\ref{P II.2.5M}), the identity (\ref{P II.2.6M}), and
Lemma~\ref{auxiliar_3}.
To establish the second claim, we use (\ref{P III.1.12})
and (\ref{P III.1.13}) to calculate
\be
w'_0 (x) = - \frac{rK}{\sigma^2 (n-m) x} \left[ \left(
\frac{x}{\zeta} \right)^n - \left( \frac{x}{\zeta} \right)^m
\right] > 0 \quad \text{for all } x \in \mbox{} ]\zeta,
\alpha[ .
\ee

To show that the functions $w_1$, $w_0$ defined by
(\ref{w0,w1,III1-1M}), (\ref{w0,w1,III1-0M}) satisfy the inequalities
associated with the HJB equation (\ref{HJBa})--(\ref{HJBb}), we
need to show that 
\begin{align}
\sigma^2 x^2 w_1'' (x) + bx w'_1 (x) - rw_1 (x) + h(x)
& \leq 0 \quad \text{for all } x < \delta_\dagger ,
\label{lem7-HJBineq1} \\
\sigma^2 x^2 w_0'' (x) + bx w'_0 (x) - rw_0 (x) & \leq 0
\quad \text{for all } x \in \mbox{} ]0, \zeta[ \mbox{} \cup
\mbox{} ]\alpha, \infty[ , \label{lem7-HJBineq2} \\
w_0 (x) - w_1 (x) - K_0 & \leq 0 \quad \text{for all } x > 0 ,
\label{lem7-HJBineq3} \\
w_1 (x) - w_0 (x) - K_1 & \leq 0 \quad \text{for all }
x \leq \alpha , \label{lem7-HJBineq4} \\
- w_1 (x) - K & \leq 0 \quad \text{for all } x > \delta_\dagger
\label{lem7-HJBineq5} \\
\text{and} \quad
- w_0 (x) - K & \leq 0 \quad \text{for all } x > \zeta .
\label{lem7-HJBineq6}
\end{align}
The inequality (\ref{lem7-HJBineq1}) is equivalent to
$h(x) \leq -rK$ for all $x\leq \delta_\dagger$, which is
true thanks to (\ref{lem7-proof-1}).
The  inequality (\ref{lem7-HJBineq2}) is trivial for $x < \zeta$
and follows immediately from (\ref{lem7-h(a)-rK1}) for
$x > \alpha$.
The inequalities (\ref{lem7-HJBineq3}) and
(\ref{lem7-HJBineq4}) for $x \leq \delta_\dagger
\wedge \zeta$ are equivalent to $K_0 \geq 0$ and
$K_1 \geq 0$, respectively, which are true by assumption,
while (\ref{lem7-HJBineq3}) for $x \geq \alpha$ is also
implied by the assumption that $K_1, K_0 > 0$.
The inequalities (\ref{lem7-HJBineq5}) and
(\ref{lem7-HJBineq6}) follow from the fact that $w_1$,
$w_0$ are increasing and the identities $w_1
(\delta_\dagger) = w_0 (\zeta) = - K$.
If $\delta_\dagger < \zeta$, then (\ref{lem7-HJBineq3})
for $x \in [\delta_\dagger, \zeta]$ is equivalent to
$w_1 (x) \geq - K - K_0$, which is true, while
(\ref{lem7-HJBineq4}) for $x \in [\delta_\dagger, \zeta]$
will follow as soon as we establish it for $x \in [\zeta,
\alpha]$ below because $w_1$ is increasing.

The inequalities (\ref{lem7-HJBineq3}) and
(\ref{lem7-HJBineq4}) for $x \in [\zeta, \alpha]$ are equivalent
to
\ben
- K_1 - K_0 \leq g_1 (x) \leq 0 \quad \text{for all } x \in [\zeta ,
\alpha] , \label{lem7-HJBineq-big}
\een
where $g_1 (x) = w_0 (x) - w_1 (x) - K_0$.
Using (\ref{P II.2.5M}) and (\ref{P III.1.12})--(\ref{P III.1.13}),
we can see that $g_1$ and $g_1'$ admit the expressions
given by (\ref{lem3-g1}) and (\ref{lem3-g1'}).
Furthermore,
\ben
g_1 (\alpha) = - K_1 - K_0 \quad \text{and} \quad
g_1' (\alpha) = 0 . \label{lem7-g1-alpha}
\een

{\em Proof of (\ref{lem7-HJBineq-big}) when $\delta_\dagger
\leq \zeta$\/} ({\em i.e., when $\pmb{h(\delta_\dagger) \geq 0}$
or $\pmb{h(\delta_\dagger) < 0}$ {\bf and} $\pmb{K_1 \geq
K_1^\dagger}$\/}).
In view of the expression for $A$ given by (\ref{P II.2.5M}),
the identity (\ref{P II.2.6M}), and Lemma~\ref{auxiliar_3},
we can see that
\begin{align}
g_1' (\zeta) & = - R'_h (\zeta) - m A \zeta^{m-1}
\nonumber \\
& = - \zeta^{m-1} \bigl[ \zeta^{-m+1} R'_h (\zeta) -
\delta_\dagger^{-m+1} R'_h (\delta) \bigr] \left.
\begin{cases} = 0 , & \text{if } \zeta = \delta_\dagger \\
< 0 , & \text{if } \delta_\dagger < \zeta \end{cases} \right\}
, \label{lem7-g1'-xhat-delta}
\end{align}
Also, since $w_1$ is increasing,
\be
g_1 (\zeta) \leq - K - w_1 (\delta_\dagger) - K_0 =
- K_0 < 0 .
\ee
Combining these inequalities with (\ref{lem7-g1-alpha})
and Lemma~\ref{auxiliar_2}, we can see that $g_1 (x)$ is
decreasing from $g_1 (\zeta) < 0$ to $- K_1 - K_0$ as
$x$ increases from $\zeta$ to $\alpha$, and
(\ref{lem7-HJBineq-big}) follows.

{\em Proof of (\ref{lem7-HJBineq-big}) when $\zeta <
\delta_\dagger$\/} ({\em i.e., when $\pmb{h(\delta_\dagger)
< 0}$ {\bf and} $\pmb{K_1 < K_1^\dagger}$\/}).
Since $w_0$ is strictly increasing in $]\zeta, \alpha[$,
\be
g_1 (\zeta) = - K_0 > - K_1 - K_0 \quad \text{and} \quad
g_1' (\delta_\dagger) = w'_0 (\delta_\dagger) > 0 .
\ee
Combining these inequalities with (\ref{lem7-g1-alpha})
and Lemma~\ref{auxiliar_2}, we can see that there exists
a unique
\ben
\hat{x} \in \mbox{} ]\delta_\dagger , \alpha[ ,
\label{lem7-xhat}
\een
such that
\ben
g_1' (x) \left. \begin{cases} > 0 & \text{for all } x \in [\zeta,
\hat{x}[ \\ < 0 & \text{for all } x \in \mbox{} ]\hat{x}, \alpha[
\end{cases} \right\} . \label{xhat-solv2}
\een
In particular, $g_1$ has a unique maximum in $[\zeta,
\alpha]$ and (\ref{lem7-HJBineq-big}) holds true if and
only if $g_1 (\hat{x}) \leq 0$.
Using the expressions (\ref{lem3-g1}), (\ref{lem3-g1'})
of $g_1$, $g_1'$ and the identity $\sigma^2 mn  = -r$,
we calculate
\begin{align}
g_1 (\hat{x}) & = - \frac{n \hat{x}^n}{r} \int _{\hat{x}}
^\alpha s^{-n-1} \bigl[ h(s) - rK_1 \bigr] \, ds - K_1 -
K_0 \nonumber \\
& = - \frac{m \hat{x}^m}{r} \int _{\hat{x}}^\alpha s^{-m-1}
\bigl[ h(s) - rK_1 \bigr] \, ds - K_1 - K_0 .
\label{lem7-g1(hatx)}
\end{align}
The second of these expressions and (\ref{P III.1.14})
imply that
\begin{align}
g_1 (\hat{x}) & = \frac{\hat{x}^m}{r} \left( m \int
_{\delta_\dagger}^{\hat{x}} s^{-m-1} \bigl[ h(s) +
rK_0 \bigr] \, ds + r(K - K_0) \delta_\dagger^{-m}
- rK \zeta^{-m} \right) \nonumber \\
& < \frac{\hat{x}^m}{r} \left( m \int _{\delta_\dagger}
^{\hat{x}} s^{-m-1} \bigl[ h(s) + rK_0 \bigr] \, ds
- rK_0 \delta_\dagger^{-m} \right) . \nonumber
\end{align}
This inequality reveals that $g_1(\hat{x}) < 0$ and
(\ref{lem7-HJBineq-big}) holds true if
\be
\pmb{0 \leq h(0) + rK_0} .
\ee
On the other hand, the first expression in
(\ref{lem7-g1(hatx)}) implies that $g_1(\hat{x})
\leq 0$ and (\ref{lem7-HJBineq-big}) holds true if
\ben
\pmb{h(0) + rK_0 < 0 \quad \text{and} \quad
K_0 \geq -K_1 - \frac{n \hat{x}^n}{r} \int _{\hat{x}}
^\alpha s^{-n-1} \bigl[ h(s) - rK_1 \bigr] \, ds =: K_0^\dagger
> 0} . \label{K0-dagger}
\een
A simple inspection of (\ref{P II.2.6M}),
(\ref{P III.1.14})--(\ref{P III.1.15}) and (\ref{lem3-g1'})
that determine $\delta_\dagger$, $\alpha$ and $\hat{x}$
reveals that these points do not depend on $K_0$.
Therefore, $K_0^\dagger$ is independent of $K_0$.
To see the last inequality in (\ref{K0-dagger}), we note
that $g_1 (\cdot) +K_0$ does not depend on $K_0$ and
has a unique global maximum in $[\delta_\dagger,
\alpha]$ at $\hat{x} \in \mbox{} ]\delta_\dagger,
\alpha[$.
Therefore,
\ben
K_0^\dagger = g_1 (\hat{x}) + K_0 > g_1 (\delta_\dagger)
+ K_0 = \frac{\delta_\dagger^n}{r} G_2 (\delta_\dagger,
\delta_\dagger, \alpha) > 0 , \label{lem7-lem8-K0=0-con}
\een
the last inequality following from (\ref{lem7-ell(a)-delta})
and the fact that $\zeta = \ell (\alpha) < \delta_\dagger$
in this part of the analysis.

For future reference, we note that (\ref{lem7-g1(hatx)})
implies that,
\begin{align}
g_1 (\hat{x}) & = - \frac{\hat{x}^n}{r} \left( n \int
_{\hat{x}}^\alpha s^{-n-1} h(s) \, ds + rK_1 \alpha^{-n}
+ rK_0 \hat{x}^{-n} \right) \nonumber \\
& = - \frac{\hat{x}^m}{r} \left( m \int _{\hat{x}}^\alpha
s^{-m-1} h(s) \, ds + rK_1 \alpha^{-m} + rK_0
\hat{x}^{-m} \right) . \nonumber
\end{align}
Comparing these identities with
(\ref{eq II.1.1})--(\ref{eq II.1.2}), we can see that,
if $\pmb{h(\delta_\dagger) < 0}$, $\pmb{K_1
\in \mbox{} ]0, K_1^\dagger]}$ and $\pmb{h(0) + rK_0
< 0}$, then
\ben
K_0 = K_0^\dagger \quad \Leftrightarrow \quad
g_1 (\hat{x}) = 0  \quad \Leftrightarrow \quad
(\hat{x}, \alpha) \text{ is the solution to 
(\ref{eq II.1.1})--(\ref{eq II.1.2})} .
\label{lem7-(xhat,a)=(b,a)}
\een
Furthermore, if we fix all other problem data
and we parametrise $\hat{x}$, $\zeta$ and $K_0^\dagger$
by $K_1$, then (\ref{lem7-G2(d,d,a,K1+)}),
(\ref{lem7-g1'-xhat-delta}) and a calculation similar to
the one in (\ref{lem7-lem8-K0=0-con}) imply that
\ben
\hat{x} (K_1^\dagger) = \zeta (K_1^\dagger) = \delta_\dagger
\text{ and } \lim _{K_1 \uparrow K_1^\dagger}
K_0^\dagger (K_1) = 0 , \label{lem7-K0dagger=0}
\een
(\ref{P II.2.6M}), (\ref{P III.1.15}) and
(\ref{lem7-(xhat,a)=(b,a)}) imply that
\ben
n \int _{\hat{x} (K_1)}^\infty s^{-n-1} \bigl[ h(s) +
rK_0^\dagger (K_1) \bigr] \, ds + r K \zeta^{-n}
(K_1) = 0 , \label{lem7-K0dagger=0-1}
\een
while (\ref{P III.1.14}) and (\ref{lem7-(xhat,a)=(b,a)})
imply that
\begin{align}
m \int _0^{\hat{x} (K_1)} s^{-m-1} \bigl[ h(s) +
rK_0^\dagger (K_1) \bigr] \, ds & \nonumber \\
- m \int _0^{\delta_\dagger} s^{-m-1} \bigl[
h(s) + rK \bigr] \, ds & - rK \zeta ^{-m} (K_1) = 0 .
\label{lem7-K0dagger=0-2}
\end{align}
\mbox{}\hfill$\Box$
\vspace{5mm}

\noindent
{\bf Proof of Lemma~\ref{lem III.2}.}
In view of Lemma~\ref{lem II.1}, the system of equations
(\ref{eq II.1.1})--(\ref{eq II.1.2}) has a
unique solution $(\alpha,\beta)$ such that $0 < \beta < \alpha$
if and only if  $\pmb{h(0) + rK_0 < 0}$, in which case,
\ben
h(x) + rK_0 < 0 \text{ for all } x \leq \beta \quad \text{and}
\quad h(x) - rK_1 > 0 \text{ for all } x \geq \alpha .
\label{lem8-h+rK0-K1-ineqs}
\een
Equations (\ref{eq II.1.1})--(\ref{eq II.1.2}) imply that
(\ref{P III.2.30'}) is equivalent to
\ben
n \int _\alpha^\infty s^{-n-1} \bigl[ h(s) - rK_1 \bigr] \, ds
+ r K \zeta^{-n} = 0 . \label{Lem8-proof.G4}
\een
In view of the second of the inequalities in
(\ref{lem8-h+rK0-K1-ineqs}), we can see that there is
no $\zeta>0$ such that (\ref{Lem8-proof.G4})
holds true unless $\pmb{K < 0}$.
On the other hand, if $K<0$, then it is straightforward to
see that there exists a unique $\zeta \in \mbox{} ]0,
\alpha[$ such that (\ref{Lem8-proof.G4}) holds true.

For future reference, we note that the first of the
inequalities in (\ref{lem8-h+rK0-K1-ineqs}) and the
assumption $K<0$ imply that  $h(0) + rK < 0$.
Therefore, there exists a unique $\delta_\dagger > 0$
such that
\ben
\int _x^\infty s^{-n-1} \bigl[ h(s) + rK \bigr] \, ds \left.
\begin{cases} < 0 , & \text{for all } x \in \mbox{}
]0, \delta_\dagger[ \\ > 0 , & \text{for all } x \in \mbox{}
]\delta_\dagger, \infty[ \end{cases} \right\} .
\label{lem8-delta-dagger}
\een

To establish the required solvability of (\ref{P III.2.20'})
and (\ref{P III.2.29'}), we first fix any $\gamma \in \mbox{}
]0, \beta[$ and we look for $\delta \in \mbox{} ]0, \gamma[$
such that $G_3 (\delta, \gamma, \beta) = 0$.
Combining the limit
\be
\lim _{\delta \downarrow 0} G_3 (\delta, \gamma, \beta)
= - \infty ,
\ee
which follows from (\ref{eq A2}) in Lemma~\ref{auxiliar_1},
with the calculation
\be
\frac{\partial G_3 (\delta, \gamma, \beta)}{\partial \delta}
= - n \delta^{-n-1} \bigl[ h(\delta) + rK_0 - r(K_0 -K) \bigr]
> 0 \quad \text{for all } \delta \in \mbox{} ]0,
\gamma[ ,
\ee
where the inequality follows thanks to
(\ref{lem8-h+rK0-K1-ineqs}) and the assumption
that $K < 0 < K_0$,
we can see that there exists $\delta \in \mbox{} ]0,
\gamma[$ such that $G_3 (\delta, \gamma, \beta)
= 0$ if and only if 
\ben
G_3 (\gamma, \gamma, \beta)
= n \int _\beta^\infty s^{-n-1} \bigl[ h(s) + rK \bigr]
\, ds - r(K_0-K) \gamma^{-n} + r(K_0-K) \beta^{-n}
> 0 . \label{lem8-G3(ggb)}
\een
In view of the calculations
\be
\lim _{\gamma \downarrow 0} G_3 (\gamma,
\gamma, \beta) = - \infty \quad \text{and} \quad
\frac{\partial G_3 (\gamma, \gamma, \beta)}
{\partial \gamma} = n r(K_0 - K) \gamma^{-n-1}
> 0 ,
\ee
there exists a unique $\hat{\gamma} \in \mbox{}
]0, \beta[$ such that (\ref{lem8-G3(ggb)}) holds
true for all $\gamma \in \mbox{} ]\hat{\gamma},
\beta]$ if and only if
\ben
H_2 (\beta) := G_3 (\beta, \beta, \beta) = n \int
_\beta^\infty s^{-n-1} \bigl[ h(s) + rK \bigr] \, ds
> 0 \quad \Leftrightarrow \quad \delta_\dagger
< \beta , \label{lem8-G3(b)}
\een
where $\delta_\dagger$ is as in
(\ref{lem8-delta-dagger}).
It follows that, if the problem data is such that
(\ref{lem8-G3(b)}) holds true, then $G_3 (\delta,
\gamma, \beta) = 0$ defines uniquely a mapping
$\ell : \mbox{} ]\hat{\gamma}, \beta] \rightarrow
\mbox{} ]0, \beta[$, such that
\begin{gather}
G_3 \bigl( \ell (\gamma), \gamma, \beta \bigr) = 0
\quad \text{and} \quad \ell (\gamma) < \gamma
\quad \text{for all } \gamma \in \mbox{} ]\hat{\gamma},
\beta[ , \nonumber \\
\ell (\hat{\gamma}) :=
\lim _{\gamma \downarrow \hat{\gamma}} \ell
(\gamma) = \hat{\gamma} \quad \text{and} \quad
\ell (\beta) = \delta_\dagger .
\nonumber
\end{gather}
Differentiating the identity $G_3 \bigl( \ell (\gamma),
\gamma, \beta \bigr) = 0$ with respect to $\gamma$,
we derive the expression 
\ben
\ell' (\gamma) = \frac{\gamma^{-n-1} \bigl[ h(\gamma)
+ rK_0 \bigr]}{\ell^{-n-1} (\gamma) \bigl[ h \bigl(
\ell (\gamma) \bigr) + rK \bigr]} > 0 . \label{lem8-ell'}
\een
Furthermore, comparing the identity
\be
\lim _{\gamma \downarrow \hat{\gamma}} G_3
(\gamma, \gamma, \beta) = n \int _\beta^\infty
s^{-n-1} \bigl[ h(s) + rK_0 \bigr]
\, ds - r(K_0-K) \hat{\gamma}^{-n} = 0
\ee
with equation (\ref{P III.2.30'}) that $\zeta > 0$ satisfies,
we can see that $K (\zeta ^{-n} - \hat{\gamma}^{-n})
= - K_0 \hat{\gamma}^{-n}$.
It follows that $\zeta < \hat{\gamma} $ because
$K < 0 < K_0$.
We conclude this part of the analysis with the
observation that
\ben
\zeta < \ell (\gamma) < \delta_\dagger \quad
\text{for all } \gamma \in \mbox{} ]\hat{\gamma},
\beta[ , \label{lem8-zeta-ghat}
\een
where we have taken into account that
$\zeta < \hat{\gamma} = \ell (\hat{\gamma})$
and the fact that $\ell$ is strictly increasing.

To determine conditions under which there exists a
unique $\gamma \in \mbox{} ]\hat{\gamma}, \beta[$
such that $G_5 \bigl( \zeta, \ell (\gamma), \gamma
\bigr) =0$ if (\ref{lem8-G3(b)}) holds true, we first
note that
\be
G_5 \bigl( \zeta, \ell (\hat{\gamma}), \hat{\gamma} \bigr)
= G_5 (\zeta, \hat{\gamma}, \hat{\gamma}) = - rK_0
\hat{\gamma}^{-m} + rK (\hat{\gamma}^{-m} - \zeta^{-m})
< 0 ,
\ee
the inequality following thanks to (\ref{lem8-zeta-ghat})
and the fact that $K < 0 < K_0$.
Combining this observation with the calculation
\be
\frac{\partial G_5 \bigl( \zeta, \ell (\gamma), \gamma)}
{\partial\gamma} & = m \gamma^{-n-1} \bigl[ h(\gamma)
+ rK_0 \bigr] \bigl( \gamma^{n-m} - \ell^{n-m} (\gamma)
\bigr) > 0 ,
\ee
where we have used (\ref{lem8-ell'}), we can see that
there exists a unique $\gamma \in \mbox{} ]\hat{\gamma},
\beta[$ such that $G_5 \bigl( \zeta, \ell (\gamma),
\gamma \bigr) =0$ if and only if 
\ben
G_5 \bigl( \zeta, \ell (\beta), \beta \bigr) = G_5
(\zeta, \delta_\dagger, \beta) > 0 .
\label{lem8-G5(b)}
\een

To derive conditions under which (\ref{lem8-G3(b)})
and (\ref{lem8-G5(b)}) hold true, we first note
that (\ref{lem8-h+rK0-K1-ineqs}) implies that
$h(\beta) < 0$.
Therefore,  (\ref{lem8-G3(b)}) can be true only if
$\pmb{h(\delta_\dagger) < 0}$.
If we fix all other problem data and we parametrise
$\alpha$, $\beta$, $\zeta$ and $\delta$ by $K_1$
and $K_0$, then (\ref{lem7-(xhat,a)=(b,a)}) and
(\ref{lem7-K0dagger=0}) imply that
\ben
\lim _{K_0 \downarrow 0} H_2 \bigl( \beta (K_1^\dagger,
K_0) \bigr) = H_2 \bigl( \hat{x} (K_1^\dagger) \bigr)
= n \int _{\delta_\dagger}^\infty s^{-n-1} \bigl[ h(s) +
rK \bigr] \, ds  = 0 . \label{lem8-G35-conds1}
\een
Furthermore, differentiating
(\ref{eq II.1.1})--(\ref{eq II.1.2}) and using the resulting
expressions, we obtain
\begin{gather}
\frac{\partial H_2 \bigl( \beta (K_1, K_0) \bigr)}
{\partial K_1} = - \frac{\sigma^2 (n-m) \bigl[
h (\beta) + rK \bigr] \alpha ^{-m}}{\bigl[ h (\beta)
+ rK_0 \bigr] (\alpha ^{n-m} - \beta^{n-m})}
< 0 \label{lem8-G35-conds2}
\intertext{and}
\frac{\partial H_2 \bigl( \beta (K_1, K_0) \bigr)}
{\partial K_0} = - \frac{\sigma^2 \bigl[ h (\beta)
+ rK \bigr] (- m \alpha^{n-m} + n \beta^{n-m})
\beta^{-n}} {\bigl[ h(\beta) + rK_0 \bigr] (\alpha^{n-m}
- \beta^{n-m})} < 0 . \label{lem8-G35-conds3}
\end{gather}
These calculations imply that
$H_2 \bigl( \beta (K_1, K_0) \bigr) < 0$ for
all $K_1 \geq K_1^\dagger$ and $K_0 > 0$.
On the other hand, (\ref{lem7-xhat}) and
(\ref{lem7-(xhat,a)=(b,a)}) imply that
$\beta (K_1, K_0^\dagger) = \hat{x} (K_1) >
\delta_\dagger$ for all $K_1 < K_1^\dagger$,
which, combined with the equivalences in
(\ref{lem7-(xhat,a)=(b,a)}), (\ref{lem8-G3(b)})
and the inequality (\ref{lem8-G35-conds3}),
implies that $H_2 \bigl( \beta (K_1, K_0) \bigr)
> 0$ for all $K_1 < K_1^\dagger$ and $K_0
\in \mbox{} ]0, K_0^\dagger]$.

In view of the results derived above, we will
conclude this part of the analysis if we show
that, given any $\pmb{K_1 < K_1^\dagger}$,
(\ref{lem8-G5(b)}) holds true if and only if
$\pmb{K_0 \in \mbox{} ]0, K_0^\dagger[}$.
In the context of the conditions $h(\delta_\dagger)
< 0$ and $K_1 < K_1^\dagger$, we can see that
a straightforward comparison of (\ref{P III.2.30'}),
which defines $\zeta$ (see the analysis in the first
paragraph of this proof), and (\ref{lem7-K0dagger=0-1})
reveals that the free-boundary point $\zeta =
\zeta (K_1, K_0)$ in this lemma identifies with
the free-boundary point $\zeta = \zeta (K_1)$
in Lemma~\ref{lem III.1} if $K_0 = K_0^\dagger$.
This observation and a comparison of
(\ref{P III.2.29'}), (\ref{lem7-K0dagger=0-2})
reveal that
\be
G_5 \bigl( \zeta (K_1, K_0^\dagger) ,
\delta_\dagger, \beta (K_1, K_0^\dagger) ;
K_0^\dagger \bigr) = 0 .
\ee
Combining this result with the calculation
\be
\frac{\partial G_5 \bigl( \zeta (K_1, K_0) ,
\delta_\dagger, \beta (K_1, K_0) ;
K_0 \bigr)}{\partial K_0} = \sigma^2 m (n-m)
\beta^{-m} \frac{\alpha^{n-m} - \zeta^{n-m}}
{\alpha^{n-m} - \beta^{n-m}} < 0 ,
\ee
we can see that (\ref{lem8-G5(b)}) holds true
if and only if $K_0 \in \mbox{} ]0, K_0^\dagger[$.

To show that $w_1$, $w_0$ are increasing, it suffices to
prove that $w_0$ is increasing in $[\zeta, \alpha]$ and
$w_1$ is increasing in $[\delta, \gamma] \cup [\beta,
\infty[$.
The first of these claims follows immediately from the calculation
\begin{align}
w'_0 (x) = - \frac{rK}{\sigma^2 (n-m) x} \left[ \left(
\frac{x}{\zeta} \right) ^n - \left( \frac{x}{\zeta} \right) ^m
\right] > 0 \quad \text{for all } x \in \mbox{} ]\zeta,
\alpha] , \nonumber
\end{align}
where we have used (\ref{III.2.D12}) and the assumption
$K<0$ that we have made above in this proof.
Using (\ref{III.2.Gam1})--(\ref{III.2.Gam2}), we calculate
\be
w'_1 (x) = \frac{m x^{m-1}}{\sigma^2 (n-m)} \int _\delta^x
s^{-m-1} \bigl[ h(s) + rK \bigr] \, ds - \frac{n x^{n-1}}
{\sigma^2 (n-m)} \int _\delta^x s^{-n-1} \bigl[ h(s) + rK
\bigr] \, ds ,
\ee
for $x \in [\delta, \gamma]$.
This expression, the fact that $w'_1(\gamma) = w_0'
(\gamma) > 0$ and Lemma~\ref{auxiliar_2'} for
$\nu = \delta$, $L = rK$ and $q = w_1'$ (see also
(\ref{lem8-delta-dagger}), (\ref{lem8-zeta-ghat})
and recall that $\delta = \ell (\gamma)$) imply that
$w_1' (x) > 0$ for all $x \in \mbox{} ]\delta, \gamma]$.
To prove that $w_1$ is increasing in $[\beta, \infty[$,
we first note that the inequality $w'_1(\beta) = w_0'
(\beta) > 0$ implies that $mA > - \beta^{-m+1} R'_h
(\beta)$.
In view of this observation, we can see that
\be
w_1' (x) = R'_h (x) + mA x^{m-1} > x^{m-1} \bigl[
x^{-m+1} R'_h (x) - \beta^{-m+1} R'_h (\beta) \bigr]
> 0 \quad \text{for all } x>\beta ,
\ee
the second inequality following by Lemma~\ref{auxiliar_3}.

To show that $w_1$ and $w_0$ satisfy the HJB equation
(\ref{HJBa})--(\ref{HJBb}), we need to prove that 
\begin{align}
\sigma^2 x^2 w_1'' (x) + bx w'_1 (x) - rw_1 (x) + h(x) & \leq 0
\quad \text{for all } x \in \mbox{} ]0, \delta[ \mbox{} \cup \mbox{}
]\gamma, \beta[ , \label{lem8-HJBineq1} \\
\sigma^2 x^2 w_0'' (x) +bx w'_0 (x)-rw_0 (x) & \leq 0 \quad
\text{for all } x \in \mbox{} ]0, \zeta[ \mbox{} \cup \mbox{}
]\alpha, \infty[ , \label{lem8-HJBineq2} \\
w_0 (x) - w_1 (x) - K_0 & \leq 0 \quad \text{for all } x \in
\mbox{} ]0, \gamma] \cup [\beta, \infty[ ,
\label{lem8-HJBineq3} \\
w_1 (x) - w_0 (x) - K_1 & \leq 0 \quad \text{for all } x \leq
\alpha , \label{lem8-HJBineq4} \\
-w_1 (x) - K & \leq 0 \quad \text{for all } x \geq \delta
\label{lem8-HJBineq5} \\
\text{and} \quad
-w_0 (x) - K & \leq 0 \quad \text{for all } x \geq \zeta .
\label{lem8-HJBineq6}
\end{align}
Inequality (\ref{lem8-HJBineq1}) for $x < \delta$ follows
immediately from (\ref{lem8-delta-dagger}),
(\ref{lem8-zeta-ghat}) and the fact that $\delta
= \ell (\gamma)$.
Inequality (\ref{lem8-HJBineq1}) for $x \in \mbox{} ]\gamma,
\beta[$ and (\ref{lem8-HJBineq2}) for $x > \alpha$ hold true
thanks to (\ref{lem8-h+rK0-K1-ineqs}), while
(\ref{lem8-HJBineq2}) for $x < \zeta$ is equivalent
to $K \leq 0$, which is true by assumption.
The inequalities (\ref{lem8-HJBineq3}) for $x \leq \zeta$
and (\ref{lem8-HJBineq4}) for $x \leq \delta$ are true
because $w_0$ is increasing and $K_1, K_0 >0$.
Also, (\ref{lem8-HJBineq3}) for $x \geq \alpha$
and (\ref{lem8-HJBineq4}) for $x \in [\gamma, \beta]$
are both equivalent to $K_1 + K_0 \geq 0$, while
(\ref{lem8-HJBineq3}) for $x \in \mbox{} ]\zeta, \delta[$
will follow as soon as we establish it for
$x \in [\delta, \gamma]$ below.
Furthermore, (\ref{lem8-HJBineq5}) and
(\ref{lem8-HJBineq6}) follow immediately from the
fact that $w_1$ and $w_0$ are increasing.

To establish (\ref{lem8-HJBineq3}) and
(\ref{lem8-HJBineq4}) for $x \in [\beta, \alpha]$, we
need to show that
\begin{equation}
-K_1 - K_0 \leq g_1 (x) \leq 0 \quad \text{for all } x \in
[\beta, \alpha] , \label{lem8-III.2.53-1}
\end{equation}
where $g_1 (x) = w_0 (x) - w_1 (x) - K_0$.
Using (\ref{eq II.1.2}), (\ref{P III.2.30'}) and
(\ref{III.2.D12})--(\ref{III.2.A}), we can verify
that $g_1$ and $g_1'$ admit the expressions given by
(\ref{lem3-g1}) and (\ref{lem3-g1'}).
These expressions, the fact that (\ref{lem8-HJBineq4})
holds with equality for $x = \beta$, and the $C^1$
continuity of $w_1$, $w_0$ at $\beta$ imply that
\be
g_1 (\beta) = g'_1 (\beta) = 0 , \quad g_1 (\alpha) =
-K_1 - K_0 \quad \text{and} \quad g'_1 (\alpha) = 0 .
\ee
In view of (\ref{lem8-h+rK0-K1-ineqs}) and
Lemma~\ref{auxiliar_2} for $\nu = \alpha$, $L = -rK_1$
and $q = g_1'$, we can see that $g_1' (x) < 0$
for all $x \in \mbox{} ]\beta, \alpha[$.
It follows that $g_1 (x)$ decreases from 0 to
$-K_1 - K_0 < 0$, and (\ref{lem8-III.2.53-1})
holds true.

Finally, the inequalities (\ref{lem8-HJBineq3}) and
(\ref{lem8-HJBineq4}) for $x \in [\delta, \gamma]$
are equivalent to
\ben
-K_1 - K_0 \leq g_2 (x) \leq 0 \quad \text{for all }
x \in [\delta, \gamma] , \label{lem8-III.2.53-2}
\een
where $g_2 (x) = w_0 (x) - w_1 (x) - K_0$.
Using (\ref{P III.2.20'})--(\ref{III.2.D12}) and
 (\ref{lem8-h+rK0-K1-ineqs}), we can verify
that $g_2$, $g_2'$ admit the expressions given by
(\ref{lem6-g2}), (\ref{lem6-g2'}), and $g_2' (x) > 0$
for all $x < \gamma$.
Combining the fact that $g_2$ is strictly increasing
in $]\delta, \gamma[$ with the identity $g_2 (\gamma)
= 0$ and the inequality $g_2 (\delta) \geq -K_1 - K_0$,
which follows from (\ref{lem8-HJBineq4}) for
$x \leq \delta$, we obtain (\ref{lem8-III.2.53-2}).
\mbox{}\hfill$\Box$

\end{document}